\documentclass{amsart}

\usepackage{amssymb,amsmath,stmaryrd,mathrsfs}
\def\definetac{\newif\iftac}    % Can't define a \newif inside another \if!
\ifx\tactrue\undefined
  \definetac
  \ifx\state\undefined\tacfalse\else\tactrue\fi
\fi
\iftac\else\usepackage{amsthm}\fi
\usepackage[all,2cell]{xy}
\UseAllTwocells
\usepackage{enumitem}
\usepackage{xcolor}
\definecolor{darkgreen}{rgb}{0,0.45,0} 
\usepackage[pagebackref,colorlinks,citecolor=darkgreen,linkcolor=darkgreen]{hyperref}
\usepackage{mathtools}          % for all sorts of things
\usepackage{tikz}
\usepackage{braket}             % for \Set, etc.

\usepackage{url}                % for citations to web sites
\usepackage{xspace}             % put spaces after a \command in text

\makeatletter
\let\ea\expandafter

\def\mdef#1#2{\ea\ea\ea\gdef\ea\ea\noexpand#1\ea{\ea\ensuremath\ea{#2}\xspace}}
\def\alwaysmath#1{\ea\ea\ea\global\ea\ea\ea\let\ea\ea\csname your@#1\endcsname\csname #1\endcsname
  \ea\def\csname #1\endcsname{\ensuremath{\csname your@#1\endcsname}\xspace}}

\DeclareRobustCommand\widecheck[1]{{\mathpalette\@widecheck{#1}}}
\def\@widecheck#1#2{%
    \setbox\z@\hbox{\m@th$#1#2$}%
    \setbox\tw@\hbox{\m@th$#1%
       \widehat{%
          \vrule\@width\z@\@height\ht\z@
          \vrule\@height\z@\@width\wd\z@}$}%
    \dp\tw@-\ht\z@
    \@tempdima\ht\z@ \advance\@tempdima2\ht\tw@ \divide\@tempdima\thr@@
    \setbox\tw@\hbox{%
       \raise\@tempdima\hbox{\scalebox{1}[-1]{\lower\@tempdima\box
\tw@}}}%
    {\ooalign{\box\tw@ \cr \box\z@}}}

\newcount\foreachcount

\def\foreachletter#1#2#3{\foreachcount=#1
  \ea\loop\ea\ea\ea#3\@alph\foreachcount
  \advance\foreachcount by 1
  \ifnum\foreachcount<#2\repeat}

\def\foreachLetter#1#2#3{\foreachcount=#1
  \ea\loop\ea\ea\ea#3\@Alph\foreachcount
  \advance\foreachcount by 1
  \ifnum\foreachcount<#2\repeat}

\def\definescr#1{\ea\gdef\csname s#1\endcsname{\ensuremath{\mathscr{#1}}\xspace}}
\foreachLetter{1}{27}{\definescr}
\def\definecal#1{\ea\gdef\csname c#1\endcsname{\ensuremath{\mathcal{#1}}\xspace}}
\foreachLetter{1}{27}{\definecal}
\def\definebold#1{\ea\gdef\csname b#1\endcsname{\ensuremath{\mathbf{#1}}\xspace}}
\foreachLetter{1}{27}{\definebold}
\def\definebb#1{\ea\gdef\csname l#1\endcsname{\ensuremath{\mathbb{#1}}\xspace}}
\foreachLetter{1}{27}{\definebb}
\def\definefrak#1{\ea\gdef\csname f#1\endcsname{\ensuremath{\mathfrak{#1}}\xspace}}
\foreachletter{1}{9}{\definefrak} % \fi is something else!
\foreachletter{10}{27}{\definefrak}
\def\definebar#1{\ea\gdef\csname #1bar\endcsname{\ensuremath{\overline{#1}}\xspace}}
\foreachLetter{1}{27}{\definebar}
\foreachletter{1}{8}{\definebar} % \hbar is something else!
\foreachletter{9}{15}{\definebar} % \obar is something else!
\foreachletter{16}{27}{\definebar}
\def\definetil#1{\ea\gdef\csname #1til\endcsname{\ensuremath{\widetilde{#1}}\xspace}}
\foreachLetter{1}{27}{\definetil}
\foreachletter{1}{27}{\definetil}
\def\definehat#1{\ea\gdef\csname #1hat\endcsname{\ensuremath{\widehat{#1}}\xspace}}
\foreachLetter{1}{27}{\definehat}
\foreachletter{1}{27}{\definehat}
\def\definechk#1{\ea\gdef\csname #1chk\endcsname{\ensuremath{\widecheck{#1}}\xspace}}
\foreachLetter{1}{27}{\definechk}
\foreachletter{1}{27}{\definechk}
\def\defineul#1{\ea\gdef\csname u#1\endcsname{\ensuremath{\underline{#1}}\xspace}}
\foreachLetter{1}{27}{\defineul}
\foreachletter{1}{27}{\defineul}

\def\autofmt@n#1\autofmt@end{\mathrm{#1}}
\def\autofmt@b#1\autofmt@end{\mathbf{#1}}
\def\autofmt@l#1#2\autofmt@end{\mathbb{#1}\mathsf{#2}}
\def\autofmt@c#1#2\autofmt@end{\mathcal{#1}\mathit{#2}}
\def\autofmt@s#1#2\autofmt@end{\mathscr{#1}\mathit{#2}}
\def\autofmt@f#1\autofmt@end{\mathsf{#1}}
\def\autofmt@u#1\autofmt@end{\underline{\smash{\mathsf{#1}}}}
\def\autofmt@U#1\autofmt@end{\underline{\underline{\smash{\mathsf{#1}}}}}
\def\autofmt@h#1\autofmt@end{\widehat{#1}}
\def\autofmt@r#1\autofmt@end{\overline{#1}}
\def\autofmt@t#1\autofmt@end{\widetilde{#1}}
\def\autofmt@k#1\autofmt@end{\check{#1}}

\def\auto@drop#1{}
\def\autodef#1{\ea\ea\ea\@autodef\ea\ea\ea#1\ea\auto@drop\string#1\autodef@end}
\def\@autodef#1#2#3\autodef@end{%
  \ea\def\ea#1\ea{\ea\ensuremath\ea{\csname autofmt@#2\endcsname#3\autofmt@end}\xspace}}
\def\autodefs@end{blarg!}
\def\autodefs#1{\@autodefs#1\autodefs@end}
\def\@autodefs#1{\ifx#1\autodefs@end%
  \def\autodefs@next{}%
  \else%
  \def\autodefs@next{\autodef#1\@autodefs}%
  \fi\autodefs@next}

\DeclareSymbolFont{bbold}{U}{bbold}{m}{n}
\DeclareSymbolFontAlphabet{\mathbbb}{bbold}

\newcommand{\bbone}{\ensuremath{\mathbbb{1}}\xspace}

\let\del\partial
\mdef\delbar{\overline{\partial}}

\mdef\hf{\textstyle\frac12 }
\mdef\thrd{\textstyle\frac13 }
\mdef\qtr{\textstyle\frac14 }

\newcommand{\op}{^{\mathrm{op}}}

\SelectTips{cm}{}
\newdir{ >}{{}*!/-10pt/@{>}}    % extra spacing for tail arrows in XYpic
\newcommand{\pushoutcorner}[1][dr]{\save*!/#1+1.2pc/#1:(1,-1)@^{|-}\restore}
\newcommand{\pullbackcorner}[1][dr]{\save*!/#1-1.2pc/#1:(-1,1)@^{|-}\restore}

\mdef\Id{\mathrm{Id}}
\mdef\id{\mathrm{id}}
\alwaysmath{ell}
\alwaysmath{infty}
\alwaysmath{odot}
\def\frc#1/#2.{\frac{#1}{#2}}   % \frc x^2+1 / x^2-1 .
\mdef\ten{\mathrel{\otimes}}

\mdef\sqten{\mathrel{\boxtimes}}

\DeclareRobustCommand\widecheck[1]{{\mathpalette\@widecheck{#1}}}
\def\@widecheck#1#2{%
    \setbox\z@\hbox{\m@th$#1#2$}%
    \setbox\tw@\hbox{\m@th$#1%
       \widehat{%
          \vrule\@width\z@\@height\ht\z@
          \vrule\@height\z@\@width\wd\z@}$}%
    \dp\tw@-\ht\z@
    \@tempdima\ht\z@ \advance\@tempdima2\ht\tw@ \divide\@tempdima\thr@@
    \setbox\tw@\hbox{%
       \raise\@tempdima\hbox{\scalebox{1}[-1]{\lower\@tempdima\box
\tw@}}}%
    {\ooalign{\box\tw@ \cr \box\z@}}}

\DeclareMathOperator\colim{colim}

\DeclareMathOperator\im{im}

\DeclareMathOperator\Ho{Ho}

\DeclareMathOperator\Map{Map}

\newcommand{\ot}{\ensuremath{\leftarrow}}

\mdef\we{\overset{\sim}{\longrightarrow}}
\mdef\leftwe{\overset{\sim}{\longleftarrow}}

\def\rightarrowtailfill@{\arrowfill@{\Yright\joinrel\relbar}\relbar\rightarrow}
\newcommand\xrightarrowtail[2][]{\ext@arrow 0055{\rightarrowtailfill@}{#1}{#2}}

\def\twoheadrightarrowfill@{\arrowfill@{\relbar\joinrel\relbar}\relbar\twoheadrightarrow}
\newcommand\xtwoheadrightarrow[2][]{\ext@arrow 0055{\twoheadrightarrowfill@}{#1}{#2}}

\def\slashedarrowfill@#1#2#3#4#5{%
  $\m@th\thickmuskip0mu\medmuskip\thickmuskip\thinmuskip\thickmuskip
   \relax#5#1\mkern-7mu%
   \cleaders\hbox{$#5\mkern-2mu#2\mkern-2mu$}\hfill
   \mathclap{#3}\mathclap{#2}%
   \cleaders\hbox{$#5\mkern-2mu#2\mkern-2mu$}\hfill
   \mkern-7mu#4$%
}
\def\rightslashedarrowfill@{%
  \slashedarrowfill@\relbar\relbar\mapstochar\rightarrow}
\newcommand\xslashedrightarrow[2][]{%
  \ext@arrow 0055{\rightslashedarrowfill@}{#1}{#2}}
\mdef\hto{\xslashedrightarrow{}}
\mdef\htoo{\xslashedrightarrow{\quad}}

\long\def\my@drawfill#1#2;{%
\@skipfalse
\fill[#1,draw=none] #2;
\@skiptrue
\draw[#1,fill=none] #2;
}
\newif\if@skip
\newcommand{\skipit}[1]{\if@skip\else#1\fi}
\newcommand{\drawfill}[1][]{\my@drawfill{#1}}

\newif\ifhyperref
\@ifpackageloaded{hyperref}{\hyperreftrue}{\hyperreffalse}
\iftac
  \let\your@state\state
  \def\state#1{\gdef\currthmtype{#1}\your@state{#1}}
  \let\your@staterm\staterm
  \def\staterm#1{\gdef\currthmtype{#1}\your@staterm{#1}}
  \let\defthm\newtheorem
  \def\currthmtype{}
  \ifhyperref
    \def\autoref#1{\ref*{label@name@#1}~\ref{#1}}
  \else
    \def\autoref#1{\ref{label@name@#1}~\ref{#1}}
  \fi
  \AtBeginDocument{%
    \let\old@label\label%
    \def\label#1{%
      {\let\your@currentlabel\@currentlabel%
        \edef\@currentlabel{\currthmtype}%
        \old@label{label@name@#1}}%
      \old@label{#1}}
  }
\else
  \ifhyperref
    \def\defthm#1#2{%
      \newtheorem{#1}{#2}[section]%
      \expandafter\def\csname #1autorefname\endcsname{#2}%
      \expandafter\let\csname c@#1\endcsname\c@thm}
  \else
    \def\defthm#1#2{\newtheorem{#1}[thm]{#2}}
    \ifx\SK@label\undefined\let\SK@label\label\fi
    \let\old@label\label
    \let\your@thm\@thm
    \def\@thm#1#2#3{\gdef\currthmtype{#3}\your@thm{#1}{#2}{#3}}
    \def\currthmtype{}
    \def\label#1{{\let\your@currentlabel\@currentlabel\def\@currentlabel%
        {\currthmtype~\your@currentlabel}%
        \SK@label{#1@}}\old@label{#1}}
    \def\autoref#1{\ref{#1@}}
  \fi
\fi

\newtheorem{thm}{Theorem}[section]

\defthm{cor}{Corollary}
\defthm{prop}{Proposition}
\defthm{lem}{Lemma}
\defthm{sch}{Scholium}
\defthm{assume}{Assumption}
\defthm{claim}{Claim}
\defthm{conj}{Conjecture}
\defthm{hyp}{Hypothesis}
\defthm{fact}{Fact}
\iftac\theoremstyle{plain}\else\theoremstyle{definition}\fi
\defthm{defn}{Definition}
\defthm{notn}{Notation}
\iftac\theoremstyle{plain}\else\theoremstyle{remark}\fi
\defthm{rmk}{Remark}
\defthm{eg}{Example}
\defthm{egs}{Examples}
\defthm{ex}{Exercise}
\defthm{ceg}{Counterexample}
\defthm{warn}{Warning}
\defthm{per}{Perspective}
\defthm{con}{Construction}
\defthm{dig}{Digression}

\def\thmqedhere{\expandafter\csname\csname @currenvir\endcsname @qed\endcsname}

\setitemize[1]{leftmargin=2em}
\setenumerate[1]{leftmargin=*}

\iftac
  \let\c@equation\c@subsection
\else
  \let\c@equation\c@thm
\fi
\numberwithin{equation}{section}

\@ifpackageloaded{mathtools}{\mathtoolsset{showonlyrefs,showmanualtags}}{}

\alwaysmath{alpha}
\alwaysmath{beta}
\alwaysmath{gamma}
\alwaysmath{Gamma}
\alwaysmath{delta}
\alwaysmath{Delta}
\alwaysmath{epsilon}
\mdef\ep{\varepsilon}
\alwaysmath{zeta}
\alwaysmath{eta}
\alwaysmath{theta}
\alwaysmath{Theta}
\alwaysmath{iota}
\alwaysmath{kappa}
\alwaysmath{lambda}
\alwaysmath{Lambda}
\alwaysmath{mu}
\alwaysmath{nu}
\alwaysmath{xi}
\alwaysmath{pi}
\alwaysmath{rho}
\alwaysmath{sigma}
\alwaysmath{Sigma}
\alwaysmath{tau}
\alwaysmath{upsilon}
\alwaysmath{Upsilon}
\alwaysmath{phi}
\alwaysmath{Pi}
\alwaysmath{Phi}
\mdef\ph{\varphi}
\alwaysmath{chi}
\alwaysmath{psi}
\alwaysmath{Psi}
\alwaysmath{omega}
\alwaysmath{Omega}

\let\la\lambda

\makeatother

\usepackage{comment}

\tikzset{lab/.style={auto,font=\scriptsize}} % arrow labels
\usetikzlibrary{arrows}
\usepackage{eucal}

\usepackage{fixme}
\fxsetup{
    status=draft,
    author=,
    layout=margin,
    theme=color
}

\definecolor{fxnote}{rgb}{1.0000,0.0000,0.0000}
\colorlet{fxnotebg}{yellow}

\newcommand{\D}{\sD}

\autodefs{\cDER\cCat\cGrp\cCAT\cMONCAT\cSet\cAb\cGrpd\nFun\cGraph\cProf\cTop\cSet\cKan\nSing\nAdj
\cPro\cMONDER\cBICAT\nInd\nCh\nDK\nAlg\nMon\nGrp\nlax\ncf\ncomp\nid\ncoker\niso\nIm\ncoh\nLan\nRan\nRKan\ncolim\nobj\nMon\nev\nSp\nPSp\ncof\nfib\sProf\ncyl\nStab\fC\fZ\fT\nhocolim\nholim\cPsNat\ncoll\csSet
\csCat\cAlg\cIm\cI\cP\cG}

\let\oldboxtimes\boxtimes
\def\boxtimes{\mathrel{\oldboxtimes}}

\newcommand{\fib}{\mathsf{fib}}
\newcommand{\cof}{\mathsf{cof}}
\newcommand{\sSet}{s\cSet}
\newcommand{\sCat}{s\cCat}
\newcommand{\sAb}{s\cAb}
\newcommand{\nFunL}{\nFun^\mathrm{L}}

\newcommand{\MonE}{\nMon_{\mathbb{E}_\infty}\!}
\newcommand{\GrpE}{\nGrp_{\mathbb{E}_\infty}\!}
\newcommand{\AlgA}{\nAlg_{\mathbb{A}_\infty}\!}
\newcommand{\AlgE}{\nAlg_{\mathbb{E}_\infty}\!}
\newcommand{\Fin}{\mathcal{F}\mathrm{in}}
\renewcommand{\la}{\langle}
\newcommand{\ra}{\rangle}

\def\ccsub{_{\mathrm{cc}}}
\def\pdh(#1,#2){\llbracket #1,#2\rrbracket}
\def\ldh(#1,#2){\llbracket #1,#2\rrbracket\ccsub}
\def\pend(#1){\pdh(#1,#1)}
\def\lend(#1){\ldh(#1,#1)}

\def\DTl#1#2#3#4#5#6#7{%
  \xymatrix@C=3pc{{#1} \ar[r]^-{#2} &
    {#3} \ar[r]^-{#4} &
    {#5} \ar[r]^-{#6} &
    {#7}
  }}

\newsavebox{\tvabox}
\savebox\tvabox{\hspace{1mm}\begin{tikzpicture}[>=latex',baseline={(0,-.18)}]
  \draw[->] (0,.1) -- +(1,0);
  \node at (.5,0) {$\scriptscriptstyle\bot$};
  \draw[->] (1,-.1) -- +(-1,0);
  \draw[->] (1,-.2) -- +(-1,0);
\end{tikzpicture}\hspace{1mm}}

\newcommand{\Mod}{\mathrm{Mod}}

\newcommand{\Aoo}{\mathbb{A}_\infty}
\newcommand{\Eoo}{\mathbb{E}_\infty}
\newcommand{\FunL}{\mathrm{Fun^L}}
\newcommand{\FunR}{\mathrm{Fun^R}}
\newcommand{\FunLL}{\mathrm{Fun^{L,L}}}
\newcommand{\FunLo}{\mathrm{Fun^{L,\otimes}}}

\newcommand{\Prl}{\mathcal{P}\mathrm{r^L}}

\newcommand{\Prlpt}{\mathcal{P}\mathrm{r^L_{Pt}}}
\newcommand{\Prlpre}{\mathcal{P}\mathrm{r^L_{Pre}}}
\newcommand{\Prladd}{\mathcal{P}\mathrm{r^L_{Add}}}
\newcommand{\Prlst}{\mathcal{P}\mathrm{r^L_{St}}}
\newcommand{\Prlo}{\mathcal{P}\mathrm{r^{L,\otimes}}}
\newcommand{\Prlpto}{\mathcal{P}\mathrm{r^{L,\otimes}_{Pt}}}

\newcommand{\Prlsto}{\mathcal{P}\mathrm{r^{L,\otimes}_{St}}}
\newcommand{\Prlsmon}{\mathcal{P}\mathrm{r^{L,smon}}}
\newcommand{\Prlstsmon}{\mathcal{P}\mathrm{r^{L,smon}_{St}}}

\title{A short course on $\infty$-categories}
\author[M. Groth]{Moritz Groth}
\address{Moritz Groth, MPIM, Vivatsgasse 7, 53111
Bonn, Germany}
\email{mgroth@mpim-bonn.mpg.de}
\urladdr{http://guests.mpim-bonn.mpg.de/mgroth/}

\date{\today}

\begin{document}

\begin{abstract}
In this short survey we give a non-technical introduction to some main ideas of the theory of $\infty$-categories, hopefully facilitating the digestion of the foundational work of Joyal and Lurie. Besides the basic $\infty$-categorical notions leading to presentable $\infty$-categories, we mention the Joyal and Bergner model structures organizing two approaches to a theory of $(\infty,1)$-categories. We also discuss monoidal $\infty$-categories and algebra objects, as well as stable $\infty$-categories. These notions come together in Lurie's treatment of the smash product on spectra, yielding a convenient framework for the study of $\Aoo$-ring spectra, $\Eoo$-ring spectra, and Derived Algebraic Geometry.
\end{abstract}

\maketitle
\setcounter{section}{-1}

\tableofcontents

\section{Introduction}
\label{sec:intro}

The aim of this short course is to give a non-technical account of some ideas in the theory of $\infty$-categories (aka.~quasi-categories, inner Kan complexes, weak Kan complexes, Boardman complexes, or quategories), as originally introduced by Boardman--Vogt~\cite[p.102]{boardman-vogt} in their study of homotopy-invariant algebraic structures. Recently, $\infty$-categories were studied intensively by Joyal~\cite{joyal:I-II,joyal:quasi-kan,joyal:barca}, Lurie~\cite{HTT,DAGI,DAGII,DAGIII}, and others, and have applications in many areas of pure mathematics. Here we try to emphasize the philosophy and some of the main ideas of $\infty$-category theory, and we sketch the lines along which the theory is developed. \emph{In particular, this means that there is no claim of originality.}
 
Category theory is an important mathematical discipline in that it provides us with a convenient language which applies whenever we put into practice the following slogan: `In order to study a collection of objects one should also consider suitably defined morphisms between such objects.' Many classes of mathematical objects like groups, modules over a ring, manifolds, or schemes can be organized into a category and from typical constructions one frequently abstracts the categorical character behind them. Let us recall that a category consists of objects, morphisms, and a composition law which is suitably associative and unital. This allows us, in particular, to speak about \emph{isomorphisms}, and all functorial constructions trivially preserve isomorphisms. 

Category theory is a very powerful and useful language, however, it also has its limitations. Namely, to put it as a slogan, in many areas of pure mathematics we would like to identify two objects which are, while possibly not isomorphic in the purely categorical sense, `equivalent' from a more homotopy theoretic perspective. To illustrate this, the desire of identifying different resolutions of objects in abelian categories leads us to study chain complexes up to quasi-isomorphisms. Similarly, in homotopy theory we would like to think of weak homotopy equivalences between topological spaces as actual isomorphisms. And even in category theory, often we do not have to distinguish two categories as long as they are equivalent.

These are only three examples for the fairly common situation that we start with a pair $(\cC,W)$ consisting of a category~$\cC$ and a class~$W$ of so-called \emph{weak equivalences}, a class of morphisms which we would like to treat as isomorphisms. In such situations, functorial constructions are only `meaningful' if they preserve weak equivalences. The search for convenient languages to study such situations has already quite some history and various different approaches have been considered. This includes triangulated categories, model categories, derivators, simplicial categories, topological categories, and $\infty$-categories. 

The three last named approaches belong to a fairly large zoo of different models all of which realize `a theory of $(\infty,1)$-categories'. While an $(\infty,1)$-category is not a well-defined mathematical notion, there is the general agreement that such a category-like concept should enjoy the following features.
\begin{enumerate}
\item As part of the structure there is a class of objects.
\item There should be morphisms between objects, $2$-morphisms between morphisms (like chain homotopies, homotopies, and natural isomorphisms), as well as $3$-morphisms, $4$-morphisms, and so forth, explaining the parameter `$\infty$' in `$(\infty,1)$-categories'.
\item Morphisms can be composed in a suitably associative and unital way.
\item Higher morphisms, i.e., $2$-morphisms, $3$-morphisms, and so on, are supposed to be invertible in a certain sense. (All morphisms above dimension one are invertible, explaining the parameter `$1$' in `$(\infty,1)$-categories'.) 
\end{enumerate}
References for survey articles on the zoo of different axiomatizations of $(\infty,1)$-categories include \cite{bergner:workshop} and the fairly recent \cite  {camerona:whirlwind}.

The aim of these notes instead is to focus essentially on one of those different models and to describe how a good deal of classical category theory can be extended to this particular model. In these notes we follow Lurie~\cite{HTT,HA} in his choice of terminology and refer to these particular models for $(\infty,1)$-categories as \emph{$\infty$-categories}. We refrain from giving a more detailed introduction here and instead refer the reader to the table of contents as well as to the short introductions of the individual sections. To conclude this introduction there are the following few remarks.

\begin{enumerate}
\item In these notes we ignore essentially all set-theoretic issues (with the exception of the discussion of locally presentable categories where some care is needed).
\item  For many of the mathematical concepts to be introduced below, there are at least two different terminologies (most frequently, one due to Joyal and one due to Lurie). Since we do not want to cause further confusion, we have to stick to one of these possible choices. As the expanded version of Lurie's thesis \cite{HTT,HA} is our main reference, we stick to Lurie's terminology.
\item In these notes we use the environment `Perspective' in order to include some sketches of `the larger picture' and also to give more references to the literature. We tried our best to write the text such that there is essentially no loss of continuity if the reader skips these `Perspectives' on a first reading. However, we hope them to be helpful on a second reading.
\item Finally, the main prerequisites for these notes are some acquaintance with key concepts from category theory \cite{MacLane} as well as basics concerning simplicial sets. References for simplicial sets include the monographs~\cite{gabriel-zisman:calculus,goerss-jardine} and the more introductory account~\cite{friedman:survey}.
\end{enumerate}

\section{Two models for $(\infty,1)$-categories}
\label{sec:basics}

In this section we define \emph{$\infty$-categories} as simplicial sets satisfying certain horn extension properties. These extension properties are a common generalization of extension properties enjoyed by singular complexes of topological spaces and nerves of ordinary categories. We indicate why this notion gives us a model for the theory of $(\infty,1)$-categories. While $\infty$-categories have mapping spaces which can be endowed with a coherently associative and unital composition law, there is also the more rigid approach to $(\infty,1)$-categories based on \emph{simplicial categories} coming with strictly associative and unital composition laws. These two approaches are respectively organized by means of model structures, the \emph{Joyal model structure} in the case of $\infty$-categories and the \emph{Bergner model structure} in the case of simplicial categories. The coherent nerve construction of Cordier can be shown to be part of a Quillen equivalence between these two approaches. 

\subsection{Basics on $\infty$-categories}

Before we give the central definition of an $\infty$-category, we consider two classes of examples, which one definitely wants to be covered by the definition. The actual definition of an $\infty$-category will then be a common generalization of these two classes of examples. The first class of examples comes from spaces.

\begin{eg}
Given a topological space $X$, recall that associated to $X$ there is the \emph{fundamental groupoid} $\pi_1(X)$ of $X.$ The objects of $\pi_1(X)$ are the points of~$X$, and morphisms from $x$ to $y$ are homotopy classes of paths from $x$ to $y$ relative to the boundary points. Note that this is a groupoid, i.e., that all morphisms are invertible, since every path admits an inverse up to homotopy. This fundamental groupoid only depends on the 1-type of~$X$, and hence discards a lot of information. A refined version is given by the \textbf{fundamental} $\infty$-\textbf{groupoid} $\pi_{\leq\infty}(X)$ which is roughly constructed as follows: objects are given by points of $X$, morphisms are paths in~$X$, 2-morphisms are homotopies between paths, and higher morphisms are given by higher homotopies. 
\end{eg}

Note that $\pi_{\leq \infty}(X)$ seems to be an $\infty$-groupoid in that all morphisms are equivalences, i.e., invertible in a certain weak sense --- this justifies that we refer to it as fundamental \emph{$\infty$-groupoid} as opposed to as fundamental \emph{$\infty$-category}. The following is a generally accepted principle of higher category theory.
\[
\text{`\emph{Spaces and $\infty$-groupoids should be the same.}'}
\] 
This principle is referred to as the \emph{Grothendieck homotopy hypothesis}. Instead of working with topological spaces, one could also consider `simplicial models for spaces', more specifically \emph{Kan complexes} (see \autoref{def:Kan}). If one formalizes $\infty$-categories in the framework of simplicial sets (by way of \autoref{def:oocat}), then the above principle tells us that $\infty$-groupoids should be the same as Kan complexes. \autoref{cor:Kangroupoid} turns this principle into an actual mathematical statement. Using this approach to $\infty$-categories, a model for the fundamental $\infty$-groupoid of a space~$X$ is the usual singular complex $\mathrm{Sing}(X)$, which is well-known to be a Kan complex.

Let us establish some basic notation related to simplicial sets. Recall that the category $\sSet$ of simplicial sets is defined as the category of functors $\Delta\op\to\cSet$ where $\Delta$ is the category of finite ordinals 
\[
[n]=(0<\ldots <n),\quad n\geq 0,
\]
and order-preserving maps. Prominent maps in $\Delta$ are the coface maps $d^k$ and codegeneracy maps $s^k.$ The map $d^k\colon[n-1]\to[n],$ $0\leq k\leq n,$ is the unique injective map which does not have $k$ in its image, while $s^k\colon[n+1]\to[n],0\leq k\leq n,$ is the unique surjective map which hits $k$ twice. We follow the usual convention and write $X_n=X([n])$ for the values of a simplicial set $X$, $d_k=X(d^k)$ for the face maps, and $s_k=X(s^k)$ for the degeneracy maps.

Obviously, ordinary categories should be examples of $\infty$-categories: we imagine all higher morphisms to be identities. To make this precise, we recall the \emph{nerve construction} which associates a simplicial set to a category.

\begin{eg}
Given a category $\cC$, one can form the simplicial set~$N(\cC)\in\sSet,$ called the \textbf{nerve} of $\cC$. By definition, we have $N(\cC)_n=\nFun([n],\cC)$, where $[n]$ also denotes the ordinal number $0<\ldots<n$ considered as a category and $\nFun(-,-)$ is the set of functors. Thus $N(\cC)_n$ is essentially the set of strings of $n$ composable morphisms in $\cC$. As a special case we have $N([n])\cong\Delta^n, n\geq 0,$ where $\Delta^n\in\sSet$ denotes the usual simplicial $n$-simplex, i.e., $\Delta^n=\hom_{\Delta}(-,[n])$ is the simplicial set represented by $[n]\in\Delta$.
\end{eg}

Let us emphasize that in the definition of the nerve we consider the categories freely generated by the ordinals $[n]$. In particular, $N(\cC)_2$ consists of functors defined on the category
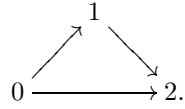
\[
\xymatrix@=1.5pc{
&1\ar[dr]&\\
0\ar[ru]\ar[rr]&&2.
}
\]
This picture makes clear that the first face map $d_1\colon N(\cC)_2\to N(\cC)_1$ is essentially given by composition, and that, strictly speaking, a $2$-simplex in $N(\cC)$ is a pair of composable morphisms together with their composition. We denote by $\cCat$ the category of categories. 

\begin{lem}
The nerve functor $N\colon\cCat\to\sSet$ is fully faithful and hence induces an equivalence onto its essential image.
\end{lem}

In order to describe the essential image one observes that nerves of categories enjoy certain \emph{horn extension properties}. Recall that the $k$-th $n$-\textbf{horn} $\Lambda^n_k\subseteq\del\Delta^n$ for $n\geq 1, 0\leq k\leq n,$ is obtained from $\del\Delta^n$ by removing the $k$-th face $\del_k\Delta^n$, i.e., the face opposite to vertex~$k$. More formally, the horn $\Lambda^n_k$ is defined as the following coequalizer
\[
\bigsqcup _{0\leq i<j \leq n}\hspace{-1.0em}\Delta ^{n-2}\rightrightarrows \bigsqcup _{i\neq k} \Delta ^{n-1} \rightarrow \Lambda^n_k.
\]
(See for example \cite[p.9]{goerss-jardine} which is a great reference for many more advanced aspects of simplicial homotopy theory.) 

In dimension $n=2$, horns $\lambda\colon\Lambda^2_k\to N(\cC)$ for $0\leq k\leq 2$ hence respectively look like
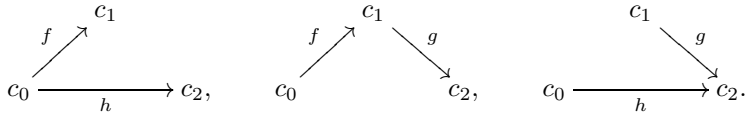
\[
\xymatrix@=1.5pc{
&c_1&&
&c_1\ar[dr]^-g&&
&c_1\ar[dr]^-g&\\
c_0\ar[ru]^-f\ar[rr]_-h&&c_2,&
c_0\ar[ru]^-f&&c_2,&
c_0\ar[rr]_-h&&c_2.
}
\]
Using the composition $h=g\circ f$ we see that one can uniquely extend any horn $\lambda\colon\Lambda^2_1\to N(\cC)$ to an entire 2-simplex $\sigma\colon\Delta^2\to N(\cC)$, i.e., there is a unique dashed arrow making the diagram
\[
\xymatrix@=1.5pc{
\Lambda^2_1\ar[r]\ar[d]&N(\cC)\\
\Delta^2\ar@{-->}[ru]_-{\exists !\sigma}
}
\]
commute. The composition is given by the new face $d_1(\sigma)\colon\Delta^1\to N(\cC)$. If instead we consider a horn $\lambda\colon\Lambda^2_0\to N(\cC)$ in the special case that $h=\id$ is an identity morphism, then the existence of an extension to a 2-simplex is equivalent to the existence of a left inverse to~$f$. Similar observations can be made for horns $\lambda\colon\Lambda^2_2\to N(\cC).$ This different behavior extends to higher dimensions and motivates the following terminology: the horns $\Lambda^n_k,0<k<n,$ are \textbf{inner horns} while the extremal cases $\Lambda^n_0$ and $\Lambda^n_n$ are \textbf{outer horns}. It turns out that these horn extension properties are suitable to describe the essential image of the nerve functor. We denote by $\cGrpd$ the category of groupoids.

\begin{prop}
Let $X$ be a simplicial set.
\begin{enumerate}
\item There is an isomorphism $X\cong N(\cC),\cC\in\cCat,$ if and only if \emph{every inner horn} $\Lambda^n_k\to X,0<k<n,$ can be \emph{uniquely} extended to an $n$-simplex $\Delta^n\to X.$ \item There is an isomorphism $X\cong N(\cG),\cG\in\cGrpd,$ if and only if \emph{every horn} $\Lambda^n_k\to X,0\leq k\leq n,$ can be \emph{uniquely} extended to an $n$-simplex $\Delta^n\to X.$
\end{enumerate}
\end{prop}

The characterization of the essential image of the nerve functor $N\colon\cGrpd\to\sSet$ reminds us of the notion of a Kan complex whose definition will be recalled here for convenience.

\begin{defn}\label{def:Kan}
A simplicial set $X$ is a \textbf{Kan complex} if \emph{every} horn $\Lambda^n_k\to X$ for $0\leq k\leq n$ can be extended to an $n$-simplex $\Delta^n\to X.$
\end{defn}

Denoting by $\cKan\subseteq\sSet$ the full subcategory spanned by the Kan complexes, we thus have the following commutative diagram of fully faithful functors
\begin{equation}
\vcenter{
\xymatrix@=1.5pc{
\cGrpd\ar[r]\ar[d]_-N&\cCat\ar[d]^-N\\
\cKan\ar[r]&\sSet.
}
}\label{eq:nerve}
\end{equation}

As a summary, Kan complexes and nerves of categories are simplicial sets with certain horn extension properties, but these properties differ in two important aspects. 
\begin{enumerate}
\item First, in a Kan complex \emph{all} horns can be extended, while in the nerve of a category this is, in general, only the case for \emph{inner} horns. 
\item Second, for Kan complexes we have a mere \emph{existence} statement while for nerves of categories the extensions are \emph{unique}. This uniqueness property is something one wants to drop for $\infty$-categories: it suffices that compositions exist and that the actual choice of compositions is `homotopically irrelevant'. This is similar to the concatenation of paths in a topological space: there is no preferred composition law for paths parametrized by the unit interval that would be associative and unital. But paths can be glued and the actual choice of the parametrization is `homotopically irrelevant', e.g., all choices are homotopic. 
\end{enumerate}

\begin{defn}\label{def:oocat}
A simplicial set $\cC$ is an \textbf{$\infty$-category} if every inner horn $\Lambda^n_k\to \cC,$ $0<k<n$, can be extended to an $n$-simplex $\Delta^n\to \cC$. 
\end{defn}

Thus, by the above, spaces and ordinary categories give rise to $\infty$-categories. We will soon see that there are interesting $\infty$-categories which are not of this form and hence provide some first `honest examples'. In particular, every simplicial model category has an underlying $\infty$-category but there are also other examples (see \autoref{cor:fibrantsimplicial}, \autoref{egs:oocats}, and \autoref{per:inftyinfty}).

We begin by introducing some basic terminology. Given an $\infty$-category $\cC$, the \textbf{objects} are the vertices $x\in\cC _0$ and the \textbf{morphisms} are the 1-simplices $f\in\cC_1$. The face map $s=d_1\colon\cC_1\to\cC_0$ is the \textbf{source map}, and $t=d_0\colon\cC_1\to\cC_0$ is the \textbf{target map}. As in ordinary category theory, we write $f\colon x\to y$ if $s(f)=x$ and $t(f)=y$. Slightly more formally, we define the set of morphisms $\hom_{\cC}(x,y)$ from~$x$ to~$y$ as the pullback
$$
\xymatrix@=1.5pc{
\hom_{\cC}(x,y)\pullbackcorner\ar[r]\ar[d]&\cC_1\ar[d]^{(s, t)}\\
\ast\ar[r]_{(x,y)}&\cC_0\times\cC_0.
}
$$
It turns out that associated to two objects in an $\infty$-category there is an entire \emph{space of morphisms} (see \autoref{rmk:principles}).

The degeneracy map $\id=s_0\colon\cC_0\to\cC_1$ is the \textbf{identity map}. It follows from the simplicial identities $d_0s_0=d_1s_0=\id_{\cC_0}$ that $\id_x$ is an endomorphism of $\cC$, i.e., that $\id_x=s_0x\colon x\to x$. We will see that every $\infty$-category~$\cC$ has an associated homotopy category. This is an ordinary category with the same objects and such that morphisms are represented by morphisms in~$\cC$. The morphism~$\id_x=s_0x$ turns out to represent the identity of~$x$ in the homotopy category, justifying the terminology and the notation (see \autoref{prop:Ho}). 

As for compositions, let us consider morphisms~$f\colon x\to y$ and~$g\colon y\to z$ in an $\infty$-category~$\cC$. These morphisms together define an inner horn in~$\cC$,
\[
\lambda=(g,\bullet,f)\colon \Lambda^2_1\to \cC,
\]
such that $d_0\lambda=g$ and $d_2\lambda=f$. Any such horn can be \emph{non-uniquely} extended to a 2-simplex $\sigma \colon\Delta^2\to \cC$. The new face $d_1(\sigma)$ opposite to vertex 1 is then a \emph{candidate composition} of $g$ and $f$. To re-emphasize, this is one of the central points in which $\infty$-category theory differs from ordinary category theory: one does not ask for uniquely determined compositions. Instead one demands only that there is a way to compose arrows and that any choice of such a composition is equally good: the space of all such choices is to be contractible (see the discussion of \autoref{thm:uniquecompo} for a precise statement).

We now want to describe the homotopy category of an $\infty$-category. This can be done in a more straightforward way, but we prefer to include a short digression in category theory as this allows us to mention a general fact which is in the background of a later construction anyhow. Nevertheless, the reader can skip this digression without loss of continuity and continue with \autoref{def:homotopic} instead.

\begin{dig}\label{dig:Yoneda} (\emph{Yoneda extension.})
Let $A$ be a small category and let us consider the associated presheaf category $\nFun(A\op,\cSet)$, i.e., the category of contravariant set-valued functors on $A$. Moreover, let~$\cC$ be a cocomplete category and let us be given a functor $Q\colon A\to\cC.$ Thus, we are in the situation
\[
\xymatrix@=1.5pc{
A\ar[r]^Q\ar[d]_y&\cC\\
\nFun(A\op,\cSet)\ar@{-->}@/_0.4pc/[ur]&
}
\]
where~$y$ denotes the Yoneda embedding of~$A.$ Recall from \cite[p.76]{MacLane} that every presheaf on a small category is canonically a colimit of representable ones (see also~\S\ref{subsec:locprescat}). The cocompleteness of~$\cC$ hence allows us to extend~$Q$ to a colimit-preserving functor 
\[
|-|_{Q}\colon\nFun(A\op,\cSet)\to\cC.
\]
(At a more conceptual level, we are thus forming the left Kan extension of~$Q$ along~$y.$) Moreover, associated to~$Q$, there is also a functor in the opposite direction
\[
\nSing_Q(-)\colon \cC\to\nFun(A\op,\cSet),
\]
which is defined by $\nSing_Q(c)_a=\hom_{\cC}(Qa,c).$ One observes now that this pair of functors defines an adjunction 
\[
(|-|_Q,\nSing_Q(-))\colon\nFun(A\op,\cSet)\rightleftarrows\cC
\]
with $|-|_Q$ as left adjoint and $\nSing_Q(-)$ as right adjoint. (Related to this see also the general discussion in \S\ref{subsec:locprescat} and in particular \autoref{thm:catSAFT}). 

But even more is true: For every cocomplete category~$\cC$ the assignment which sends $Q$ to the adjunction $(|-|_Q,\nSing_Q(-))$ defines an equivalence of categories 
\[
\nFun(A,\cC)\stackrel{\simeq}{\to}\nAdj\big(\nFun(A\op,\cSet),\cC\big).
\]
Here, $\nAdj(-,-)$ is the category of adjunctions where objects are adjunctions and morphisms are, say, natural transformations of left adjoints. This construction is sometimes referred to as \textbf{Yoneda extension} since it is essentially given by left Kan extension along the Yoneda embedding. A more detailed discussion of this can for example be found in \cite[pp.62-64]{Schapira}.
\end{dig}

In what follows we are only interested in the special case where $A=\Delta$, the category of finite ordinals. Thus, we conclude that a cosimplicial object $\Delta\to\cC$ in a cocomplete category is equivalently specified by an adjunction $\sSet\rightleftarrows\cC$. The notation employed in \autoref{dig:Yoneda} is, of course, motivated by the following example. 

\begin{eg}
Let $\cC=\cTop$ be the category of topological spaces and let us consider the standard cosimplicial space $|\Delta^\bullet|\colon\Delta\to\cTop.$ The associated adjunction is the usual adjunction given by the geometric realization and the singular complex functor, 
\[
(|-|,\nSing)\colon\sSet\rightleftarrows\cTop.
\]
\end{eg}

\begin{eg}
Let~$\cC=\cCat$ be the cocomplete category of small categories (see \autoref{rmk:VCatbicomplete} for the fact that $\cCat$ is cocomplete). The inclusion $\Delta\to\cCat$ obtained by considering a finite ordinal as a category induces by \autoref{dig:Yoneda} an adjunction
\[
(\tau_1,N)\colon\sSet\rightleftarrows\cCat
\] 
with right adjoint the nerve functor. The left adjoint~$\tau_1$ is the \textbf{fundamental category functor} or the \textbf{categorical realization functor}. The motivation for the first terminology and the notation~$\tau_1$ stems from the following. A composition with the groupoidification $\cCat\to\cGrpd$, i.e., the left adjoint to the forgetful functor $\cGrpd\to\cCat$ yields the usual \emph{fundamental groupoid functor} $\pi_1\colon\sSet\to\cGrpd.$ Thus associated to the cosimplicial object $\Delta\to\cGrpd$ which sends $[n]$ to the free groupoid on $[n]$ we obtain the adjunction 
\[
(\pi_1,N)\colon\sSet\rightleftarrows\cGrpd.
\]
\end{eg}

By the above abstract nonsense, the functor $\tau_1\colon\sSet\to\cCat$ sends a simplicial set $X$ to the category
\[
\tau_1(X)=\colim_{(\Delta/X)}[-]\circ p.
\]
Here, $(\Delta/X)$ is the \textbf{simplex category} of $X$ and $p\colon(\Delta/X)\to\Delta$ the canonical projection functor sending a simplex $([n],\Delta^n\to X)$ to~$[n].$ This description is, of course, not very explicit. However, using the skeletal filtration of~$X$, it follows that~$\tau_1(X)$ is obtained by the following two steps. First, one constructs the free category~$FX$ generated by the non-degenerate 1-simplices of~$X.$ The $2$-simplices of~$X$ define a congruence relation on~$FX$ which is generated as follows: Given a $2$-simplex $\sigma\colon\Delta^2\to X$ with boundary $\del \sigma =(g,h,f)$ as depicted in
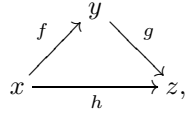
\[
\xymatrix@=1.5pc{
&y\ar[dr]^g&\\
x\ar[ur]^f\ar[rr]_h&&z,
}
\]
the morphisms $g\circ f$ and~$h$ of~$FX$ are congruent. The fundamental category~$\tau_1(X)$ is obtained from~$FX$ by dividing out this congruence relation (for more details we refer to~\cite{gabriel-zisman:calculus}). In particular, a morphism in~$\tau_1(X)$ can be represented by a \emph{finite chain} of 1-simplices of~$X.$ 

If the simplicial set happens to be an $\infty$-category, then there is a further simplification of the description of the fundamental category. It turns out that morphisms can be represented by actual 1-simplices as we shall discuss now (\autoref{prop:Ho}). For that purpose, we need the following definition.

\begin{defn}\label{def:homotopic}
Two morphism $f,g\colon x\to y$ in an $\infty$-category~$\cC$ are \textbf{homotopic} (notation: $f \simeq g$) if there is a 2-simplex $\sigma\colon\Delta^2\to \cC$ with boundary $\del \sigma =(g,f,\id_x)$, i.e., the boundary looks like
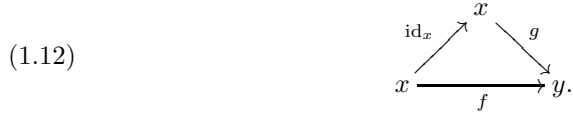
\begin{equation}
\vcenter{
\xymatrix@=1.5pc{
&x\ar[dr]^g&\\
x\ar[ur]^{\id _x}\ar[rr]_f&&y.
}
}
\label{eq:htpy}
\end{equation}
Any such 2-simplex $\sigma$ is a \textbf{homotopy} from $f$ to $g$, denoted $\sigma\colon f\to g.$
\end{defn}

There is a similar notion of homotopies where the identity morphism sits on the face opposite to vertex zero. However, using the inner horn extension property both notions can be shown to be the same. Moreover, there is the following result.

\begin{prop}
Let $\cC$ be an $\infty$-category and let $x,y\in\cC.$ The homotopy relation is an equivalence relation on $\hom_{\cC}(x,y)$. The \textbf{homotopy class} of a morphism $f\colon x\to y$ is denoted by~$[f]$.
\end{prop}

We include a partial proof in order to give an idea of how this works. Associated to a morphism $f\colon x\to y$, let us consider $\kappa_f=s_0f\colon\Delta^2\to\cC$. The simplicial identities imply that  $d_0\kappa_f=d_1\kappa_f=f$ and also $d_2\kappa_f=d_2s_0f=s_0d_1f=\id_x$. Thus, the boundary of $\kappa_f$ is $\partial\kappa_f=(f,f,\id_x)$ and we hence have a homotopy $\kappa_f\colon f\to f,$ the \textbf{constant homotopy} of~$f$, which establishes the reflexivity of $\simeq$. For the symmetry one way to proceed is as follows. Given a homotopy $\sigma\colon f\to g$, let us form the inner horn 
\begin{equation}
(\sigma, \kappa_g, \bullet, \kappa_{\id_x})\colon\Lambda^3_2\to \cC
\label{eq:inversehtpy}
\end{equation}
in $\cC.$ By definition of an $\infty$-category this horn can be extended to a 3-simplex $\tau\colon\Delta^3\to\cC$. The new face $\tilde{\sigma}=d_2\tau\in \cC_2$ defines a homotopy $\tilde{\sigma}\colon g\to f$, an \textbf{inverse homotopy} of~$\sigma$, showing that $\simeq$ is symmetric. 

With the homotopy relation at our disposal, we would like to define the \emph{homotopy category} $\Ho(\cC)$ of an $\infty$-category~$\cC$ by passing to homotopy classes of morphisms. The composition law in $\Ho(\cC)$ is obtained by representing homotopy classes by morphisms in $\cC$, \emph{choosing candidate compositions} of the representatives, and then passing to homotopy classes again. Of course, in order to get a well-defined category there are a lot of things to be checked, but we content ourselves by showing that all candidate compositions are homotopic. Let us consider morphisms $f\colon x\to y$ and $g\colon y\to z$ together with 2-simplices $\sigma_1,\sigma_2\colon\Delta^2\to \cC$ witnessing that $h_1=d_1(\sigma_1)$ and $h_2=d_1(\sigma_2)$ are candidate compositions of $g$ and $f$. Then we can form the inner horn 
\[
(\sigma_1,\sigma_2,\bullet,\kappa_f)\colon\Lambda^3_2\to \cC
\]
in~$\cC$. Again, we can find an extension to a 3-simplex $\tau\colon\Delta^3\to\cC$ and the new face $d_2\tau\colon\Delta^2\to\cC$ gives us the desired homotopy $h_2\to h_1$. Using similar arguments, one can establish the following result which can already be found in \cite{boardman-vogt}.

\begin{prop}\label{prop:Ho}
Let $\cC$ be an $\infty$-category. There is an ordinary category~$\Ho(\cC)$, the \textbf{homotopy category} of $\cC$, with the same objects as $\cC$ and morphisms the homotopy classes of morphisms in~$\cC$. Composition and identities are given by
\[
[g]\circ[f]:=[g\circ f]\qquad 
\text{and}\qquad\id_x:=[\id_x]=[s_0x],
\]
where $g\circ f$ is an \emph{arbitrary} candidate composition of $g$ and $f.$ Furthermore, there is a natural isomorphism of categories $\Ho(\cC)\cong\tau_1(\cC).$
\end{prop}

\begin{rmk}\label{rmk:principles}
\begin{enumerate}
\item One guiding principle for the theory of $\infty$-categories is that there should be a way to compose arrows and that the space of all such choices is contractible. Using the extension property for inner 2-horns, it is immediate that the space is non-empty. By means of the extension property for inner horns up to dimension three, we just checked that two candidate compositions are homotopic, i.e., that the space of all choices is connected. But this is only a $\pi_0$-statement of something much stronger: The extension property with respect to higher-dimensional inner horns can be thought of as guaranteeing the higher connectivity of the space of all such choices, giving finally that it is weakly contractible. See the discussion of \autoref{thm:uniquecompo} for a precise statement.
\item A second guiding principle for the theory of $\infty$-categories is that there should be morphisms of arbitrary dimensions. Let $\cC$ be an $\infty$-category and let $x,y$ be objects in $\cC.$ Then a morphism $f\colon x\to y$ is given by 
$$
f\colon \Delta^1\to \cC\qquad \text{such that}\qquad  f\hspace{-0.3em}\mid_{\Delta^{\{0\}}}=x\qquad\text{and}\qquad f\hspace{-0.3em}\mid_{\Delta^{\{1\}}}=y.
$$
Here and in the sequel, the notation is as follows: for vertices $i_0,\dots,i_k$ in $\Delta ^n$, $\Delta ^{\{i_0,\dots,i_k\}}\subseteq \Delta ^n$ denotes the $k$-simplex of $\Delta ^n$ spanned by the given vertices. A homotopy between two parallel morphisms $x\to y$ in $\cC$ can be interpreted as a \textbf{2-morphism} from $x$ to $y$. Recall that a homotopy is given by
$$
\sigma\colon\Delta ^2\to\cC\qquad \text{such that}\qquad \sigma\hspace{-0.3em}\mid_{\Delta^{\{0,1\}}}=x\qquad\text{and}\qquad \sigma\hspace{-0.3em}\mid_{\Delta^{\{2\}}}=y.
$$
This can be generalized to higher dimensions: an \textbf{$n$-morphism} from $x$ to $y$ is a map of simplicial sets 
$$
\tau\colon\Delta^{n+1}\to\cC\qquad\text{such that}\qquad \tau\hspace{-0.3em}\mid_{\Delta^{\{0,\dots,n\}}}=x\qquad\text{and}\qquad \tau\hspace{-0.3em}\mid_{\Delta^{\{n+1\}}}=y.
$$ 
For varying $n$, the sets of $n$-morphisms can be assembled in a \textbf{space of morphisms} $\Map_\cC^R(x,y)\in \sSet$ which can be shown to be a Kan complex.
 
We already mentioned that there is a variant to our definition of a homotopy (obtained by choosing the identity morphism in \eqref{eq:htpy} to sit opposite to vertex zero). More generally, there is an obvious dual way to define a space of morphisms $\Map^L_\cC(x,y)$ which turns out to be a weakly equivalent Kan complex. Thus, the homotopy type of the mapping space is well-defined. We will come back to this in the next subsection and give a conceptual explanation why these two models are weakly equivalent Kan complexes (see \autoref{per:mapping}).
\item A third guiding principle for the theory of $\infty$-categories is that they should give a model for $(\infty,1)$-categories, i.e., all higher morphisms should be invertible in some weak sense. To indicate that we succeeded in establishing such a framework, let us consider a homotopy~$\sigma$ in an $\infty$-category $\cC$,
$$
\sigma\colon f\simeq g\colon x\to y.
$$
In order to establish the symmetry of the homotopy relation, we considered the inner horn \eqref{eq:inversehtpy} which can be extended to a 3-simplex $\tau \colon\Delta ^3\to \cC.$ The new face $\widetilde{\sigma}=d_2\tau$ then gives us the intended homotopy $\widetilde{\sigma}\colon g\simeq f.$ Note that $\tau$ satisfies $\tau\!\!\mid_{\Delta^{\{0,1,2\}}}=x$ and $\tau$ is hence a \emph{3-morphism}, which can be interpreted as a 2-homotopy 
$$
\tau\colon\kappa _g\simeq \widetilde{\sigma}\circ \sigma.
$$ 
Thus, every homotopy has (up to a 2-homotopy) a left inverse and a similar observation can be made for right inverses. Taking for granted that the horn extension property for higher dimensional horns allows us to deduce similar observations for higher homotopies, we are reassured that $\infty$-categories really provide a model for $(\infty,1)$-categories.
\end{enumerate}
\end{rmk}

The following theorem due to Joyal makes precise that we succeeded in finding an axiomatic framework for categories with compositions determined up to contractible choices. Let $i\colon\Lambda^2_1\to\Delta^2$ be the obvious inclusion. Moreover, let us denote by $\Map(-,-)\colon\sSet\op\times\sSet\to\sSet$ the simplicial mapping space functor,
\begin{equation}
\Map(X,Y)_\bullet=\hom_{\sSet}(\Delta^\bullet\times X,Y),
\label{eq:ssetenriched}
\end{equation}
so that vertices are maps, edges are homotopies, and higher dimensional simplices are `higher homotopies' (see e.g.~\cite[p.20]{goerss-jardine}).

\begin{thm}\label{thm:uniquecompo}
A simplicial set $X$ is an $\infty$-category if and only if the restriction map $i^\ast\colon\Map(\Delta^2,X)\to \Map(\Lambda^2_1,X)$ is an acyclic Kan fibration. 
\end{thm}

We can think of $\Map(\Lambda^2_1,X)$ as the \emph{space of composition problems} and similarly of $\Map(\Delta^2,X)$ as the \emph{space of solutions to composition problems}. The theorem then tells us that the \emph{defining feature} of an $\infty$-category is that these two spaces are the same from a homotopical perspective. Let $f\colon x\to y$ and $g\colon y\to z$ be a pair of composable arrows in an $\infty$-category $\cC$ and let $\lambda=(g,-,f) \colon\Lambda^2_1\to \cC$ be the associated horn, i.e., a vertex $\lambda\colon\Delta^0\to\Map (\Lambda^2_1,\cC)$. The fiber $F_\lambda$ of $i^\ast$ over this vertex, i.e., the upper left corner in the pullback diagram
\[
\xymatrix@=1.5pc{
F_\lambda\pullbackcorner\ar[r]\ar[d] & \Map (\Delta^2,X)\ar[d]^-{i^\ast} \\
\Delta ^0  \ar[r]_-\lambda & \Map (\Lambda ^2_1,X)
}
\]
\noindent
can be regarded as the \emph{space of all possible compositions of $g$ and $f$}. By the theorem this space is a \emph{contractible Kan complex} which tells us that
\begin{enumerate}
\item a composition of $g$ and $f$ exists (the space of compositions is non-empty),
\item any two choices of compositions are homotopic (the space of compositions has trivial~$\pi_0$),
\item the homotopies comparing two such compositions are unique up to homotopy (the space of compositions has trivial $\pi_1$),
\item and so on in higher dimensions (also all $\pi_i,i\geq 2,$ vanish).
\end{enumerate}
This motivates us to henceforth suppress the `candidate' in `candidate composition'. 

We now turn to equivalences in an $\infty$-category.

\begin{defn}
A morphism $f\colon x\to y$ in an $\infty$-category $\cC$ is an \textbf{equivalence} if $[f]\colon x\to y$ is an isomorphism in $\Ho(\cC).$ 
\end{defn}

It is immediate that identities are equivalences and that for two homotopic morphisms $f_1\simeq f_2$ we have that $f_1$ is an equivalence if and only if $f_2$ is one. Moreover, it turns out that a morphism $f\colon x\to y$ in $\cC$ is an equivalence if and only if there is a morphism $g\colon y\to x$ in~$\cC$ such that there are 2-simplices with boundaries as in
\[
\xymatrix@=1.5pc{
&y\ar[dr]^g& && &x\ar[dr]^f&\\
x\ar[ur]^f\ar[rr]_{\id_x}&&x, && y\ar[ur]^g\ar[rr]_{\id_y}&&y.
}
\]
In a way one could be surprised that we can characterize equivalences in an $\infty$-category by these two conditions. Since $\infty$-category theory is some sort of `coherent category theory' one might have expected that also higher coherence data would be necessary to characterize equivalences in an $\infty$-category. For a precise statement and a proof of the equivalence of these two potentially different invertibility conditions we refer to \cite[Corollary~1.6]{joyal:quasi-kan} and \cite[Proposition~2.2]{dugger-spivak:mapping}.

We mentioned already the accepted principle that all $\infty$-\emph{groupoids} should come from spaces. In order to make this explicit let us give the following definition.

\begin{defn}
An $\infty$-category is an \textbf{$\infty$-groupoid} if the homotopy category is a groupoid.
\end{defn}

Thus an $\infty$-category is an $\infty$-groupoid if and only if all morphisms are equivalences. In the motivation of the definition of an $\infty$-category, we saw that, in general, one should only demand the horn extension property for inner horns in order to obtain a good generalization of arbitrary categories (and not just of groupoids!). Joyal established the following result, saying that outer horns can be extended as soon as certain maps are equivalences.

\begin{prop}
Let $\cC$ be an $\infty$-category. Any horn $\lambda\colon \Lambda^n_0\to \cC,\;n\geq 2,$ such that $\lambda \!\mid _{\Delta^{\{0,1\}}}$ is an equivalence can be extended to a simplex $\Delta^n\to\cC.$
\end{prop}

There is of course a similar statement using the horns $\Lambda^n_n$ instead. This allows us to turn the principle that all $\infty$-groupoids should be given by spaces into the following precise statement (see \cite[Corollary~1.4]{joyal:quasi-kan} or \cite[p.35]{HTT}).

\begin{cor}\label{cor:Kangroupoid}
An $\infty$-category is an $\infty$-groupoid if and only if it is a Kan complex.
\end{cor}

With this result at hand, diagram \eqref{eq:nerve} consisting of fully faithful functors can be refined to
\begin{equation}
\vcenter{
\xymatrix@=1.5pc{
\cGrpd\ar[r]\ar[d]_-N&\cCat\ar[d]^-N\ar@/^0.8pc/[rd]^-N&\\
\cKan=\cGrpd_\infty\ar[r]&\cCat_\infty\ar[r]&\sSet.
}
}
\end{equation}

\begin{per}
Similar to \autoref{cor:Kangroupoid}, for low values of $n\in\mathbb{N}$ there are statements using $n$\emph{-types of spaces} and $n$\emph{-groupoids} in a certain precise sense. The statements that higher homotopy types should be classified by higher groupoids is frequently also referred to as the \emph{homotopy hypothesis}. 

In the case of $n=1$, a precise statement can for example be found in \cite[Theorem 3.9]{hollander:htpy-thy-stacks} where it is shown that such a classification is induced by the adjunction $(\pi_1,N)\colon\sSet\rightleftarrows\cGrpd.$ In fact, this adjunction can be seen to be a Quillen adjunction with respect to the Kan--Quillen model structure on~$\sSet$ and the so-called \emph{natural model stucture} on~$\cGrpd$ (see \autoref{per:fibrations}). The slogan that `groupoids do not carry any higher homotopical information' can be made precise as follows: the Quillen adjunction $(\pi_1,N)$ induces a Quillen equivalence between the~$S^2$-nullification of~$\sSet$ and $\cGrpd.$ Related results in the cases of $n=2,3,$ more precisely, in the context of \emph{bicategories} \cite{benabou:intro} and \emph{Gray categories} \cite{gray:formal,leinster:higher}, are made explicit in \cite[\S6]{lack:gray-categories}.
\end{per}

\subsection{Simplicial categories and the relation to $\infty$-categories}
\label{subsec:joyalbergner}

There are many alternative approaches to a theory of $(\infty,1)$-categories including simplicial categories \cite{bergner:scat}, Segal categories \cite{hirschowitz-simpson}, and complete Segal spaces \cite{rezk:model}. Besides in the original references, more details can for example be found in \cite{bergner:survey,bergner:workshop,simpson:higher} and the nice recent \cite{camerona:whirlwind}. Here we include a short discussion of \textbf{simplicial categories} or, more precisely, \emph{simplicially enriched categories}. Given two objects $x,y$ in a simplicial category~$\cC,$ we write $\Map_{\cC}(x,y)$ for the associated simplicial mapping space. This more rigid approach --- coming with a specified strictly associative and unital composition law --- gives us, by definition, a notion of a category with morphisms of arbitrary dimensions. Building on work of Joyal and Bergner, Lurie has shown that this approach and the one using $\infty$-categories are equivalent in a very precise sense (see \autoref{thm:JoyalBergner} and \autoref{per:models}).

We begin by describing a relation between simplicial sets and simplicial categories. First, let us recall that the nerve $N(\cC)$ of an ordinary category $\cC$ is the simplicial set 
\[
N(\cC)_\bullet=\hom_{\cCat}([\bullet],\cC),
\]
where $[\bullet]\colon\Delta\to\cCat$ is obtained by considering the finite ordinals $[n]$ as categories. Given a \emph{simplicial} category $\cC$, we could simply forget the simplicial enrichment and form the nerve of the underlying ordinary category. More precisely, if we denote by $\sCat$ the category of (small) simplicial categories and simplicial functors, then there is the forgetful functor $\sCat\to\cCat$ which we could compose with the ordinary nerve functor $\cCat\to\sSet$. But this approach obviously discards too much information and instead one should proceed differently. 

A better way is given by replacing $[n]\in\cCat$ by \emph{simplicially thickened versions} $C[\Delta^n]\in\sCat$ and then building the simplicial set 
\[
N\hspace{-0.2em}_\Delta\hspace{-0.1em}(\cC)_\bullet=\hom_{\sCat}(C[\Delta^\bullet],\cC),
\] 
where $\hom_{\sCat}(-,-)$ denotes the set of simplicial functors. The idea behind this simplicial thickening is that $C[\Delta^n]$ encodes as objects the vertices of the standard simplex $\Delta^n$, as morphisms all paths in increasing direction, as $2$-morphisms all homotopies, and so on in higher dimensions. More conceptual comments about this construction can be found in \autoref{per:thickening}. Before we give a precise definition of $C[\Delta^\bullet]$ let us describe what we expect to obtain in low dimensions.

\begin{eg}\label{eg:fat-ordinals}
In dimensions 0 and 1 nothing new happens, and the simplicial categories $C[\Delta^0]$ and $C[\Delta^1]$ are just the ordinary categories $[0]$ and $[1]$, respectively, considered as simplicial categories with discrete mapping spaces. Thus, the pictures we have in mind are
\[
C[\Delta^0]\colon\quad 0\quad\qquad \text{and}\qquad C[\Delta^1]\colon\quad 0\to 1.
\] 
But from dimension 2 on the simplicial picture is richer. In~$\Delta^2$, there are two ways to pass from~0 to~2, namely the straight path and the path passing through 1. These paths should be encoded in~$C[\Delta^2]$ together with a homotopy between them. The simplicial category~$C[\Delta^2]$ can hence be depicted by
\begin{equation}
\vcenter{
\UseAllTwocells
\xymatrix@=1.5pc{
&1\ar[rd]&\\
0 \ar[ru] \ar[rr] & & 2.\lltwocell\omit{<2.5>}
}}\label{eq:cohtwo}
\end{equation}
\end{eg}

We now give a precise definition of $C[\Delta^n]$. The objects of $C[\Delta^n]$ are the numbers $0,1,\ldots,n$. The strategy behind the definition of the simplicial mapping spaces is the following. Given objects~$i\leq j$ we encode a path from $i$ to $j$ by specifying the vertices of the corresponding path. Thus, let $P_{i,j}$ be the poset
\[
P_{i,j}=\big\{I\subseteq[i,j]\;\mid\;i,j\in I\big \}
\]
ordered by inclusion where $[i,j]$ is short hand notation for $\{i,i+1,\ldots,j-1,j\}$. Considering these posets as categories, we can define the simplicial mapping spaces in $C[\Delta^n]$ by 
\[
\Map_{C[\Delta^n]}(i,j)=\left\{ \begin{array}{c@{\;,\;}l}
				NP_{i,j} & \quad i\leq j,\\ \emptyset & \quad  i > j.
				\end{array} \right. 
\]
The composition is induced by the union of subsets, which fits fine with the strategy to encode a path by specifying the vertices one passes along. It is also immediate that identities are given by the singletons $\{i\}.$ This concludes the definition of $C[\Delta^n]\in\sCat$.

One easily checks, that this definition specializes to the pictures we had in mind in low dimensions (\autoref{eg:fat-ordinals}). For example in dimension $n=2$, there is the following table of non-degenerate $k$-simplices in the mapping spaces $\Map_{C[\Delta^2]}(i,j)$:
\[
\renewcommand{\arraystretch}{1.5}
\begin{array}{c||c|c|c}
k & i=j=0 & i=0, j=1 & i=0, j=2\\ \hline 
0 & \{0\} & \{0,1\} & \{0,2\},\;\{0,1,2\} \\
1 & &  &\{0,2\}\subseteq\{0,1,2\}\\[0.5ex]
\end{array}
\]
By definition $\{0,1,2\}=\{1,2\}\circ\{0,1\}$ is the composition and we see that the non-degenerate 1-simplex in $\Map_{C[\Delta^2]}(0,2)$ encodes the homotopy 
\[
\{0,2\}\to\{1,2\}\circ\{0,1\}
\]
we were aiming for in \eqref{eq:cohtwo}.

It is straightforward to check that the assignment $[n]\mapsto C[\Delta^n]$ defines a cosimplicial object $C[\Delta^\bullet]\colon\Delta\to\sCat.$ This allows us to give the following definition which appears to be due to Cordier \cite{cordier:coherent}.

\begin{defn}\label{defn:cohnerve}
The \textbf{coherent nerve} $N\hspace{-0.2em}_\Delta\hspace{-0.1em}(\cC)$ of a simplicial category~$\cC$ is the simplicial set 
\[
N\hspace{-0.2em}_\Delta\hspace{-0.1em}(\cC)_\bullet=\hom_{\sCat}(C[\Delta^\bullet],\cC).
\] 
\end{defn}

Thus we have a coherent nerve functor $N_\Delta\colon\sCat\to\sSet.$ By the very definition, this coherent nerve construction takes into account the higher structure on a simplicial category given by the mapping spaces. For example a 2-simplex in such a coherent nerve is given by a homotopy as depicted in
\[
\UseAllTwocells
\xymatrix@=1.5pc{
&y\ar[rd]&\\
x \ar[ru] \ar[rr] & & z.\lltwocell\omit{<2.5>}
}
\]
Note that such a 2-simplex is, in general, not determined by its restriction to the horn $\Lambda^2_1\subset\Delta^2.$

Using the observation that $\sCat$ is cocomplete (see \autoref{rmk:VCatbicomplete}) we can extend the cosimplicial object $\Delta\to\sCat\colon[n]\mapsto C[\Delta^n]$ to a colimit-preserving functor $C[-]\colon \sSet\to\sCat.$ More explicitly, for $X\in\sSet$ we make the definition
\begin{equation}
C[X]=\colim_{(\Delta/X)}C[-]\circ p\label{eq:C}
\end{equation}
where $(\Delta/X)$ is the category of simplices of~$X$ and $p\colon(\Delta/X)\to\Delta$ is the canonical functor. It turns out that this extension defines a left adjoint to the coherent nerve~$N_\Delta,$
\begin{equation}\label{eq:CN-adj}
(C[-],N\hspace{-0.2em}_\Delta\hspace{-0.1em})\colon\sSet\rightleftarrows\sCat.
\end{equation}
\noindent
(In fact, this can be considered as an example of \emph{Yoneda extensions} in the sense of \autoref{dig:Yoneda}.) We observe that the notation $C[-]$ is not in conflict with the notation $C[\Delta^n]$ for the simplicial thickening of $[n]$ since the colimit-preserving extension $C[-]$ applied to $\Delta^n$ is isomorphic to what we just defined. This follows because the category of simplices $(\Delta/\Delta^n)$ has $([n],\id\colon\Delta^n\to\Delta^n)$ as terminal object so that the defining colimit in \eqref{eq:C} simplifies accordingly.

\begin{rmk}\label{rmk:VCatbicomplete}
The cocompleteness of $\sCat$ is a special instance of a more general result. Given a symmetric monoidal category~$\cM$ let us denote by~$\cCat_{\cM}$ the category of (small) $\cM$-enriched categories and $\cM$-enriched functors. Thus in this notation we have $\sCat=\cCat_{\sSet}.$ It turns out that if~$\cM$ is complete and cocomplete then so is~$\cCat_{\cM}.$ The harder part is the cocompleteness and was established by Wolff in \cite{wolff:v-cat}. As important examples we deduce that the categories of categories, simplicial categories, topological categories, spectral categories, differential-graded categories, and $2$-categories are complete and cocomplete.
\end{rmk}

It is a result due to Lurie that the adjunction \eqref{eq:CN-adj} is in fact a Quillen equivalence with respect to the \emph{Joyal model structure} on $\sSet$ and the \emph{Bergner model structure} on $\sCat$. Since we will not make an intensive use of the Bergner model structure, we do not go too much into detail and instead refer to \cite{bergner:scat}. Let us recall that every simplicial category $\cC$ has an underlying \textbf{path component category} or \textbf{homotopy category}~$\pi_0\cC$. This is an ordinary category with the same objects while the sets of morphisms are obtained by applying $\pi_0$ to the simplicial mapping spaces (see \autoref{rmk:basechange}). For example, if we endow $\sSet$ with the usual simplicial enrichment given by \eqref{eq:ssetenriched}, then $\pi_0\sSet$ is the naive homotopy category with all simplicial sets as objects and simplicial homotopy classes as morphisms. 

We can now define the weak equivalences in the Bergner model structure.

\begin{defn}\label{defn:DK}
A simplicial functor $F\colon\cC\to\cD$ is a \textbf{weak equivalence} if 
\begin{enumerate}
\item the induced functor $\pi_0 F\colon\pi_0\cC\to\pi_0\cD$ is essentially surjective and 
\item for all objects $x,y\in \cC$ the map $\Map_{\cC}(x,y) \to \Map_{\cD}(Fx,Fy)$ is a weak equivalence of simplicial sets (i.e., it induces a weak equivalence on geometric realizations).
\end{enumerate}
\end{defn}

Recall that a functor between ordinary categories is an equivalence if and only if it is essentially surjective and fully faithful. The definition of a weak equivalence between simplicial categories can be read as a higher categorical generalization of equivalences since it is asking that the simplicial functor is \emph{homotopically essentially surjective} and \emph{homotopically fully faithful.} Such a functor is also called a \textbf{Dwyer--Kan equivalence}, attributing credit to \cite{dwyer-kan:function}. Obviously, such a weak equivalence $\cC\to\cD$ induces an equivalence $\pi_0\cC\to\pi_0\cD$ but having a weak equivalence is, in general, a much stronger statement.

\begin{rmk}\label{rmk:basechange}
We want to include a short remark about enriched category theory (see \cite[p.313-316]{borceux:2} and as a general reference for that subject \cite{kelly:enriched}). Given a monoidal functor $G \colon\cM\to\cN,$ one obtains a \textbf{change of base functor} 
\[
G _*\hspace{-0.2em}=\cCat_{G}\colon\cCat_{\cM}\to \cCat_{\cN}
\] 
between the corresponding categories of enriched categories (see \autoref{rmk:VCatbicomplete} for the notation). This is simply a many-object-version of the fact that monoidal functors send monoid objects to monoid objects. In slightly more detail, if $\cC$ is an $\cM$-enriched category, then the $\cN$-enriched category $G_*\cC$ has the same objects as~$\cC$ while the mapping objects are given by 
\[
\Map_{G_*\cC}(x,y)=G\Map_{\cC}(x,y)\in\cN.
\]
It is an easy exercise to define compositions and identities in $G_\ast\cC.$ Moreover, there is a similar assignment for enriched functors (and also for enriched transformations but let us ignore that there are 2-categories in the background). 

As an example, the path component functor $\pi_0\colon\sSet\to\cSet$ is monoidal with respect to the Cartesian monoidal structures on both categories and the induced functor $\pi_0 F\colon\pi_0\cC\to\pi_0\cD$ used in \autoref{defn:DK} is simply $(\pi_0)_*F.$ Strictly speaking, also the formation of the underlying category of a simplicial category and the discrete simplicial category associated to an ordinary category are examples of this pattern (when applied to the corresponding product preserving functors $\sSet\to\cSet$ and $\cSet\to\sSet)$.
\end{rmk}

Building on work of Dwyer and Kan, Bergner \cite{bergner:scat} established the following result.

\begin{thm}\label{thm:bergner}
The category $\sCat$ carries a left proper combinatorial model structure with the Dwyer--Kan equivalences as weak equivalences. With respect to this model structure, a simplicial category is fibrant if and only if it is \emph{locally fibrant}, i.e., if all simplicial mapping spaces are Kan complexes.
\end{thm}

This model structure is referred to as the \textbf{Bergner model structure} and provides us with an example of a \emph{homotopy theory of homotopy theories}. For more details about it we refer to \cite{bergner:scat}.

\begin{per}\label{per:fibrations}
We want to include a brief discussion of the fibrations in the Bergner model structure since there are many other examples of model categories of, say, certain enriched categories in which the classes of fibrations are defined similarly. For this purpose let us begin by recalling that there is the `canonical' or `natural' model structure on $\cCat$ itself. This was established in a broader generality in \cite{joyal-tierney:stacks} while a nice short account can be found in \cite{rezk:natural-model}. The above adjectives refer to the fact that the weak equivalences are precisely the equivalences in the $2$-category of small categories. 

The fibrations in the natural model structure are the so-called \emph{isofibrations}. Recall that a functor $F\colon A\to B$ is an \textbf{isofibration} if every isomorphism $Fa\to b$ in~$B$ can be lifted to an isomorphism in~$A$ with domain given by~$a$. In other words, a functor is an isofibration if and only if it has the right lifting property with respect to the inclusion of the object $0$ in the groupoid $J$ generated by $[1]=(0\to 1),$
\[
\xymatrix@=1.5pc{
[0]\ar[r]\ar[d]_-0&A\ar[d]^-F\\
J\ar[r]\ar@{-->}[ru]_-\exists&B.
}
\]
The groupoid~$J$ is sometimes referred to as the \textbf{freely living isomorphism}. Note that there is an analogy to topology where fibrations have the path lifting property. Since paths in topological spaces are invertible we insist on isomorphisms in the categorical picture, and the definition of an isofibration can hence be read as asking for a lifting property with respect to invertible paths.

Now, a simplicial functor $F\colon\cC\to\cD$ is a \textbf{fibration} in the Bergner model structure if 
\begin{enumerate}
\item the ordinary functor $\pi_0F\colon\pi_0\cC\to\pi_0\cD$ is an isofibration and 
\item for all objects $x,y\in\cC$ the map $\Map_{\cC}(x,y)\to\Map_{\cD}(Fx,Fy)$ is a Kan fibration.
\end{enumerate}
Since the discrete simplicial category $C[\Delta^0]$ is the terminal object in $\sCat$ we immediately obtain the description of the fibrant objects as in \autoref{thm:bergner}. 

There are quite a few additional examples of model structures on $\cCat_{\cM}$ for nice monoidal model categories~$\cM$ which are similar in spirit. In these examples, an enriched functor is a weak equivalence if and only if it is a `Dwyer--Kan type equivalence'. Similarly, an enriched functor is a fibration if and only if it satisfies some `homotopical isofibration condition' and if it induces fibrations on all mapping objects. Specific results along these lines were established for differential-graded categories by Tabuada~\cite{tabuada:dgcats}, for spectral categories by Stanculescu~\cite{stanculescu:spectral} and Tabuada~\cite{tabuada:spectral}, for topological categories by Ilias~\cite{ilias:topological}, and in an axiomatic way by Lurie \cite[Appendix~A.3]{HTT}, Muro \cite{muro:dwyer-kan}, and Berger--Moerdijk \cite{berger-moerdijk:enriched}.

And this is not the end of the story. At least `morally similar model structures' have been established by Lack for 2-categories in~\cite{lack:model-2-categories}, for bicategories in~\cite{lack:model-bicategories}, and for Gray categories in~\cite{lack:gray-categories}. An extension in a further direction was given by Weiss in his thesis~\cite{weiss:thesis} who constructed a `natural model structure' on the category of small multicategories (aka.~colored operads or simply operads); see also \cite{weiss:folk}. This opens the door for attempts to establish model structures on certain categories of enriched multicategories which `mimic' this natural model structure. Related results were established in the simplicial context by Robertson~\cite{robertson:simplicial}, Cisinski--Moerdijk~\cite{cisinski-moerdijk:sop}, and Stanculescu~\cite{stanculescu:simplicial}, and in an axiomatic way by Caviglia~\cite{caviglia:enriched}. 
\end{per}

\begin{per}\label{per:thickening}
We have already seen that the functor $C[-]\colon\sSet\to\sCat$ is essentially determined by the cosimplicial object~$C[\Delta^\bullet]\colon\Delta\to\sCat$ given by the `simplicial thickenings' of the finite ordinals.  These simplicial thickenings arise more conceptually as follows. If we consider the categories $[n]$ as discrete simplicial categories, then we obtain a cosimplicial object $[\bullet]\colon\Delta\to\sCat$. The Bergner model structure of \autoref{thm:bergner} induces a Reedy model structure on the category $\nFun(\Delta,\sCat)$ of cosimplicial objects in $\sCat.$ It turns out that $C[\Delta^\bullet]\colon\Delta\to\sCat$ gives us a Reedy cofibrant replacement of $[\bullet]\colon\Delta\to\sCat.$  For an additional perspective on these simplicial thickenings we refer to \autoref{per:coherent}.
\end{per}

With a view towards the Joyal model structure on $\sSet$, we make the following definition.

\begin{defn}\label{defn:catequivalence}
A map $f\colon X\to Y$ in $\sSet$ is a \textbf{categorical equivalence} if the induced simplicial functor $C[f]\colon C[X]\to C[Y]$ is a Dwyer--Kan equivalence.
\end{defn}

This terminology is not the original one of Joyal. The maps in this definition are called \emph{weak categorical equivalences} by Joyal \cite{joyal:I-II}, while he has a stronger notion of categorical equivalence. However, his notions of categorical equivalence and weak categorical equivalence coincide when only maps between $\infty$-categories are considered.

For simplicity, a categorical equivalence between $\infty$-categories is called an \textbf{equivalence of $\infty$-categories} and we say that the $\infty$-categories are \textbf{equivalent}. The fact that it suffices to consider direct equivalences as opposed to more complicated zig-zags is a consequence of the following important theorem of Joyal~\cite{joyal:I-II}.

\begin{thm}\label{thm:joyal}
The category $\sSet$ carries a left proper combinatorial model structure with the monomorphisms as cofibrations and the categorical equivalences as weak equivalences. Moreover, a simplicial set is fibrant with respect to this model structure if and only if it is an $\infty$-category.
\end{thm}

The model structure of this theorem is referred to as the \textbf{Joyal model structure}. We want to mention that the model structure is \emph{Cartesian} which implies that there is an entire $\infty$-category of functors $\nFun(K,\cC)$ where~$K$ is a simplicial set and~$\cC$ an $\infty$-category (see \S\ref{subsec:functors}). On the negative side, this model structure is neither right proper nor simplicial. The latter drawback implies that one has to work harder in order to obtain an \emph{$\infty$-category of $\infty$-categories} (see \autoref{per:inftyinfty}).

We add some comments on the fibrations in the Joyal model structure. For this we first recall that in the Kan--Quillen model structure on $\sSet$ the fibrations are the Kan fibrations. A Kan fibration~$p\colon X\to S$ gives us a family of Kan complexes, i.e., $\infty$-groupoids, namely the fibers $X_s$ defined by the pullbacks
\[
\xymatrix@=1.5pc{
X_s\pullbackcorner\ar[r]\ar[d]&X\ar[d]^p\\
\Delta^0\ar[r]_s&S.
}
\]
Similarly, there is the following class of maps giving families of $\infty$-categories.

\begin{defn}\label{defn:inner-fibration}
A morphism of simplicial sets $p\colon X\to S$ is an \textbf{inner fibration} if it has the right lifting property with respect to $\Lambda^n_k\to\Delta^n,\:0<k<n.$
\end{defn}
\noindent
Joyal uses the term \emph{mid-fibration} instead of inner fibration. Since any class of morphisms defined by a right lifting property is closed under pullbacks this implies that the fibers of an inner fibration are $\infty$-categories. The dependence of the fiber on the base point is only functorial in a very weak sense, namely in the sense of \textbf{correspondences} (see \cite[p.97]{HTT}), which are also known as \textbf{distributors, profunctors} or \textbf{bimodules} in the classical setting (see \cite[\S7.8]{borceux:1}). 

Now, the \textbf{categorical fibrations}, i.e., the fibrations in the Joyal model structure happen to be a bit difficult to describe. However, for a morphism $p\colon X\to S$ such that the target~$S$ is an $\infty$-category, Joyal gave the following characterization: such a map is a categorical fibration if and only if it is an inner fibration and an \emph{isofibration}. More precisely, since $p$ is an inner fibration and $S$ is an $\infty$-category this is also the case for~$X$, and a map $X\to S$ between $\infty$-categories is an \textbf{isofibration} if and only if every equivalence $p(x)\to s$ in $S$ can be lifted to an equivalence in $X$ with domain~$x$. 

Given a functor $F\colon A\to B$ of ordinary categories, then $N(F)$ is automatically an inner fibration; thus the notion of inner fibrations does not have a classical analogue (in particular, a functor can always be thought of as a family of categories parametrized by the objects of the target category). In the $\infty$-categorical world however in many definitions this condition has to be imposed (this is the case for the categorical fibrations but also for further classes of fibrations as we will see later). 

As already mentioned, the original definition of categorical equivalences due to Joyal~\cite{joyal:I-II} is different. He gives a definition without reference to simplicial categories and his proof of the existence of the Joyal model structure is purely combinatorial. Lurie  gives this alternative definition because he is heading for the following comparison result \cite[p.89]{HTT}.

\begin{thm}\label{thm:JoyalBergner}
The adjunction $(C[-],N_{\Delta}\hspace{-0.1em})\colon\sSet\rightleftarrows\sCat$ is a Quillen equivalence with respect to the Joyal model structure and the Bergner model structure,
\[
(C[-],N_{\Delta}\hspace{-0.1em})\colon\sSet\stackrel{\stackrel{Q}{\sim}}{\to}\sCat.
\]
\end{thm}
A similar result can also be obtained by a combination of results due to Bergner, Joyal, Rezk, and Tierney (see \autoref{per:models}). \autoref{thm:JoyalBergner} makes precise in which sense the approaches to a theory of $(\infty,1)$-categories given by $\infty$-categories and simplicial categories are equivalent. The proof of this theorem is actually hard work including a deep `rigidification or straightening result' and can be found in~\cite{HTT}. An alternative proof was given by Dugger and Spivak in \cite{dugger-spivak:mapping, dugger-spivak:rigid}.

As a corollary, we have the following result which can also be obtained directly and without any mention of model structures (see for example the proof of \cite[Theorem 2.1]{cordier-porter:vogt}). However, with \autoref{thm:JoyalBergner} in mind the result is put into perspective. Recall that a simplicial category is called locally fibrant if all mapping spaces are Kan complexes.

\begin{cor}\label{cor:fibrantsimplicial}
The coherent nerve of a locally fibrant simplicial category is an $\infty$-category.
\end{cor}

We now turn to some typical examples of $\infty$-categories which are neither nerves of categories nor singular complexes of spaces. So, these are somehow our first honest examples of $\infty$-categories. It turns out that the first class of examples is generic in a sense which is made precise by \autoref{thm:comparison}.

\begin{egs}\label{egs:oocats}
\begin{enumerate}
\item Let $\cM$ be a simplicial model category \cite[\S II.2]{quillen} and let $\cM_{\ncf}\subseteq\cM$ be the full simplicial subcategory spanned by the fibrant and cofibrant objects. Then it is an immediate consequence of Quillen's axiom (SM7), that $\cM_{\ncf}$ is a locally fibrant simplicial category. Thus, via the coherent nerve construction, we obtain the $\infty$-category $N_\Delta(\cM_{\ncf}),$ the \textbf{underlying $\infty$-category} of the simplicial model category $\cM.$
\item As a more specific example let us consider $\sSet$ endowed with the usual Kan--Quillen model structure. This is a simplicial model category and with respect to this model structure we have $\sSet_{\ncf}=\cKan$, where~$\cKan$ is the full simplicial subcategory spanned by the Kan complexes. The \textbf{$\infty$-category}~$\cS$ \textbf{of spaces} is given by 
\[
\cS=N_\Delta(\cKan).
\]
Let us mention that this is only one model for the $\infty$-category of spaces and that there are others, for example the underlying $\infty$-category of topological spaces. For all purposes of~$\infty$-category theory it turns out that any $\infty$-category equivalent to~$\cS$ is equally good and that the precise model is irrelevant. All these $\infty$-categories satisfy the same universal property of being the `free cocomplete $\infty$-category on a single generator' (see \autoref{cor:spaces}).
\item Another class of examples is induced by additive categories. Given an additive category~$\cA$, the category~$\nCh(\cA)$ of chain complexes in~$\cA$ can be enriched over the category~$\nCh(\mathbb{Z})$ of chain complexes of abelian groups. The enrichment is set up in a way that if we take the 0-cycles of all the mapping complexes then we obtain the usual category of chain complexes, while cycles of positive dimensions give chain homotopies and higher chain homotopies.

Let us recall that the \emph{Dold--Kan correspondence} gives us an equivalence of categories $\nDK\colon\nCh(\mathbb{Z})\to\sAb,$ where~$\sAb$ is the category of simplicial abelian groups. An inverse to $\nDK$ is given by the \emph{normalized chain complex functor} which can be shown to be lax comonoidal with respect to the levelwise tensor product on $\sAb$ and the usual one on $\nCh(\mathbb{Z})$ (in fact, this lax comonoidal structure is induced by the Alexander--Whitney maps). It follows by abstract nonsense that~$\nDK$ carries canonically a lax \emph{monoidal} structure (see~\cite{kelly:doctrinal} for the abstract framework and \cite{schwede-shipley:equivalences} for precisely this context). Thus, we can apply~$\nDK$ to the morphism complexes in~$\nCh(\cA)$ in order to obtain the category~$\nDK_\ast(\nCh(\cA))$ enriched in simplicial abelian groups (see~\autoref{rmk:basechange}). Since simplicial abelian groups are Kan complexes this gives us a locally fibrant simplicial category, and we can define the \textbf{$\infty$-category} $\nCh(\cA)$ \textbf{of chain complexes} in~$\cA$ by 
\[
\nCh(\cA)=N_\Delta\big(\nDK_\ast(\nCh(\cA))\big).
\] 
A smaller model for these $\infty$-categories can be obtained by means of the \emph{differential-graded nerve construction} (see \cite{HA}).
\end{enumerate}
\end{egs}

We conclude this section with three more perspectives, one on conceptual remarks on mapping spaces in $\infty$-categories, one on enriched $\infty$-categories, and one on questions related to uniqueness issues for a theory of $(\infty,1)$-categories.

\begin{per}\label{per:mapping}
We saw in \autoref{rmk:principles} that associated to two objects in an $\infty$-category there is an entire space of morphisms $\Map_\cC^R(x,y)$ based on `right cones', and there is a dually defined space $\Map_\cC^L(x,y)$ given by `left cones'. Moreover, using prisms instead of cones, Lurie gave a third model $\Map_\cC(x,y)$ for such a space in \cite{HTT}, and shows that all these spaces are weakly equivalent Kan complexes. Using the Joyal model structure, this result can be put into perspective. 

To make this precise, let us recall that there is the theory of homotopy function complexes in an arbitrary model category~$\cM$, developed by Dwyer and Kan in the series of papers \cite{dwyer-kan:simplicial, dwyer-kan:calculating, dwyer-kan:function}. These function complexes $\Map_\cM(x,y)$ are simplicial sets which can be calculated by cosimplicial resolutions of~$x$, by simplicial resolutions of~$y,$ or by a combination of these resolutions (for a detailed account see e.g.\ \cite{hirschhorn:model}). Here, a cosimplicial resolution of~$x$ is a cofibrant replacement of the constant cosimplicial object~$\kappa_x\colon\Delta\to\cM$ in the Reedy model structure on~$c\cM$, and dually for simplicial resolutions of~$y$.

We denote by $\sSet_{K/}$ the category of simplicial sets under~$K\in\sSet$. This category comes with a canonical model structure such that the (co)fibrations and weak equivalences are created by $\sSet_{K/}\to\sSet$. In particular, there is the Joyal model structure on the category $\sSet_{\ast,\ast}=\sSet_{\del\Delta^1/}$ of simplicial sets with two distinguished vertices. Given an $\infty$-category~$\cC$ and two objects $x,y\in\cC$ we obtain the classifying map $(x,y)\colon\del\Delta^1\to\cC$, and hence an object $\cC_{x,y}\in\sSet_{\ast,\ast}.$ Moreover, the boundary inclusion $\del\Delta^1\to\Delta^1$ defines a further object $\Delta^1\in\sSet_{\ast,\ast}.$ Now, it can be shown that the various spaces of maps 
\[
\Map_{\cC}^L(x,y),\qquad\Map_{\cC}^R(x,y),\qquad\text{and}\qquad\Map_{\cC}(x,y)
\]
arise as different models for the homotopy function complexes $\Map_{\sSet_{\ast,\ast}}(\Delta^1,\cC_{x,y})$ with respect to the Joyal model structure on $\sSet_{\ast,\ast}$ (see \cite[Proposition 4.4]{dugger-spivak:mapping}). In fact, they are obtained by using different cosimplicial resolutions of $\Delta^1\in\sSet_{\ast,\ast}.$ Thus, using the fact that $\infty$-categories are fibrant objects in the Joyal model structure, this explains conceptually why these spaces have to be weekly equivalent Kan complexes.

All these models for mapping spaces are Kan complexes, so that it is `relatively easy' to calculate invariants like homotopy groups. As a drawback they do not admit a strict composition law. Using the simplicial category $C[\cC]$ associated to an $\infty$-category~$\cC$, we obtain still a different model for such spaces, namely the simplicial mapping spaces $\Map_{C[\cC]}(x,y)$ coming with a strictly associative and unital composition law. Besides the fact that these spaces might be very complicated to calculate there is the additional drawback that these spaces, in general, are \emph{not} Kan complexes. In particular, it is harder to calculate invariants. More precisely, Riehl established the following result \cite[Theorem 5.1]{riehl:scat-qcat}: Given an $\infty$-category~$\cC$ such that all mapping spaces in~$C[\cC]$ are Kan complexes, then~$\cC$ lies in the essential image of the usual nerve functor. More generally, Riehl shows that this already follows when all the $\Map_{C[\cC]}(x,y)$ are $\infty$-categories. 
\end{per}

\begin{per}\label{per:enriched}
Recall that a $2$-category is simply a category enriched over categories. Similarly, a strict $n$-category is just a category enriched over $(n\text{-}1)$-categories. This definition via iterated enrichments suggests that $(\infty,1)$-categories should be categories enriched over $(\infty,0)$-categories, i.e., Kan complexes by \autoref{cor:Kangroupoid}. And as we just saw, simplicial categories with the property that all mapping spaces are Kan complexes play a special role (see \autoref{thm:bergner}).

However, if we take the philosophy of $\infty$-categories serious, then the good notion of enrichment in that context should be that of a `weak enrichment'. We saw that for every pair of objects $x,y$ in an $\infty$-category~$\cC$, there is a space of morphisms $\Map_\cC(x,y)$. As opposed to the case of simplicial categories, for $\infty$-categories there is not a specified strictly associative and strictly unital composition law $\Map_\cC(y,z)\times\Map_\cC(x,y)\to\Map_\cC(x,z)$. In fact, the idea of having a weak enrichment should be expressed by the existence of composition laws which are only \emph{coherently associative}. Taking this perspective of $\mathbb{A}_\infty$-multiplications seriously, an abstract theory of enriched $\infty$-categories has recently been proposed by Gepner and Haugseng \cite{gepner-haugseng:enriched}. Their framework even allows one to study weak enrichments in not necessarily Cartesian contexts. For related rectification results see \cite{haugseng:rectification}.
\end{per}

\begin{per}\label{per:models}
\autoref{thm:JoyalBergner} is a comparison result between two approaches to a theory of $(\infty,1)$-categories. A similar result was also obtained by a combination of results due to Rezk, Joyal, Tierney, and Bergner. In \cite{rezk:model}, Rezk introduces Segal spaces as an alternative model for a theory of $(\infty,1)$-categories. These are certain bisimplicial sets and Rezk shows that there is an adapted model structure on the category of bisimplicial sets. Joyal and Tierney \cite{joyal-tierney:quasi-segal} show that the model category for Segal spaces and the Joyal model structure for $\infty$-categories are Quillen equivalent. Finally, Bergner \cite{bergner:three} shows that the model category for Segal spaces and the Bergner structure are Quillen equivalent through a zig-zag of Quillen equivalences. Thus, these results taken together also give a proof of the Quillen equivalence of the Joyal and the Bergner model structures through a zig-zag of Quillen equivalences. Nice introductions to some of these different models can be found in \cite{bergner:workshop} and in \cite{camerona:whirlwind}.  (We also want to mention that a part of this rich picture was generalized by Cisinski and Moerdijk in the series of papers \cite{cisinski-moerdijk:htpy-operads, cisinski-moerdijk:dendroidal-segal, cisinski-moerdijk:sop} to the framework of multicategories (or colored operads).)

Since by now there are many different approaches to a theory of $(\infty,1)$-categories one might wonder whether such a theory itself can be axiomatized. Motivated by previous work of Simpson \cite{simpson:n-categories}, a result in that direction was given by To{\"e}n in \cite{toen:axioms}. In that paper To{\"e}n axiomatizes a theory of $(\infty,1)$-categories as being a model category satisfying a list of seven axioms. He shows that any such model category is Quillen equivalent to the model category of complete Segal spaces. Moreover, the automorphisms of this theory are calculated to be $\mathbb{Z}/2\mathbb{Z},$ the non-trivial element being the passage to \emph{opposite $(\infty,1)$-categories}.

Similar calculations have also been carried out for $(\infty,n)$-categories. A first such calculation based on $\Upsilon_n$-spaces was given by Barwick--Schommer-Pries in \cite{barwick-schommer-pries:unicity}, while a calculation using $\Theta_n$-spaces was carried out by Ara, Guti\'errez, and the author in \cite{agg:autos}. In \cite{agg:autos} one can also find calculations of the automorphisms of the theory of $(\infty,1)$-operads and related theories.
\end{per}

\section{Categorical constructions with $\infty$-categories}
\label{sec:con}

The aim of this section is to extend some key constructions and notions from ordinary category theory to the world of $\infty$-categories (or more general simplicial sets) in a way that the following principles are satisfied.
\begin{enumerate}
\item The concepts are extensions of the ordinary concepts in that everything is compatible with the fully faithful nerve functor $N\colon\cCat\to\sSet.$
\item The notions are \emph{coherent} variants of the classical notions, i.e., $\infty$-category theory realizes some kind of homotopy coherent category theory.
\item The extensions are often defined for arbitrary simplicial sets, and when applied to $\infty$-categories we want these extensions to again give rise to $\infty$-categories. 
\item All concepts are \emph{invariant concepts}, i.e., an application of these constructions to equivalent input $\infty$-categories yields equivalent output $\infty$-categories.
\end{enumerate}
We are mostly interested in a robust theory of \emph{limits} and \emph{colimits} (see \S\ref{subsec:limits}) giving us an $\infty$-categorical version of the more well-known homotopy (co)limits in model categories, but this needs some preparation. The reader who is less inclined towards abstract categorical constructions is asked to consider this as a justification for the discussion of the constructions in \S\ref{subsec:join} and \S\ref{subsec:slice}. 

\subsection{Functors}
\label{subsec:functors}

Since $\infty$-categories are simply particular simplicial sets we can make the following definition.

\begin{defn}\label{defn:functors}
Let $K$ be a simplicial set and let~$\cC$ be an $\infty$-category. A \textbf{functor} $F\colon K\to\cC$ is a map of simplicial sets. Similarly, a \textbf{natural transformation} is a map $\Delta^1\times K\to\cC.$ More generally, the \textbf{space of functors} $\nFun(K,\cC)$ is 
\[
\nFun(K,\cC)_\bullet=\Map_{\sSet}(K,\cC)_\bullet=\hom_{\sSet}(\Delta^\bullet\times K,\cC)\in\sSet.
\]
\end{defn}
 
This extends the classical concept of a functor because $N\colon\cCat\to\sSet$ is fully faithful. More generally, there is the following refinement of this observation.

\begin{lem}
For categories $A,B$ there is a natural isomorphism of simplicial sets
\[
N(\nFun(A,B))\cong\nFun(NA,NB).
\]
\end{lem}
\begin{proof}
For $[n]\in\Delta$, there are the following natural bijections
\begin{align}
N(\nFun(A,B))_n&=\hom_{\cCat}([n],\nFun(A,B))\\
&\cong\hom_{\cCat}([n]\times A,B)\\
&\cong\hom_{\sSet}(N([n]\times A),NB)\\
&\cong\hom_{\sSet}(\Delta^n\times NA,NB)\\
&=\nFun(NA,NB)_n,
\end{align}
given by the exponential laws, the fully faithfulness of $N$, the fact that $N$ preserves products, and the isomorphism $N([n])\cong\Delta^n$. 
\end{proof}

This seemingly naive definition of a functor is actually the good one, as we want to indicate now (but see also \autoref{per:coherent}).

\begin{eg}
Let $A\in\cCat$ be an ordinary category and let $\cM$ be a locally fibrant simplicial category. By \autoref{cor:fibrantsimplicial} we know that $N_\Delta(\cM)$ is an $\infty$-category and we want to unravel a bit what it means to have a functor $F\colon NA\to N_\Delta(\cM).$ The behavior on 0-simplices and 1-simplices amounts to saying that associated to each arrow $x\to y$ in $A$ there is a morphism $Fx\to Fy$ in~$\cM.$ Moreover, given two composable arrows $x\stackrel{f}{\to} y \stackrel{g}{\to} z$ in $A$, i.e., a 2-simplex $\sigma\colon\Delta^2\to NA,$ we obtain a 2-simplex $F(\sigma)\colon\Delta^2\to N_\Delta(\cM).$ By \autoref{defn:cohnerve} this means that we are given a simplicial functor $C[\Delta^2]\to\cM$ which boils down to having a diagram in $\cM$ of the form
\[
\UseAllTwocells
\xymatrix@=1.5pc{
&Fy\ar[rd]^-{Fg}&\\
Fx \ar[ru]^-{Ff} \ar[rr]_-{F(g\circ f)} & & Fz.\lltwocell\omit{<2.5>}
}
\]
Thus, the functor $F$ preserves compositions up to specified homotopies. 

However, there is still by far more information encoded by $F$, namely all the higher simplices obtained from longer sequences of composable arrows in $A$. These encode the idea that $F$ is not only a `functor up to homotopy' but gives us a `functor up to \emph{coherent} homotopy', i.e., a \emph{homotopy coherent diagram}. For a precise statement see \cite{cordier:coherent} and also \autoref{per:coherent}. In order to at least give a partial justification for this claim let us consider a 3-simplex~$\tau\colon\Delta^3\to NA$, i.e., a chain consisting of three composable arrows $f,g,$ and $h.$ By the above we know that the four faces of the 3-simplex $F(\tau)\colon\Delta^3\to N_\Delta(\cM)$ altogether give us two different composite homotopies from $F(h\circ g\circ f)$ to $F(h)\circ F(g)\circ F(f)$ as in the boundary of
\begin{equation}
\vcenter{
\xymatrix@=1.5pc{
F(h\circ g\circ f)\ar[r]\ar[d]& F(h)\circ F(g\circ f)\ar[d]\\
F(h\circ g)\circ F(f)\ar[r]&F(h)\circ F(g)\circ F(f).
}
}
\label{eq:htpycoh}
\end{equation}
The 3-simplex $F(\tau)$ specifies one more homotopy $F(h\circ g\circ f)\to F(h)\circ F(g)\circ F(f)$ together with two 2-homotopies in $\cM$ relating this additional homotopy to the two compositions in \eqref{eq:htpycoh}. This follows from the fact that $\Map_{C[\Delta^3]}(0,3)$ is isomorphic to the product $\Delta^1\times\Delta^1$. (More generally, we have $\Map_{C[\Delta^n]}(i,j)\cong(\Delta^1)^{\times (j-i-1)}$ for $0\leq i< j\leq n$.)
\end{eg}

Given two ordinary categories $A$ and $B$, it is straightforward that the functors from $A$ to $B$ together with the natural transformations assemble to a category $\nFun(A,B).$ Moreover, if we have equivalences $A\simeq A'$ and $B\simeq B',$ then there is a canonical equivalence $\nFun(A,B)\simeq\nFun(A',B').$ Similar results also hold true in the world of $\infty$-categories, but there --- compared to the classical context --- a proof is already quite some work. One proof uses certain stability properties of the class of categorical equivalences and the so-called \emph{inner anodyne maps} (see for example \cite{joyal:barca}). Using these properties, one is able to deduce the following result.

\begin{prop}\label{prop:functor}
Let $\cC$ and $\cD$ be $\infty$-categories and let $K$ and $L$ be simplicial sets.
\begin{enumerate}
\item The simplicial set $\nFun(K,\cC)$ is an $\infty$-category.
\item If $\cC\to\cD$ is an equivalence of $\infty$-categories, then also the induced map $\nFun(K,\cC)\to \nFun(K,\cD)$ is an equivalence of $\infty$-categories.
\item If $K\to L$ is a categorical equivalence of simplicial sets, then the induced map $\nFun(L,\cC)\to\nFun(K,\cC)$ is an equivalence of $\infty$-categories.
\end{enumerate}
\end{prop}

Thus this proposition tells us that the formation of functor $\infty$-categories is an invariant notion --- as it should be the case for all categorical constructions in the world of $\infty$-categories. On a more conceptual side, this proposition is a consequence of the fact that with respect to the Joyal model structures the categorical product $\times\colon\sSet\times\sSet\to\sSet$ is a left Quillen functor of two variables. 

\begin{per}\label{per:coherent}
The theory of $\infty$-categories should really be thought of as \emph{homotopy coherent category theory}. In particular, the basic notion of a functor is to model the idea of having a \emph{homotopy coherent diagram}. Homotopy coherent category theory has quite some history and references include \cite{vogt:homotopy-limits, cordier:coherent, cordier-porter:vogt, cordier-porter:coherent}. Here we content ourselves by including a short detour only.

To begin with, whenever we have an adjunction $(L,R)\colon\cC\rightleftarrows\cD$ we obtain a comonad $C$ on $\cD$ with $C=LR\colon\cD\to\cD$ as underlying functor. The structure morphisms of the comonad are induced by the unit and counit of the adjunction, and the necessary relations follow directly from the triangular identities of an adjunction. These comonads can be used to form simplicial resolutions of objects in~$\cD$ (which play, for example, some role in relative homological algebra; see e.g.~\cite[\S8]{weibel:homological}). Slightly more precisely, given an object $d\in\cD$, the associated simplicial object $C_\ast(d)$ is given by $C_n(d)=C^{n+1}(d)$. (Of course, there is also a dual part of the story leading to monads and cosimplicial resolutions.) 

Now, the adjunction which is in the background of our context is the adjunction between graphs and categories given by the forgetful functor~$U$ sending a category to the underlying reflexive graph (`graph with identities') and the free category functor~$F$,
\begin{equation}
(F,U)\colon\cGraph\rightleftarrows\cCat.\label{eq:graph}
\end{equation}
Note that both the free category functor and the underlying graph functor preserve the set of objects. Thus, given a category~$A$, the resulting simplicial resolution $C_\ast A\in\nFun(\Delta\op,\cCat)$ also has a constant set of objects and hence defines an object $C_\ast A\in\sCat$, i.e., a simplicially enriched category as opposed to merely a simplicial object in categories. 

This allows us to make the following definition which goes back to Vogt~\cite{vogt:homotopy-limits}. Given a small category $A$ and a simplicial category $\cM\in\sCat,$ a \textbf{homotopy coherent diagram} in $\cM$ of shape $A$ is a simplicial functor $C_\ast A\to \cM.$ Thus, the set $\ncoh(A,\cM)$ of such homotopy coherent diagrams is defined by
\begin{equation}
\ncoh(A,\cM)=\hom_{\sCat}(C_\ast A,\cM).\label{eq:coherent}
\end{equation}

If as a special case we consider $A=[n]$, then the resulting simplicial category $C_\ast[n]$ is isomorphic to the `simplicial thickening'~$C[\Delta^n]$ introduced at the beginning of \S\ref{subsec:joyalbergner}. Thus, the coherent nerve $N_\Delta(\cM)$ of $\cM$ is precisely obtained by considering \emph{homotopy coherent chains of composable arrows} in the given $\cM$. 

More generally, Riehl~\cite[Theorem 6.7]{riehl:scat-qcat} showed that for every $A\in\cCat$ there is a natural isomorphism of simplicial categories $C[NA]\cong C_\ast(A).$ In particular, for $A\in\cCat$ and $\cM\in\sCat$ we have a bijection
\[
\ncoh(A,\cM)\cong\hom_{\sCat}(C[NA],\cM)\cong\hom_{\sSet}(NA,N_\Delta\cM).
\]
This shows that the seemingly naive \autoref{defn:functors} subsumes the concept of homotopy coherent diagrams, and justifies to think more generally of functors in the sense of that definition as homotopy coherent diagrams.

A final remark is in order: From a conceptual perspective one might argue that it is not very nice that we defined homotopy coherent diagrams in simplicial categories by means of \emph{specific simplicial resolutions} --- namely the simplicial comonad resolutions associated to \eqref{eq:graph}. In fact, the actual choice of resolution should not matter, and in some modern treatments like \cite[\S8]{goerss-jardine} homotopy coherent diagrams are defined by means of \emph{more general} simplicial resolutions.
\end{per}

\begin{rmk}
\begin{enumerate}
\item \autoref{prop:functor} reveals one of the technical advantages of $\infty$-categories over model categories: $\infty$-categories are stable under the formation of functor categories without any further assumption. In this respect, model categories are less well-behaved since one has to impose certain conditions on the model categories involved to obtain this stability property: for cofibrantly-generated model categories, associated diagram categories always admit the \emph{projective} model structure \cite[p.224]{hirschhorn:model}, whereas in the case of combinatorial model categories the \emph{projective} and the \emph{injective} structure both always exist on the diagram categories \cite[p.824]{HTT}. Note however that \autoref{prop:functor} is significantly more general since these results only tell us something about model structures on $\cM^A$, $A\in\cCat.$ For example, we never dared to ask for the existence of a `canonical' model category of functors between two given model categories.
\item A further technical advantage of $\infty$-categories over model categories is the following one. The `correct' notion of equivalence for model categories is the notion of \emph{Quillen equivalence}. Since, in general, a Quillen equivalence can not be inverted, the equivalence relation generated by this notion is quite complicated: frequently model categories are only Quillen equivalent through a zig-zag of Quillen equivalences pointing in different directions. The appropriate notion of equivalence for $\infty$-categories is the notion of \emph{categorical equivalence} (\autoref{defn:catequivalence}). Since $\infty$-categories are precisely the fibrant and cofibrant objects with respect to the Joyal model structure where categorical equivalences are the weak equivalences, it follows that a zig-zag of categorical equivalences can always be replaced by a single categorical equivalence.
\item A third technical advantage of $\infty$-categories was already mentioned, but will be repeated here for the sake of completeness. The notion of homotopy coherent diagrams is quite easily established in the world of $\infty$-categories since it is simply a map of simplicial sets with the domain given by the nerve of an ordinary category. There will be further advantages of this flavor, i.e., where `higher coherences' are easily encoded in the setting of $\infty$-categories. For example the notions of $\Aoo$- and $\Eoo$-algebras in monoidal and symmetric monoidal $\infty$-categories are conveniently introduced in this setting as specific sections of certain $\infty$-categorical versions of Grothendieck opfibrations, called coCartesian fibrations by Lurie. We will come back to this in \S\ref{sec:monoidal}.
\end{enumerate}
\end{rmk}

We conclude this subsection with a perspective on how to get an $\infty$-category of $\infty$-categories and also an underlying $\infty$-category for arbitrary model categories.

\begin{per}\label{per:inftyinfty}
We just introduced the notion of a functor between $\infty$-categories and hence obtain a category of $\infty$-categories. But to stick more seriously to the $\infty$-categorical framework, we would like to have an \emph{$\infty$-category~$\cCat_{\infty}$ of $\infty$-categories}. Since non-invertible natural transformations play an important role in classical category theory it would be even nicer to have an~$(\infty,2)$-category of $\infty$-categories but let us content ourselves with such an $\infty$-category.  Recall that the Joyal model structure is not simplicial so that we cannot apply  \autoref{cor:fibrantsimplicial} directly in order to get such a gadget. Luckily, there is the Quillen equivalent simplicial model category~$\sSet^+$ of so-called \emph{marked simplicial sets} which serves for this purpose (among many others), and which we want to describe very briefly. 

A \textbf{marked simplicial set} is a pair $(X,\mathcal{E}_X)$ consisting of a simplicial set~$X$ together with a subset $\mathcal{E}_X\subseteq X_1$ of edges containing the degenerate ones. A morphism of marked simplicial sets $(X,\mathcal{E}_X)\to(Y,\mathcal{E}_Y)$ is a morphism $f\colon X\to Y$ of simplicial sets such that $f(\mathcal{E}_X)\subseteq\mathcal{E}_Y$. This defines the category $\sSet^+$ of marked simplicial sets, a category playing a key role in the relative context \cite{HTT} (see also \autoref{per:Grothendieck} for further motivational remarks).

Degenerate edges, i.e., identity morphisms in an $\infty$-category are equivalences and it follows that every $\infty$-category~$\cC$ gives us a marked simplicial set~$\cC^\natural$ by marking the equivalences. It can be shown that there is a simplicial model structure on~$\sSet^+$ such that the fibrant and cofibrant objects are precisely the marked simplicial sets of the form~$\cC^\natural$ for some $\infty$-category~$\cC.$ Thus, the \textbf{$\infty$-category of $\infty$-categories}~$\cCat_\infty$ can be defined as the underlying $\infty$-category of~$\sSet^+,$
\[
\cCat_\infty=N_\Delta(\sSet^+_\ncf).
\]

A further use of these marked simplicial sets is that they allow us to define the \textbf{underlying $\infty$-category of an \emph{arbitrary} model category}. Given a model category $\cM$, by forgetting the cofibrations and the fibrations one obtains a category with weak equivalences $(\cM,W)$ (this is sometimes also called a \emph{relative category} \cite{barwick-kan:relative}). The \emph{ordinary} nerve construction then gives us a marked simplicial set $N(\cM,W)$. The underlying $\infty$-category of $\cM$ is defined to be a fibrant and cofibrant replacement of this marked simplicial set --- the idea being that this replacement amounts to inverting the marked edges in a universal way. It can be shown that the underlying $\infty$-category enjoys the expected universal property as the localization at~$W$; see \cite[\S1]{HA}.
\end{per}

\subsection{Join construction}
\label{subsec:join}

It seems to be fair to say that one main reason for studying the join construction in this context is that it allows us to define the slice construction (see \S\ref{subsec:slice}) which in turn is fundamental to the theory of (co)limits (see \S\ref{subsec:limits}). Before giving the definition of the join at the level of simplicial sets, we recall the classical situation in category theory.

Given categories $A$ and $B$, one can form a new category $A\star B$, the \textbf{join} of $A$ and~$B,$ as follows. The class of objects $\nobj(A\star B)$ is given by the disjoint union of $\nobj(A)$ and $\nobj(B).$ For the morphisms, there are the following four different cases
\[
\hom_{A\star B}(x,y)=\left\{ \begin{array}{c@{,\;}l}
				\hom_A(x,y)&\quad  x,y\in A,\\
				\hom _B(x,y)&\quad x,y\in B,\\
				\ast&\quad x\in A,y\in B,\\
				\emptyset &\quad x\in B,y\in A,
				\end{array} \right.
\]
and the composition is completely determined by requiring that $A$ and $B$ are full subcategories of $A\star B$ in the obvious way. Note that the construction is not symmetric in $A$ and $B$. To illustrate the join construction, we mention a few basic examples.

\begin{egs} 
\begin{enumerate}
\item If $A\in\cCat$ is arbitrary and if $B=\bbone$ is the terminal category, then $A^\rhd=A\star \bbone$ is the \textbf{right cone} or \textbf{cocone} on $A.$ It is obtained by adjoining a new terminal object $\infty$ to $A,$ and plays a central role in the study of colimits.
\item Dually, if $A=\bbone$ is terminal and $B\in\cCat$ is arbitrary, then $B^\lhd=\bbone\star B$ is the \textbf{left cone} or \textbf{cone} on $B.$ This category is obtained from $B$ by adjoining a new initial object $-\infty$ and is central to the theory of limits.
\item As a more specific example, let $A$ be the category occurring in the study of pushout diagrams,
\[
\xymatrix@=1.0pc{
(0,0)\ar[r]\ar[d]&(1,0)\\
(0,1).&  
}
\]
Then the cocone on $A$ is given by the commutative square $A^\rhd\cong [1]\times[1]$, which can of course be depicted as
\[
\xymatrix@=1.0pc{
(0,0)\ar[r]\ar[d]\ar@{}[rd]|{=}&(1,0)\ar[d]\\
(0,1)\ar[r]&(1,1).
}
\]
In particular, a diagram of this shape consists of four morphisms in the target category such that the two compositions agree. Similarly, if $B$ is the diagram occurring in the study of pullback diagrams, 
\[
\xymatrix@=1.0pc{
&(1,0)\ar[d]\\
(0,1)\ar[r]&(1,1),
}
\]
then also the cone $B^\lhd$ is the square.
\end{enumerate}
\end{egs}

The join construction can be extended to simplicial sets. There is a very conceptual approach to this construction as described by Joyal in \cite{joyal:barca} (basically as a Day convolution construction \cite{day:convolution} applied to the ordinal addition). In \cite{joyal:barca} one can also find many `elementary relations' about this join construction. Since we only want to rush through the theory of this notion, we instead give the following more direct `definition'.

\begin{defn}\label{def:join}
Let $K$ and $L$ be simplicial sets. The \textbf{join} $K\star L$ of $K$ and $L$ is the simplicial set defined by
\[
(K\star L)_n=K_n\cup L_n \cup \bigcup_{i+1+j=n}K_i\times L_j,\quad n\geq 0.
\]
\end{defn}

We leave it to the reader to define the structure maps of $K\star L$ and to check that this defines a functor $\star\colon\sSet\times\sSet\to\sSet$. It follows from those details that $K\star L$ comes with canonical inclusions $K\to K\star L$ and $L\to K\star L.$ The join operation for simplicial sets is in fact characterized by the following two properties.

\begin{prop}\label{prop:join}
\begin{enumerate}
\item The partial join functors $K\star (-)\colon\sSet\to\sSet_{K/}$ and $(-)\star L\colon\sSet\to\sSet_{L/}$ preserve colimits.
\item For the standard simplices we find $\Delta ^i \star \Delta ^j \cong \Delta ^{i+1+j}, i,j\geq 0,$ and these isomorphisms are compatible with the obvious inclusions of $\Delta^i$ and $\Delta^j.$
\end{enumerate}
\end{prop}

The following lemma is immediate and the proof is recommended as an exercise to those readers who just saw these notions for the first time. The solution to this exercise also suggests how to define the structure maps in \autoref{def:join}.

\begin{lem}\label{lem:nervejoin}
The nerve is compatible with the join constructions in that there is a natural isomorphism $N(A)\star N(B)\to N(A\star B), A,B\in\cCat.$
\end{lem}

To give some examples, we consider the (co)cone constructions and again pushout and pullback diagrams.

\begin{egs}\label{egs:join}
\begin{enumerate}
\item If $K\in\sSet$ is arbitrary and $L=\Delta^0,$ then $K^\rhd=K\star\Delta^0$ is the \textbf{right cone} or the \textbf{cocone} on $K$. Dually, if $L\in\sSet$ is arbitrary then $L^\lhd=\Delta^0\star L$ is the \textbf{left cone} or \textbf{cone} on $L.$
\item Let $K=\Lambda^2_0$. Using the description of this horn as a pushout, it follows immediately from \autoref{prop:join} that the cocone $(\Lambda^2_0)^\rhd$ is isomorphic to the square $\square=\Delta^1\times\Delta^1,$ i.e., we have the picture
\[
\xymatrix@=1.5pc{
(0,0)\ar[r]\drtwocell<\omit>{<2>}\ar[d]\ar[dr]&(1,0)\ar[d]\\
(0,1)\ar[r]&(1,1).\ultwocell<\omit>{<2>}
}
\]
Note that if $\cM$ is a simplicial category then a diagram of this shape in~$N_\Delta(\cM)$ consists of five morphisms and two homotopies as depicted by the diagram. In particular, in general, we do not have a commutative square but only a coherent version thereof. This example can be dualized and we obtain a similar isomorphism $(\Lambda^2_2)^\lhd\cong\Delta^1\times\Delta^1$.
\end{enumerate}
\end{egs}

A careful analysis of the join construction allows one to establish the following important properties.

\begin{prop}
\begin{enumerate}
\item If $\cC$ and $\cD$ are $\infty$-categories, then the join $\cC\star\cD$ is again an $\infty$-category.
\item If $F\colon\cC\to\cC'$ and $G\colon\cD\to\cD'$ are equivalences of $\infty$-categories, then also the induced map $F\star G\colon\cC\star\cD\to\cC'\star\cD'$ is an equivalence.
\end{enumerate}
\end{prop}

\subsection{Slice construction}
\label{subsec:slice}

We again begin by recalling the more classical situation of ordinary category theory. Given a category $B$ and an object $x\in B$, one can form the \textbf{overcategory} $B_{/x}$ where objects are morphisms $x'\to x$ in $B$ with target $x$. Given two such objects $x'\to x$, $x''\to x$, a morphism $(x'\to x)\to(x''\to x)$ in $B_{/x}$ is simply a morphism $x'\to x''$ in $B$ making the obvious triangle
\begin{equation}
\vcenter{
\xymatrix@=1.5pc{
x' \ar[rd] \ar[rr] \ar@{}[rrd]|{=}& & x''\ar[ld]\\
&x&  
}\label{eq:overcat}
}
\end{equation}
commute. 

One generalization of this notion is obtained by replacing the object $x\colon\bbone\to B$ by a more general diagram in~$B$. More precisely, if we are given a functor $p\colon A\to B,$ then we can form the \textbf{slice category} $B_{/p}$ of \emph{objects over~$p$} or \emph{cones on~$p$}. Thus, an object in $B_{/p}$ is simply a cone on $p$, i.e., an object $b\in B$ together with a natural transformation from the constant functor on~$b$ to the given~$p$. Morphisms in this category are simply morphisms in~$B$ which are compatible with the natural transformations. Using the join construction one can see that the slice construction satisfies a universal property: for any category $C,$ there is a natural bijection
\[
\nFun(C,B_{/p})\cong\nFun_p(C\star A ,B),
\]
where the right-hand side denotes all functors $C\star A\to B$ making the following triangle commute
\[
\xymatrix@=1.5pc{
&A\ar[dl]\ar[dr]^-p&\\
C\ast A\ar[rr]&&B.
}
\]
To emphasize a bit more that the right-hand side in the above universal property takes certain structure maps into account, we rewrite this as
\[
\hom_{\cCat}(C,B_{/p})\cong\hom_{\cCat_{A/}}(A\to C\star A, A\stackrel{p}{\to} B),\quad C\in\cCat.
\] 
Of course this sounds like an unnecessarily complicated reformulation of something very simple. But the point of this reformulation is that it gives us an idea on how to extend these notions to $\infty$-categories --- as it was done by Joyal in \cite{joyal:quasi-kan}.

\begin{prop}
Let $p \colon L\to\cC$ be a map of simplicial sets with $\cC$ an $\infty$-category. There is an \emph{$\infty$-category $\cC_{/p}$} characterized by the following universal property: For every simplicial set $K,$ there is a natural bijection
\[
\hom_{\sSet}(K,\cC_{/p})\cong\hom_{\sSet_{L/}}(L\to K\star L,L\to\cC).
\]
The $\infty$-category $\cC_{/p}$ is the \textbf{$\infty$-category of cones on~$p$}.
\end{prop}

To check that there is such a \emph{simplicial set} $\cC_{/p}$, one can use the universal property as a definition. More precisely, the special cases of the standard simplices $K=\Delta^n$ give us by the Yoneda lemma a description of the $n$-simplices of $\cC_{/p}$ as 
\begin{equation}
(\cC_{/p})_n\cong\hom_{\sSet_{L/}}(L\to\Delta^n\ast L,L\to \cC).\label{eq:slice}
\end{equation}
To show that one actually obtains an $\infty$-\emph{category} requires more work and will not be done here. 

\begin{rmk}
\begin{enumerate}
\item One can show that the slice construction is an \emph{invariant} notion in a certain precise sense.
\item In both the classical and the $\infty$-categorical situation the constructions can be dualized. For example, there is the \textbf{$\infty$-category $\cC_{p/}$ of cocones on $p$}. Any of the $\infty$-categories $\cC_{p/}$ and $\cC_{/p}$ is also referred to as a \textbf{slice $\infty$-category}.
\end{enumerate}
\end{rmk}
 
There is the following result about the compatibility of the nerve $N\colon\cCat\to\sSet$ with slice constructions.

\begin{lem}
If $p\colon A\to B$ is a functor, then there is a natural isomorphism of simplicial sets
\[
N(B_{/p})\cong N(B)_{/N(p)}.
\]
\end{lem}
\begin{proof}
In simplicial dimension~$n$ we have the following chain of natural bijections
\begin{align}
N(B_{/p})_n&=\hom_{\cCat}([n],B_{/p})\\
&\cong\hom_{\cCat_{A/}}(A\to[n]\star A,A\to B)\\
&\cong\hom_{\sSet_{N(A)/}}(N(A)\to N([n]\star A),N(A)\to N(B))\\
&\cong\hom_{\sSet_{N(A)/}}(N(A)\to \Delta^n\star N(A),N(A)\to N(B))\\
&\cong\hom_{\sSet}(\Delta^n,N(B)_{/N(p)}),
\end{align}
where in step three we used \autoref{lem:nervejoin} and it is obvious which canonical isomorphisms were used in the remaining steps.
\end{proof}

\begin{eg}\label{eg:undercategory}
Let $\cC$ be an $\infty$-category and let $x\in\cC$ be an object, classified by the map $\kappa_x\colon\Delta^0\to\cC.$ Then the $\infty$-category $\cC_{/\kappa_x}$ is called the \textbf{$\infty$-category of objects over $x$}, and is simply denoted by $\cC_{/x}$. Dually, the $\infty$-category $\cC_{\kappa_x/}$ is called the \textbf{$\infty$-category of objects under} $x,$ and is denoted by $\cC_{x/}.$

Let us spell out some of the details about $\cC_{/x}$ in order to convince ourselves that this is the expected coherent version of the corresponding classical construction. By \eqref{eq:slice}, an object in $\cC_{/x}$ is a map $f\colon\Delta^1\cong\Delta^0\star\Delta^0\to\cC$ such that $d_0(f)=x$, i.e., we are given a map in~$\cC$ with target~$x.$ Similarly, \eqref{eq:slice} tells us that a map in $\cC_{/x}$ is given by a 2-simplex $\sigma\colon\Delta^2\cong\Delta^1\star\Delta^0\to\cC$ such that $\sigma(2)=x.$ In the case that $\cC$ is the coherent nerve of a (locally fibrant) simplicial category, the picture of~$\sigma$ is
\[
\UseAllTwocells
\xymatrix@=1.5pc{
x' \ar[rd] \ar[rr]  & & x''\ar[ld]\\
&x,&\lltwocell\omit{<3.5>}
}
\]
indicating that we obtained the coherent version of \eqref{eq:overcat} we were aiming for.
\end{eg}

\subsection{Final and initial objects}
\label{subsec:final}

We now come to the $\infty$-categorical variant of final and initial objects. In classical category theory, final objects are characterized by the property that for all objects there is a unique morphism to the final object. In these notes, we take a slightly different approach to the $\infty$-categorical generalization (but see \autoref{prop:final}). The following definition is not the original one of Joyal \cite{joyal:quasi-kan}, but it is an equivalent one as shown in \cite{joyal:quasi-kan}.

\begin{defn}
An object $x\in\cC$ of an $\infty$-category~$\cC$ is a \textbf{final object}, if the canonical map $\cC_{/x}\to\cC$ is an acyclic fibration of simplicial sets.
\end{defn}

Thus, an object is final if `forgetting that we lived above it does not result in a loss of information'. This notion can be reformulated in the following more usual way.

\begin{prop}\label{prop:final}
The following are equivalent for an object $x$ of an $\infty$-category~$\cC.$
\begin{enumerate}
\item The object $x$ is final.
\item The mapping spaces $\Map_{\cC}(x',x)$ are acyclic Kan complexes for all $x'\in\cC.$
\item Every simplicial sphere $\alpha\colon\del\Delta^n\to\cC$ such that $\alpha(n)=x$ can be filled to an entire $n$-simplex $\Delta^n\to\cC.$
\end{enumerate}  
\end{prop}

\autoref{cor:choices} is an immediate consequence of this proposition, but let us first make precise the notion of a \emph{full} subcategory of an $\infty$-category. There is a more general notion of subcategories of an $\infty$-category but since we will not need this additional generality we stick to full subcategories. Given objects $\cD_0\subseteq\cC_0$ in an $\infty$-category~$\cC$, let $\cD\subseteq \cC$ be the simplicial subset consisting precisely of those simplices $\Delta^n\to\cC$ which have the property that all vertices belong to $\cD_0$. Then it is immediate that $\cD$ is again an $\infty$-category, the \textbf{full subcategory of $\cC$ spanned by $\cD_0.$} Obviously, $\cD$ comes with an inclusion $\cD\to\cC.$

\begin{cor}\label{cor:choices}
Let $\cC$ be an $\infty$-category and let $\cD\subseteq \cC$ be the full subcategory spanned by the final objects of $\cC.$ Then $\cD$ is empty or a contractible Kan complex.
\end{cor}
\begin{proof}
This is immediate from \autoref{prop:final}.
\end{proof}

\begin{rmk}\label{rmk:unique}
The conclusion of \autoref{cor:choices} is typical for uniqueness statements in $\infty$-category theory: It states that a space `parametrizing universal objects' is empty or a contractible Kan complex. In classical category theory, if universal objects exist then they are unique up to unique isomorphism. This can be reformulated by saying that if universal objects exist then the category of such is a \emph{contractible groupoid}: the existence of comparison maps shows that it is a connected category, while the uniqueness of these maps tells us that it is a groupoid with no nontrivial endomorphisms. Thus, the possibly non-trivial $\pi_0$ and $\pi_1$ both vanish. Now, with \autoref{cor:Kangroupoid} in mind we see that the two uniqueness statements are morally very similar. 

One example of such a uniqueness statement was already obtained in the discussion of \autoref{thm:uniquecompo}, namely the space of compositions of two composable arrows is a contractible Kan complex. With \autoref{cor:choices} we have a further example of such a statement. In \S\ref{subsec:limits} we will define (co)limits as universal objects of certain slice categories and will hence again obtain such uniqueness statements. Let us mention that uniqueness results of this kind are also ubiquitous in the theory of model categories. Compare for example to \cite{hirschhorn:model} where many categories of choices (for example, cofibrant replacements) are shown to have contractible nerves.
\end{rmk}

We conclude by a statement about the compatibility with the nerve construction, but leave the proof to the reader.

\begin{lem}
Let $A$ be a category. An object $a\in A$ is final if and only if $a\in N(A)$ is final.
\end{lem}

\subsection{Limits and colimits}
\label{subsec:limits}

Having the basic categorical notions at our disposal, we are ready to talk about (co)limits in the framework of $\infty$-categories. Recall that the colimit of an ordinary functor $p\colon A\to B$ consists of an object $\colim_A p$ in $B$ together with a universal cocone. Said differently, such a pair is equivalently specified by an initial object in the category $B_{p/}$ of cocones on~$p$. This definition can now readily be extended to $\infty$-categories (see \cite{joyal:quasi-kan}).

\begin{defn}\label{defn:limits}
Let $K$ be a simplicial set and let $\cC$ be an $\infty$-category. A \textbf{colimit} of a diagram $p \colon K\to\cC$ is an initial object in $\cC_{p/}.$ An $\infty$-category is \textbf{cocomplete} if it admits colimits of all small diagrams. Dually, we define \textbf{limits} of diagrams and \textbf{complete} $\infty$-categories. 
\end{defn}

It is an immediate consequence of \autoref{cor:choices} that for every $p\colon K\to\cC$ the full subcategory $\cD\subseteq\cC _{p/}$ spanned by the colimits of $p$ is empty or a contractible Kan complex, i.e., that if a colimit exists, then it is unique up to a contractible choice.

\begin{rmk}
\begin{enumerate}
\item The nerve functor $N\colon\cCat\to\sSet$ is compatible with the notion of (co)limits.
\item The notions of limits and colimits are invariant notions. 
\end{enumerate}
\end{rmk}

One justification for \autoref{defn:limits} is that it extends the classical theory of (co)limits. Although this is certainly a convenient aspect of the notion, there is an additional justification which is more relevant for our purposes (as discussed by Lurie in \cite[Theorem 4.2.4.1]{HTT}). Namely, there is a precise meaning in that these notions of (co)limits coincides with the notion of \emph{homotopy (co)limits} in simplicial categories. Joyal was fully aware of this in \cite{joyal:quasi-kan}. References to various aspects of the theory of homotopy (co)limits in other contexts include \cite{bousfield-kan:htpy-limits}, \cite{vogt:homotopy-limits}, \cite{hirschhorn:model}, \cite{dhks:homotopy}, \cite{shulman:homotopy-limits}, and \cite{dugger:primer}. We follow Joyal and Lurie and simply speak of \emph{(co)limits} instead of \emph{homotopy (co)limits}, leading to a less heavy terminology.

As a special case we have the following definition; see again \autoref{egs:join}.

\begin{defn}\label{defn:pull-push}
Let $\cC$ be an $\infty$-category and let $q\colon\square\to\cC$ be a square.
\begin{enumerate}
\item The square~$q$ is a \textbf{pushout} if $q\colon(\Lambda^2_0)^\rhd\to\cC$ is a colimiting cocone.
\item The square~$q$ is a \textbf{pullback} if $q\colon(\Lambda^2_2)^\lhd\to\cC$ is a limiting cone.
\end{enumerate}
\end{defn}

Having a good definition of (co)limits available, one would now like to establish batteries of techniques to play with these notions. Many classical facts from the calculus of ordinary (co)limits can be extended to the context of $\infty$-categories, although the proofs of even fundamental statements are significantly harder in this framework. Nevertheless, such an extension was achieved to an impressive extent by Lurie in \cite{HTT}.

Instead of pursuing this any further, we conclude this section by two more perspectives, one dealing with the fact that the underlying $\infty$-category of a (simplicial) combinatorial model category is (co)complete, and the other one with \emph{derivators} as a convenient framework encoding the calculus of Kan extensions available in (co)complete $\infty$-categories.

\begin{per}
A model category is by definition a category admitting (at least finite) limits and colimits. At least equally interesting, model categories allow for a calculus of homotopy (co)limits, and it is thus natural to wonder if the underlying $\infty$-category of a model category is (co)complete. Morally, this should be true for \emph{any} model category (see \autoref{per:derivator} for some justification of this statement), but until now proofs were only written up in the case of \emph{combinatorial} model categories (see \S\ref{subsec:combmod} for some basics on combinatorial model categories although no details are needed for the understanding of this perspective). The main steps in the case of \emph{simplicial combinatorial} model categories are as follows. 
\begin{enumerate}
\item First, one introduces \emph{cofinal morphisms} of simplicial sets \cite[p.224]{HTT}. In the case that the target of the morphism is an $\infty$-category, one can give a definition which is the `expected analogue' of the classical variant. For this, let $v\colon K\to\cC$ be a diagram and let $\cC$ be an $\infty$-category. Then for every object $x\in\cC,$ consider the following pullback diagram 
\[
\xymatrix@=1.5pc{
\cC_{x/}\underset{\cC}{\times}K\pullbackcorner\ar[r]\ar[d]&\cC_{x/}\ar[d]\\
K\ar[r]&\cC,
}
\]
in which $\cC_{x/}\to\cC$ sends an object $x\to x'$ to its target. The morphism~$v$ is \textbf{cofinal} if the above pullback is \emph{contractible} for every $x\in\cC.$ Recall from classical category theory that the cofinality of a functor can be characterized by the \emph{connectivity} of similar undercategories. The case of morphisms of arbitrary simplicial sets is more subtle and was first treated by Joyal \cite{joyal:I-II}.
\item A justification for the terminology of a cofinal morphism $v\colon K\to K'$ is given by the following result: if an $\infty$-category admits colimits of shape $K$ then it also admits colimits of shape~$K'$ (see \cite[\S4.1]{HTT} for further and more precise statements).
\item Lurie shows next that for every simplicial set $K$ there is a category $A$ and a cofinal functor $N(A) \to K$ \cite[p.255]{HTT}. Thus, concerning existence questions for colimits, one can restrict attention to diagrams defined on nerves of categories.
\item A further ingredient is that homotopy coherent diagrams (\autoref{per:coherent}) in the $\infty$-category underlying a simplicial model category can be rigidified, and that the two possibly different notions of (homotopy) colimits coincide.
\end{enumerate}
All these results can be combined to establish the (co)completeness of the $\infty$-category underlying a simplicial combinatorial model category. Obviously, this lies outside the scope of these notes, but is treated in detail in \cite[\S4.1-4.2]{HTT}. Lurie also establishes a similar result for combinatorial model categories in \cite[\S1]{HA}.
\end{per}

\begin{per}\label{per:derivator}
Having a notion of (co)limits at our disposal, one might wonder if there is also a way to extend the theory of \emph{Kan extension} to the context of $\infty$-categories and hence to obtain a `relative theory of (co)limits'. Joyal and Lurie showed this to be possible, and many of the key formal properties of Kan extension are still true in this $\infty$-categorical framework.

We begin by recalling the concept of Kan extensions from ordinary category theory. Let us focus on left Kan extensions but everything can be dualized to right Kan extensions. Given a category~$\cC$ and a functor $u\colon A\to B$ between small categories, then we can consider the precomposition functor $u^\ast\colon\cC^B\to\cC^A.$ A \textbf{left Kan extension} $\nLan_u\colon\cC^A\to\cC^B$ along~$u$ is a left adjoint to $u^\ast$,
\[
(\nLan_u,u^\ast)\colon\cC^A\rightleftarrows\cC^B.
\]
As a special case, associated to a small category $A$ there is the unique functor $\pi_A\colon A\to\bbone$ to the terminal category, and in this case left Kan extensions along $\pi_A$ reduce to ordinary colimits of shape~$A$. One can show, more generally, that arbitrary left Kan extensions always exist if the target category $\cC$ is cocomplete. 

This sounds like a fairly abstract concept but many important constructions are special cases of (combinations of) \textbf{homotopy Kan extensions}, i.e., suitably derived versions of ordinary Kan extensions. This is for example the case for cones, fibers, suspensions, loops, Baratt--Puppe sequences, pushout products, homotopy tensor products of functors, homotopy orbits, homotopy fixed points, homotopy (co)limits, and spectrifications of prespectrum objects. 

A key feature of the calculus of (homotopy) Kan extensions is that they can be calculated \emph{pointwise}. In the case of ordinary left Kan extensions in cocomplete categories this means the following. Given a diagram $X\colon A\to\cC$ and an object $b\in B,$ then a certain canonical map
\[
\colim_{(A/b)}X\circ p\to\nLan_u(X)_b
\]
is an isomorphism. As mentioned above, the theory of Kan extensions can be extended to the framework of $\infty$-categories, and it turns out that Kan extensions for $\infty$-categories can also be calculated pointwise as well as homotopy Kan extensions in model categories.

There is an alternative approach to abstract homotopy theory which focuses on these Kan extensions and which can be motivated as follows. Recall that, despites their usefulness in many contexts, ordinary homotopy categories or derived categories are categorically rather badly behaved. In typical situations the only categorical (co)limits which exist are final or initial objects and (finite) products and coproducts. Moreover, the powerful calculus of homotopy Kan extensions is not visible to the homotopy category \emph{alone}. As a punchline, the transition from a model category or complete and cocomplete $\infty$-category to the homotopy category results in a serious loss of information.

The theory of \textbf{derivators} going back to Heller~\cite{heller:htpythies}, Grothendieck~\cite{grothendieck:derivators}, Keller~\cite{keller:epivalent}, Franke~\cite{franke:adams}, and probably others addresses this issue (see \cite{grothendieck:derivators} for many references): a derivator can be considered as a minimal extension of an ordinary homotopy category to a framework in which this calculus is available. Roughly speaking, a derivator \D consists of homotopy categories $\D(A)$ of $A$-shaped coherent diagrams together with restriction maps $u^\ast\colon\D(B)\to\D(A).$ The key defining \emph{property} of a derivator is that each $u^\ast$ admits a left adjoint and a right adjoint, and that these Kan extensions can be calculated pointwise. In particular, in a derivator homotopy (co)limits can be characterized by \emph{ordinary categorical universal properties}, and the theory is hence accessible by elementary categorical methods.

Cisinski showed that \emph{every} model category has an underlying derivator (see \cite{cisinski:direct} for the general case and \cite{groth:ptstab} for an easy proof in the case of combinatorial model categories) --- and this is one reason why we expect that every model category has an underlying (co)complete $\infty$-category. Similarly, in \cite{gps:stable} there is a sketch proof that every complete and cocomplete $\infty$-category has an underlying derivator. Thus, derivators can be considered as a convenient, purely categorical framework for a formal study of the calculus of homotopy Kan extensions which is available in complete and cocomplete $\infty$-categories and model categories and hence in various contexts arising in algebra, geometry, and topology. Derivators were recently studied by Ayoub \cite{ayoub:I,ayoub:II}, Cisinski \cite{cisinski:derived-kan}, Maltsiniotis \cite{maltsiniotis:k-theory,maltsiniotis:htpy-exact}, Tabuada \cite{tabuada:universal-invariants} and others. See the webpage \cite{grothendieck:derivators} for more references.
\end{per}

\section{Presentable $\infty$-categories and the relation to model categories}
\label{sec:presentable}

In this section we give a short introduction to \emph{presentable $\infty$-categories}, a class of $\infty$-categories having very good formal properties. For example, the adjoint functor theorem of Freyd takes in this context the form that a functor is a left adjoint if and only if it preserves colimits, and there is also a result for right adjoint functors. Besides these good formal properties, it turns out that a good deal of mathematics is encoded by presentable $\infty$-categories. Many typical $\infty$-categories showing up in nature, e.g., in algebra, topology, or (derived) algebraic geometry are presentable. 

The main reason for us to include a short discussion of presentable $\infty$-categories is that they play an essential role in Lurie's treatment of the \emph{smash product} on the $\infty$-category of spectra (see \S\ref{sec:stable}). Moreover, they allow us to indicate a more precise relation between $\infty$-categories and model categories (see \S\ref{subsec:pres} and, in particular, \autoref{thm:comparison}). Nevertheless, the reader can skip this section on a first reading since this will only affect the understanding of some of the statements in \S\ref{sec:stable}.

The theory of presentable $\infty$-categories has (at least) two precursors; \emph{locally presentable categories} in the classical context as well as \emph{combinatorial model categories} in the homotopical framework. We emphasize that, in all these three cases, there are two main ideas in the background, which help organize the material.

\begin{enumerate}
\item The idea of passing from a small category to a cocomplete category in a universal way is realized by categories of presheaves in the classical case, and by model categories and $\infty$-categories of simplicial presheaves in the remaining two. One can think of this passage as some kind of \emph{free generation}.
\item All locally presentable categories can be obtained as certain localizations of presheaf categories, and similarly in the homotopical contexts using suitable notions of (Bousfield) localizations. The passage to (Bousfield) localizations can be thought of as a way of \emph{imposing relations}. 
\end{enumerate}

\subsection{Locally presentable categories}
\label{subsec:locprescat}

In this subsection we give a short review of the theory of locally presentable categories. For more details we refer to \cite{gabriel-ulmer:lokal, makkai-pare, adamek-rosicky} or \cite{borceux:1, borceux:2}. 

To begin with, let us again take up the idea that categories of presheaves (with values in sets) can be considered as cocompletions. If $A$ is a small category, then we denote by $y\colon A\to \nFun(A\op,\cSet)$ the Yoneda embedding which sends $a\in A$ to the represented presheaf $\hom_A(-,a)$. Associated to each presheaf $X\colon A\op\to\cSet$ there is the comma category $(y/X)$. Objects of $(y/X)$ are pairs $(a,\alpha)$ consisting of an object $a\in A$ and a natural transformation $\alpha\colon y(a)\to X$, and a morphism $f\colon (a,\alpha)\to (a',\alpha')$ is a morphism $f\colon a\to a'$ in $A$ such that the diagram
\[
\xymatrix@=1.5pc{
y(a)\ar[rr]^-{y(f)}\ar[dr]_-\alpha&&y(a')\ar[dl]^-{\alpha'}\\
&X&
}
\]
commutes. (Note that the Yoneda lemma implies that this category is isomorphic to the \emph{category of elements} of $X$.) The category $(y/X)$ comes with a projection functor $p\colon (y/X)\to A$ sending a pair $(a,\alpha)$ to $a$, and we can hence consider the composition
\[
(y/X)\stackrel{p}{\to} A\stackrel{y}{\to} \nFun(A\op,\cSet).
\]
Using these diagrams one can make precise that presheaves on small categories are canonically colimits of representable ones (see for example \cite[p.76]{MacLane}).

\begin{prop}\label{prop:Yonedadense}
Let $A$ be a small category and let $X\colon A\op\to\cSet$ be a set-valued presheaf. There is a canonical isomorphism
\[
\colim_{(y/X)}y\circ p\cong X.
\]
\end{prop}

Recall that if $\cC$ is a cocomplete category then also the associated functor categories $\nFun(A,\cC)$ are cocomplete (and colimits are defined pointwise). In particular, given a small category $A$, then the presheaf category $\nFun(A\op,\cSet)$ is cocomplete (actually complete and cocomplete). The proposition tells us in a certain vague sense that $\nFun(A\op,\cSet)$ should be the universal approximation of $A$ by a cocomplete category: since every presheaf is canonically a colimit of representable presheaves we did not add too much by passing to this cocomplete category. Thus, we think of the passage to presheaves as a \textbf{cocompletion} or, more informally, as a \emph{free generation}. In order to make a more precise statement, let $\nFunL(-,-)$ denote the category of colimit-preserving functors. The following result seems to go back to Ulmer \cite[Rmk.~2.29]{ulmer:dense}.

\begin{thm}\label{thm:yoneda}
Let $A$ be a small category and let $\cC$ be a cocomplete category. The restriction along the Yoneda embedding $y\colon A\to\nFun(A\op,\cSet)$ induces an equivalence of categories 
\[
y^\ast\colon\nFunL(\nFun(A\op,\cSet),\cC)\stackrel{\sim}{\to}\nFun(A,\cC).
\]
\end{thm} 

To motivate the notion of a locally presentable category we recall Freyd's Adjoint Functor Theorem. Let $F\colon\cC\to\cD$ be a functor between cocomplete categories. It is obvious that if $F$ is a left adjoint, then $F$ preserves all colimits, but, in general, the converse is not true. There is the celebrated \textbf{Adjoint Functor Theorem} of Freyd (see \cite[pp.~84-86]{freyd:abelian} or \cite[p.~121]{MacLane}) which gives \emph{necessary and sufficient} conditions for the existence of a right adjoint. Slightly more precisely, the following are equivalent for a functor $F\colon\cC\to\cD$ between cocomplete categories.
\begin{enumerate}
\item The functor $F$ is a left adjoint.
\item The functor $F$ preserves colimits and the \emph{solution-set-condition} is satisfied.
\end{enumerate}
Without going into detail, let us only mention that the solution-set-condition states that a certain \emph{class} of arrows turns out to be small enough to actually form a \emph{set}. Hence, one can imagine that this condition is automatically satisfied if we impose some `smallness conditions' on the categories. Recall that a \emph{small, cocomplete} category is necessarily a poset (\cite[p.114]{MacLane}). Thus, in order to not rule out interesting examples, it is essential that these `smallness conditions' are chosen in a smart way.

\begin{defn}
A category is \textbf{locally presentable} if it is cocomplete and accessible.
\end{defn}

The accessibility assumption in this definition is the smallness assumption alluded to above. The idea is that a category $\cC$ is accessible if it admits certain filtered colimits and if it is formally determined by some small subcategory $\cD$ consisting of small objects. To make this more precise, let us consider a regular cardinal number~$\kappa$. A category $\cC$ is \textbf{$\kappa$-accessible} if $\cC$ admits $\kappa$-filtered colimits and if there is a small subcategory $\cD\subseteq\cC$ such that the following two conditions are satisfied.
\begin{enumerate}
\item Every object of $\cC$ can be canonically written as a $\kappa$-filtered colimit of objects in~$\cD$ (more technically, the inclusion $\cD\to\cC$ is \emph{dense}). 
\item The set-valued functors $\hom_{\cC}(d,-)\colon\cC\to\cSet,\;d\in\cD,$ preserve $\kappa$-filtered colimits --- expressing the idea that objects in~$\cD$ are small. 
\end{enumerate}
We say that a category is \textbf{accessible} if it is $\kappa$-accessible for some~$\kappa$. Similarly, a functor $F\colon\cC\to\cD$ is \textbf{$\kappa$-accessible} if $\cC$ and $\cD$ admit $\kappa$-filtered colimits and if they are preserved by $F$. Finally, a functor is \textbf{accessible} if it is $\kappa$-accessible for some~$\kappa$. 

For more details on the rich theory of accessible and locally presentable categories we again refer to \cite{gabriel-ulmer:lokal, makkai-pare, adamek-rosicky} or \cite{borceux:1, borceux:2}. But let us emphasize that every locally presentable category is also a \emph{complete} category. To indicate the ubiquity of locally presentable categories, we include the following list of examples.

\begin{egs}\label{egs:locpres}
\begin{enumerate}
\item The category $\cSet$ of sets is locally presentable.
\item If $A$ is a small category, then the category $\nFun(A\op,\cSet)$ of presheaves on $A$ is locally presentable. In particular, the category $\sSet$ of simplicial sets is locally presentable.
\item For a ring $R$ the categories $\Mod(R)$ of $R$-modules and $\mathrm{Ch}(R)$ of chain complexes over $R$ are locally presentable. 
\item Recall that an abelian category with exact filtered colimits is \textbf{Grothendieck abelian} (see \cite{grothendieck:tohoku} and \cite[\S14]{faith}) if it admits a generator. It can be shown that an abelian category with exact filtered colimits is Grothendieck if and only if it is locally presentable. Important examples are given by categories of quasi-coherent $\mathcal{O}_X$-modules on any scheme~$X$.
\item Categories of modules of (multi-sorted) algebraic theories are locally presentable \cite[\S6]{adamek-rosicky-vitale}.
\item If $T\colon\cC\to\cC$ is an accessible monad on a locally presentable category, then the category of $T$-algebras is locally presentable~\cite{adamek-rosicky}.
\item Every Grothendieck topos~\cite{SGA41,borceux:3,maclane-moerdijk} is locally presentable.
\item The category $\cCat$ of small categories is locally presentable. More generally, if $\cM$ is a locally presentable, symmetric monoidal category, then the category $\cCat_\cM$ of $\cM$-enriched categories is locally presentable (see~\cite{kelly-lack:vcat}). 
\item The category $\cTop$ of topological spaces is \emph{not} locally presentable, but this can be fixed by passing to the Quillen equivalent category of $\Delta$-generated spaces. It is shown in \cite{fajstrup-rosicky} that this category is locally presentable.
\end{enumerate}
\end{egs}

Here is the simplified form of the adjoint functor theorem for locally presentable categories.

\begin{thm}\label{thm:catSAFT}
\begin{enumerate}
\item A functor between locally presentable categories is a left adjoint if and only if it preserves colimits.
\item A functor between locally presentable categories is a right adjoint if and only if it preserves limits and is accessible.
\end{enumerate}
\end{thm}

There are variants of the first part of the theorem in slightly different settings as the classical Brown representability theorem in stable homotopy theory \cite{brown:cohomology}, the Brown representability results for triangulated categories \cite[\S8]{neeman:triangulated}, Watt's theorems in homological algebra \cite[\S5.3]{rotman}, representability theorems for Grothendieck categories \cite[p.186]{Schapira}, as well as versions of Watt's theorem in homotopical algebra~\cite{hovey:watts}.

It turns out that up to equivalence \emph{all} locally presentable categories can be obtained as certain localizations of presheaf categories. To state this more precisely, let us begin by recalling that a \textbf{reflective localization} is an adjunction $(L,R)\colon\cC\rightleftarrows\cD$ such that the right adjoint $R$ is fully faithful (see for example \cite[\S3.5 and~\S5.3]{borceux:1}). In this situation it follows that $\cD$ is equivalent to the localization $\cC[S^{-1}]$ where $S$ is the class of morphisms in $\cC$ which are sent to isomorphisms by~$L.$ This is nicely illustrated by the following example.

\begin{eg}
Let $\delta\colon\cSet\to\sSet$ be the discrete simplicial set functor. The adjunction 
$(\pi_0,\delta)\colon\sSet\rightleftarrows\cSet$
is a reflective localization. Thus, if in $\sSet$ we invert all maps inducing isomorphisms on $\pi_0$ then we simply get (discrete simplicial) sets.
\end{eg} 

We want to emphasize that the typical terminology from \emph{Bousfield localization theory} \cite{bousfield:localization} makes already perfectly well sense in this classical context, and that reflective localizations can be nicely described using that terminology (see \autoref{prop:Bousfield}). Since this does not seem to be over-emphasized in the literature, we even include a proof of \autoref{prop:Bousfield}.

\begin{defn}\label{def:catSlocal}
Let $\cC$ be a category and let $S$ be a class of morphisms in $\cC$.
\begin{enumerate}
\item An object $c$ in $\cC$ is \textbf{$S$-local} if $f^\ast\colon\hom_{\cC}(c_2,c)\to\hom_{\cC}(c_1,c)$ is a bijection for all $f\colon c_1\to c_2$ in $S$.
\item A morphism $f\colon c_1\to c_2$ in $\cC$ is an \textbf{$S$-local equivalence} if for all $S$-local objects $c\in\cC$ the map $f^\ast\colon\hom_{\cC}(c_2,c)\to\hom_{\cC}(c_1,c)$ is a bijection.
\end{enumerate}
\end{defn}

\begin{prop}\label{prop:Bousfield}
Let $(L,R)\colon\cC\rightleftarrows\cD$ be a reflective localization and let $S$ be the morphisms in $\cC$ which are inverted by $L$.
\begin{enumerate}
\item The essential image of $R$ consists precisely of the $S$-local objects.
\item The $S$-local equivalences are precisely the maps in $S$.
\end{enumerate}
\end{prop}
\begin{proof}
We begin with a proof of the first statement and consider $c=Rd$ in the image of $R$ together with a morphism $f\colon c_1\to c_2$ in $S$. The naturality of the adjunction isomorphism
\begin{equation}
\vcenter{
\xymatrix@=1.5pc{
\hom_{\cC}(c_2,Rd)\ar[r]^-\cong\ar[d]_-{f^\ast}&\hom_{\cD}(Lc_2,d)\ar[d]^-{(Lf)^\ast}\\
\hom_{\cC}(c_1,Rd)\ar[r]_-\cong&\hom_{\cD}(Lc_1,d)
}
}
\label{eq:natsq}
\end{equation}
implies that $Rd$ is $S$-local since $Lf$ is an isomorphism by definition of $S$. 

The converse implication is more involved. Let $c\in\cC$ be an $S$-local object. We want to show that $c$ lies in the essential image of $R$. By abstract nonsense this is the case if and only if the adjunction unit $\eta_c\colon c\to RLc$ is an isomorphism. Since the right adjoint functor $R$ is fully faithful, the adjunction counit $\epsilon\colon LR\to \id$ is an isomorphism. It follows from the triangular identity $\id=\epsilon_L\circ L\eta\colon L\to LRL\to L$ that $L\eta$ is an isomorphism, i.e., that the components of $\eta$ lie in $S$. Since $c$ is $S$-local we obtain a bijection
\[
\eta_c^\ast\colon\hom_{\cC}(RLc,c)\to\hom_{\cC}(c,c),
\]
and hence a unique $g\colon RLc\to c$ such that $g\circ\eta_c=\id$. To conclude this converse direction it suffices to show that also $\eta_c\circ g$ is the identity. First, the relation $g\circ\eta_c=\id$  implies by an additional application of the above triangular identity that
$Lg\circ L\eta_c=\id_{Lc}=\epsilon_{Lc}\circ L\eta_c$, and hence that $Lg=\epsilon_{Lc}$ since $L\eta_c$ is an isomorphism. This equality together with the commutativity of
\[
\xymatrix@=1.5pc{
RLc\ar[r]^-g\ar[d]_-{\eta_{RLc}}&c\ar[d]^-{\eta_c}\\
RLRLc\ar[r]_-{RLg}&RLc
}
\]
yields the following chain of equalities
\begin{align}
\eta_c\circ g&=RLg\circ \eta_{RLc}\\
&=R\epsilon_{Lc}\circ\eta_{RLc}\\
&=(R\epsilon\circ\eta_R)_{Lc}\\
&=\id_{RLc},
\end{align}
where the last step is again given by a triangular identity. This concludes the proof of the first statement. The second statement follows now easily from \eqref{eq:natsq} and the Yoneda lemma.
\end{proof}

With this preparation, we now state the `classification result of locally presentable categories' (see for example \cite{adamek-rosicky}). Let us recall that an \textbf{accessible, reflective localization} is a reflective localization $(L,R)$ such that the right adjoint~$R$ is accessible, i.e., preserves $\kappa$-filtered colimits for a sufficiently large regular cardinal~$\kappa$. 

\begin{thm}\label{thm:locpres}
A category is locally presentable if and only if it is equivalent to an accessible, reflective localization of $\nFun(A\op,\cSet)$ for some small category~$A$.
\end{thm}

Thus, a category is locally presentable if and only if it can be obtained from a small category by a free generation (cocompletion) followed by imposing relations in a suitable way (accessible reflective localization). In \S\ref{subsec:combmod} and \S\ref{subsec:pres} we will see that there are variants of \autoref{thm:yoneda} and \autoref{thm:locpres} valid in the context of model categories or $\infty$-categories.

\subsection{Combinatorial model categories}
\label{subsec:combmod}

One reason why \emph{set-valued presheaves} play such a central role in classical category theory is that questions about the existence of universal constructions can be reformulated as representability questions for certain set-valued presheaves. In higher category theory, the representable functors take values in the category of simplicial sets. So one might expect that the central role is now taken by \emph{simplicial presheaf categories}. In the world of model categories, these were intensively studied by Jardine (see for example \cite{jardine:presheaves}). 

Recall that the category $\sSet$ of simplicial sets can be endowed with the Kan--Quillen model structure. Since this model structure is cofibrantly-generated it follows that, for every small category~$A$, the category $\nFun(A\op,\sSet)$ inherits the projective model structure~\cite{bousfield-kan:htpy-limits}. Let us recall that a morphism $f\colon X\to Y$ in $\nFun(A\op,\sSet)$ is a
\begin{enumerate}
\item \textbf{projective fibration} if each $f_a\colon X_a\to Y_a$ is a Kan fibration,
\item \textbf{projective weak equivalence} if each $f_a\colon X_a\to Y_a$ is a weak equivalence,
\item and a \textbf{projective cofibration} if it has the left lifting property with respect to all maps which are simultaneously projective fibrations and projective weak equivalences.
\end{enumerate}
Dugger shows in \cite{dugger:universal} that these projective model structures on simplicial presheaves are `universal homotopy theories', i.e., can be thought of as `homotopical cocompletions'. In order to make this more precise let us recall that the Yoneda embedding $y_A\colon A\to \nFun(A\op,\cSet)$ is \textbf{dense}, which simply means that every presheaf is canonically a colimit of representable ones (as made precise by \autoref{prop:Yonedadense}). 

Of course, since $\nFun(A\op,\sSet)$ is itself isomorphic to the presheaf category $\nFun((A\times\Delta)\op,\cSet)$ it follows that $A\times\Delta$ is densely embedded in $\nFun(A\op,\sSet)$, essentially by the Yoneda embedding $y_{A\times\Delta}$. But this is not what we want to do here. Instead, we want to combine the discrete simplicial set functor $\delta\colon\cSet\to\sSet$ with the Yoneda embedding $y_A\colon A\to \nFun(A\op,\cSet)$ in order to obtain
\begin{equation}
\delta_\ast\circ y_A\colon A\to \nFun(A\op,\cSet)\to\nFun(A\op,\sSet).
\label{eq:homotopicallydense}
\end{equation} 

This can be used to make precise that every simplicial presheaf is canonically a \emph{homotopy colimit} of simplicial presheaves in the image of $\delta_\ast\circ y_A$. The universality of the homotopy theory
\[
U(A)=\nFun(A\op,\sSet)_{\mathrm{proj}}
\]
is expressed by the following theorem, but we refer the reader to \cite{dugger:universal} for more precise formulations of both the theorem and the `homotopical density' of \eqref{eq:homotopicallydense}. The theorem is a model category theoretic variant of \autoref{thm:yoneda}.

\begin{thm}[\cite{dugger:universal}]
Let $A$ be a small category and let $\cM$ be a model category. Every functor $Q\colon A\to\cM$ factors over the model category $U(A)$, i.e., there is a left Quillen functor $L\colon U(A)\to\cM$ together with a natural weak equivalence
\[
\xymatrix@=1.5pc{
A\ar[r]^-{\delta_\ast\circ y_A}\ar@/_0.8pc/[dr]_-Q\drtwocell\omit{}&U(A)\ar[d]^-L\\
&\cM.
}
\]
Moreover, the category of such factorizations is contractible.
\end{thm}

The idea of the proof is as follows. Since the functor $\delta_\ast\circ y_A$ is homotopically dense, every simplicial presheaf $X$ is canonically a homotopy colimit of things coming from~$A$. A left Quillen functor necessarily has to preserve these homotopy colimits, and for the building blocks of the corresponding diagrams the functor $Q$ already dictates the values in~$\cM$. Thus, one formally realizes the canonical homotopy colimits of simplicial presheaves in the target model category~$M$ by means of the `test diagram' $Q$. The fact that we only expect factorizations up to natural weak equivalences is related to cofibrancy issues; see \cite{dugger:universal}. For an alternative proof of this universality see also \cite{rosicky-tholen:left-determined}.

Having a basic understanding of `homotopical cocompletions', we now turn to an analogue of locally presentable categories in the framework of model categories. The following definition is due to J.~Smith (unpublished). Reference for this notion include \cite{beke:sheafifiable,dugger:combinatorial,HTT,rosicky:combinatorial}. 

\begin{defn}
A model category is \textbf{combinatorial} if the model structure is cofibrantly generated and if the underlying category is locally presentable.
\end{defn}

This class of model categories has good technical properties. For example diagram categories $\nFun(A,\cM)$ with values in a combinatorial model category can always be endowed with both the projective and the injective model structure. This makes the discussion of homotopy colimits and homotopy limits particularly tractable. A further convenient property is that (left proper) combinatorial model categories admit a good theory of (Bousfield) localizations. Bousfield localizations are a homotopical analogue of the reflective localizations in ordinary category theory. They are based on homotopical versions of the notions of $S$-local objects and $S$-local equivalences (\autoref{def:catSlocal}) which in turn are obtained by replacing the morphism sets by homotopy function complexes. For more details on these homotopy function complexes see \cite{dwyer-kan:simplicial, dwyer-kan:calculating, dwyer-kan:function} and \cite{hovey:model,hirschhorn:model}, and for details on Bousfield localization theory see \cite{bousfield:localization} and \cite[Part~1]{hirschhorn:model}.

In \cite{dugger:combinatorial} Dugger shows that `combinatorial model categories have presentations'. By this he means that up to Quillen equivalence, every combinatorial $\cM$ can be obtained as a left Bousfield localization of a universal homotopy theory in the sense of \cite{dugger:universal}. 

\begin{thm}[\cite{dugger:combinatorial}]
Every combinatorial model category~$\cM$ has a presentation. More precisely, there is a small category~$A$, a set $S$ of morphisms in $U(A)$, and a left Quillen equivalence $U(A)[S^{-1}]\to\cM$, where $U(A)[S^{-1}]$ denotes the left Bousfield localization.
\end{thm}

This theorem is a model category theoretic analogue of \autoref{thm:locpres}.

\subsection{Presentable $\infty$-categories}
\label{subsec:pres}

Similarly to the case of model categories also in the context of $\infty$-categories a key role is played by $\infty$-categories of simplicial presheaves. Recall that we denote the $\infty$-category of spaces by $\cS=N_\Delta(\cKan)$ (see \autoref{egs:oocats}). Given a small simplicial set~$K$, the $\infty$-category $\cP(K)$ of \textbf{(simplicial) presheaves} on $K$ is defined by
\[
\cP(K)=\nFun(K\op,\cS).
\]
It follows from \autoref{prop:functor} that $\cP(K)$ is an $\infty$-category. 

In order to define the Yoneda embedding let us recall from \S\ref{subsec:joyalbergner} that we have the adjunction $(C[-],N_\Delta)\colon\sSet\rightleftarrows\sCat$; see \eqref{eq:CN-adj}. Given an arbitrary simplicial set~$K$, then, in general, the simplicial category $C[K]$ is not locally fibrant (\autoref{per:mapping}). This can be fixed by choosing a product-preserving fibrant replacement functor for $\sSet$, like the one induced by the Quillen equivalence $\sSet\rightleftarrows\cTop$ or Kan's $\mathrm{Ex}^\infty$-functor; see \cite{kan:css}. Composing the mapping space functor of $C[K]$ with such a replacement functor, we obtain a simplicial functor $C[K]\op\times C[K]\to\cKan$. Combining this with the canonical map $C[K\op\times K]\to C[K]\op\times C[K]$ and passing to adjoints, this yields a map $K\op\times K\to N_\Delta(\cKan)=\cS$. The exponential law finally gives us the \textbf{Yoneda embedding}
\begin{equation}\label{eq:Yoneda}
y\colon K\to\nFun(K\op,\cS)=\cP(K),
\end{equation}
which can be shown to be fully faithful (\cite[Prop.~5.1.3.1]{HTT}).

The Yoneda embedding provides a model for the cocompletion. In order to make this precise let us introduce the following notation. Given $\infty$-categories $\cC$ and $\cD$, we denote by 
\[
\nFunL(\cC,\cD)\subseteq\nFun(\cC,\cD)
\]
the full subcategory spanned by the \emph{colimit-preserving} functors. Lurie establishes the following result (\cite[Thm.~5.1.5.6]{HTT}) which is an $\infty$-categorical version of \autoref{thm:yoneda}.

\begin{thm}\label{thm:ooyoneda}
Let $K$ be a small simplicial set and let $\cC$ be a cocomplete $\infty$-category. Restriction along the Yoneda embedding \eqref{eq:Yoneda} induces an equivalence of $\infty$-categories
\[
\nFunL(\cP(K),\cC)\stackrel{\sim}{\to}\nFun(K,\cC).
\]
\end{thm}

In particular, the \emph{$\infty$-category $\cS$ of spaces is freely generated by $\Delta^0\in\cS$} in the sense of the following corollary.

\begin{cor}\label{cor:spaces}
For any cocomplete $\infty$-category~$\cC$ the evaluation on the 0-simplex $\Delta^0\in\cS$ induces an equivalence of $\infty$-categories
\[
\nFunL(\cS,\cC)\stackrel{\sim}{\to}\cC.
\]
\end{cor}

In \cite[\S4]{HTT} Lurie establishes batteries of techniques which allow us to manipulate (co)limits in the context of $\infty$-categories. In \cite[\S5]{HTT} more advanced notions are introduced, including filtered colimits and small objects. Once the basic notions are in place, more advanced concepts from classical category theory can be formally extended to the $\infty$-categorical framework. We want to emphasize once more that although impressively many \emph{statements} are still true in this more general framework, at least at present the \emph{proofs}, in general, are more involved. 

\begin{defn}
An $\infty$-category is \textbf{presentable} if it is cocomplete and accessible.
\end{defn}

Similarly to ordinary category theory, it can be shown that a presentable $\infty$-category is automatically also complete (\cite[Cor.~5.5.2.4]{HTT}). There are also $\infty$-categorical variants of Freyd's special adjoint functor theorem (see \autoref{thm:catSAFT}). A first step towards this of course consists of making precise the notion of an adjunction between two $\infty$-categories. We will not get into this here and instead refer the reader to \cite[p.337]{HTT} and \cite[\S3]{DAGII}. With this concept at hand, there is the following result due to Lurie (\cite[Cor.~5.5.2.9]{HTT}).

\begin{thm}\label{thm:SAFT}
\begin{enumerate}
\item A functor between presentable $\infty$-categories is a left adjoint if and only if it preserves colimits.
\item A functor between presentable $\infty$-categories is a right adjoint if and only if it preserves limits and is accessible.
\end{enumerate}
\end{thm}

Lurie also establishes a classification result for presentable $\infty$-categories similar to the one given by \autoref{thm:locpres}. Using the concept of an adjunction of $\infty$-categories, one can extend additional concepts from ordinary category theory to the context of $\infty$-categories. A functor between two $\infty$-categories is a \textbf{reflective localization} if it admits a fully faithful right adjoint, and there is also the notion of an \textbf{accessible, reflective localization}. With these definitions at hand, the classification result for presentable $\infty$-categories takes the following form as proved by Lurie as \cite[Thm.~5.5.1.1]{HTT}. In that section, Lurie attributes the result to Simpson~\cite{simpson:giraud}.

\begin{thm}\label{thm:normal}
For an $\infty$-category $\cC$ the following are equivalent.
\begin{enumerate}
\item The $\infty$-category $\cC$ is presentable.
\item There is a small $\infty$-category $\cD$ such that $\cC$ is an accessible, reflective localization of~$\mathcal{P}(\cD)$.
\end{enumerate}
\end{thm}

In light of \autoref{thm:normal}, if one wants to understand presentable $\infty$-categories then it seems to be important to have good control over accessible localizations of $\infty$-categories of presheaves or, more generally, of presentable $\infty$-categories. Given a reflective localization $(L,R)\colon\cC\rightleftarrows\cD$ let us also write $L\colon\cC\to\cC$ for the \textbf{localization functor} $\cC\to\cD\to\cC,$ and let $S_L$ be those maps in $\cC$ which are sent to equivalences by~$L.$ An object $x\in\cC$ is \textbf{$S_L$-local} if all maps $f^\ast\colon\Map_{\cC}(z,x)\to\Map_{\cC}(y,x)$ induced by $f\in S_L$ are weak equivalences. One can then show that the essential image of $L$ consists precisely of the $S_L$-local objects (see \cite[Prop.~5.5.4.2]{HTT}).

In particular, an accessible localization~$L$ of a presentable $\infty$-category~$\cC$ is thus completely determined by the class $S_L.$ In that case, the class $S_L$ is closed under the formation of colimits in $S_L$ as a subcategory of $\cC^{[1]}$, is stable under the formation of retracts, contains the equivalences, satisfies the 2-out-of-3 property (with respect to 2-simplices), and is stable under cobase change. Lurie calls a class of morphisms with these closure properties \textbf{strongly saturated}. Intersections of strongly saturated classes are again strongly saturated as is the class of all morphisms. Thus, for each class~$T$ of morphisms in $\cC$ there is a smallest strongly saturated class $\bar{T}$ which contains~$T.$ We call a strongly saturated class $S$ \textbf{of small generation} if there is a \emph{subset} $T\subseteq S$ such that $S=\bar{T}.$ Lurie establishes the wonderful fact that a strongly saturated class $S$ in a presentable $\infty$-category $\cC$ is of small generation if and only if there is an accessible localization $L\colon\cC\to\cC$ such that $S=S_L$ (see \cite[\S5.5.4]{HTT}).

Lurie then shows that in the case of simplicial presheaves the localization theory of $\infty$-categories interacts nicely with the localizaton theory of certain associated model categories (see \cite[Appendix~3.7]{HTT}). Having established all this theory, he is then able to build on Dugger's work \cite{dugger:combinatorial} in order to deduce the following result (see \cite[Appendix~3]{HTT}).

\begin{thm}\label{thm:comparison}
An $\infty$-category $\cC$ is presentable if and only if there is a combinatorial, simplicial model category $\mathcal{M}$ such that $\cC$ is equivalent to the underlying $\infty$-category $N_{\Delta}(\mathcal{M}_{\ncf}).$
\end{thm}

In this section we considered the `theme of locally presentable categories' in three different frameworks, namely in ordinary category theory, in model category theory, and in $\infty$-category theory. We conclude this section with a short perspective in which we give an outlook on similar pictures for \emph{Grothendieck topoi} and \emph {algebraic categories}.

\begin{per}\label{per:topoi}
\begin{enumerate}
\item Let us recall that a \emph{Grothendieck topos} can be defined as a category equivalent to a category of set-valued sheaves on a Grothendieck site; references for this vast subject include the original \cite{SGA41,SGA42,SGA43} and the monographs~\cite{borceux:3,maclane-moerdijk,moerdijk:classifying}. It turns out that Grothendieck topoi admit a different characterization as suitable localizations of presheaf categories. In fact, they are precisely the \emph{left exact, reflective, accessible localizations} of presheaf categories. In more detail, a category~$\cC$ is a Grothendieck topos if and only if there is small category~$A$ and an adjunction
\[
(L,R)\colon \nFun(A\op,\cSet)\rightleftarrows\cC,
\]
such that the following three properties hold.
\begin{enumerate}
\item The right adjoint~$R$ is fully faithful, i.e., we have a reflective localization.
\item The right adjoint~$R$ is accessible.
\item The left adjoint~$L$ is left exact, i.e., it preserves finite limits.
\end{enumerate}
Thus, Grothendieck topoi are, in particular, locally presentable categories.

The `theme of sheaves and topoi' was taken up again in homotopical frameworks, both using the language of model categories (see for example \cite{jardine:presheaves,jardine:boolean,dsi:hyper,toen-vezzosi:hag,rezk:toposes}) as well as in the $\infty$-categorical picture \cite[\S\S6-7]{HTT}.

\item A different theme that was taken up in all these frameworks is the `theme of algebraic theories and algebraic categories'. Let us recall from the more categorical approach~\cite{lawvere-algebraic} to universal algebra~\cite{manes:algebraic,grätzer:universal,borceux:2,adamek-rosicky-vitale} that an algebraic theory (or a multi-sorted one) is simply a small category~$\cT$ with finite products. An algebra for a theory~$\cT$ is given by a product preserving functor $\cT\to\cSet$, yielding the category $\cT\text{-}\nAlg$ of $\cT$-algebras. An \emph{algebraic category} is a category which is equivalent to $\cT\text{-}\nAlg$ for some algebraic theory~$\cT$. Similarly to the case of locally presentable categories and Grothendieck topoi, it turns out that algebraic categories can be characterized as suitable localizations of presheaf categories. In fact, they are precisely the \emph{reflective, sifted localizations} of presheaf categories. In more detail, a category $\cC$ is algebraic if and only if there is small category~$A$ and an adjunction
\[
(L,R)\colon \nFun(A\op,\cSet)\rightleftarrows\cC,
\]
such that the following two properties hold (see \cite[Theorem~6.18]{adamek-rosicky-vitale}).
\begin{enumerate}
\item The right adjoint~$R$ is fully faithful, i.e., we have a reflective localization.
\item The right adjoint~$R$ preserves \emph{sifted} colimits.
\end{enumerate}
We do not want to get into the notion of sifted colimits here and only mention that sifted colimits include filtered colimits and reflective coequalizers; for more details see for example~\cite{adamek-rosicky:sifted,adamek-rosicky-vitale:sifted}. However, this already implies that a functor which preserves sifted colimits is, in particular, accessible and it hence follows that algebraic categories are locally presentable.

In the spirit of \cite{boardman-vogt,segal:categories}, the `theme of algebraic categories and algebraic theories' was reconsidered in homotopical frameworks, aiming for a study of homotopy coherent versions of such algebraic structures. References in the case of model categories include \cite{badzioch:algebraic,bergner:rigidification,rosicky:homotopy-varieties,schwede:algebraic-theories}, while in the framework of $\infty$-categories there are \cite[\S5.5.8]{HTT} and \cite{cranch:algebraic,gepner-groth-nikolaus}.
\end{enumerate}
\end{per}

\section{Monoidal and symmetric monoidal $\infty$-categories}
\label{sec:monoidal}

In this section we give an introduction to the theory of monoidal $\infty$-categories. Let us recall that a monoidal structure on a category $\cM$ consists of a monoidal pairing $\otimes\colon\cM\times\cM\to\cM$ and a monoidal unit $\lS\in\cM$ together with three natural isomorphisms, namely the associativity and left and right unitality constraints. Moreover, to obtain the notion we have in mind we have to impose certain compatibility assumptions. In particular, we have to ask axiomatically that the two ways of comparing four-fold products $\big((X\otimes Y)\otimes Z\big)\otimes W$ and $X\otimes \big(Y\otimes (Z\otimes W)\big)$ as described by the boundary in 
\begin{equation}
\vcenter{
\xymatrix@=0.5pc{
\big((X\otimes Y)\otimes Z\big)\otimes W\ar[rr]\ar[d]&&(X\otimes Y)\otimes(Z\otimes W)\ar[d]\\
\big(X\otimes(Y\otimes Z)\big)\otimes W\ar[rd]&&X\otimes \big(Y\otimes (Z\otimes W)\big)\\
&X\otimes\big((Y\otimes Z)\otimes W\big)\ar[ru]&
}
\label{eq:monoidal-intro}
}
\end{equation}
coincide. This leads to the \emph{classical presentation} of monoidal categories (\cite{maclane:natural,kelly:natural}).

In the $\infty$-categorical setting this presentation will not model the good notion anymore. Instead one expects a monoidal structure on an $\infty$-category to be some kind of a monoidal pairing which is coherently associative and unital in the sense of $\Aoo$-multiplications. In particular, it is hence insufficient to consider Mac~Lane's pentagon ~\eqref{eq:monoidal-intro} and instead one expects that all Stasheff associahedra with their complicated combinatorics play a key role (see \cite{stasheff-Aoo} or \cite[\S I.1.6 and \S II.1.6]{markl-shnider-stasheff}). 

Luckily, all this structure does not have to be made explicit if one chooses a \emph{different presentation} of ordinary monoidal categories, namely, as suitable Grothendieck opfibrations over~$\Delta\op$ (the same observation but in the context of $\Aoo$-spaces motivated Adams to refer to the category $\Delta$ as a `storehouse of formulas' \cite{adams:infinite}). This different presentation extends more readily to $\infty$-categories and is obtained by combining two main ideas, namely

\begin{enumerate}
\item the Segal perspective on $\mathbb {A}_\infty$-monoids as `special simplicial objects' and
\item the Grothendieck construction applied to category-valued functors.
\end{enumerate}

The passage to \emph{symmetric} monoidal $\infty$-categories then essentially amounts to a change of combinatorics, i.e., one replaces the category $\Delta\op$ by a skeleton $\Fin$ of the category of finite pointed sets.

\subsection{Monoidal categories via Grothendieck opfibrations}
\label{subsec:monoidal}

Let us consider the following classical situation. Let $p\colon\cC\to\cD$ be a functor between ordinary categories and let $d\in\cD$ be an object. We denote by $\cC_d$, the \textbf{fiber} of~$p$ over $d,$ defined by the following pullback diagram
\[
\xymatrix@=1.5pc{
\cC_d\pullbackcorner\ar[r]\ar[d]&\cC \ar[d]^p\\
[0]\ar[r]_d&\cD.
}
\]
Thus, $\cC_d\subseteq\cC$ is the (in general not full) subcategory given by the objects $c\in\cC$ such that $p(c)=d$ and those morphisms in $\cC$ which are sent to $\id_d.$ A functor $p\colon\cC\to\cD$ can always be thought of as a collection of categories $\cC_d$ parametrized by the objects in $\cD.$ 

Our first aim is to find conditions which ensure that the fiber~$\cC_d$ \emph{depends covariantly} on the object $d\in\cD.$

\begin{defn}\label{defn:coCartesianarrow}
Let $p\colon\cC\to\cD$ be a functor and let $f\colon c_1\to c_2$ be a morphism in~$\cC$ with image $p(f)=\alpha\colon d_1\to d_2.$ The morphism $f$ is \textbf{$p$-coCartesian} or a \textbf{$p$-coCartesian lift} of $\alpha$ if it has the following property: For every $h\colon c_1\to c_3$ in~$\cC$ with image $\gamma=p(h)\colon d_1\to d_3$ and every $\beta\colon d_2\to d_3$ such that $\gamma=\beta\circ\alpha$ there is a unique $g\colon c_2\to c_3$ in $\cC$ such that 
\[
\beta=p(g)\qquad\text{and}\qquad h=g\circ f.
\]
\end{defn}

Thus, the defining property of a $p$-coCartesian arrow can be described by the following diagram
\[
\xymatrix@=1.5pc{
c_1\ar[r]\ar@{-}[dd]\ar@/_1.4pc/[rrd]_-\forall&c_2\ar@{-}'[d][dd]\ar@{-->}[rd]^-{\exists !}&\\
&&c_3\ar@{-}[dd]\\
d_1\ar[r]\ar@/_1.4pc/[rrd]\ar@{}[rrd]_=&d_2\ar[rd]^-\forall&\\
&&d_3\\
}
\]
in which the vertical lines indicate the effect of an application of $p\colon\cC\to\cD.$ Note that a morphism $f\colon c_1\to c_2$ is $p$-coCartesian if and only if the diagram 
\begin{equation}
\vcenter{
\xymatrix@=1.5pc{
\hom_\cC(c_2,c_3)\ar[r]^-{f^\ast}\ar[d]_-p&\hom_\cC(c_1,c_3)\ar[d]^-p\\
\hom_\cD(p(c_2),p(c_3))\ar[r]_-{p(f)^\ast}&\hom_\cD(p(c_1),p(c_3))
}\label{eq:coCartesian}
}
\end{equation}
is a pullback diagram for all objects $c_3\in\cC.$ This slightly cryptic  reformulation of \autoref{defn:coCartesianarrow} will have its uses when it comes to extending these notions to the setting of $\infty$-categories (see \autoref{lem:coCartesian} and \autoref{defn:oocoCartesian}). To develop some feeling for the notion we recommend the reader to give the easy proof of the following lemma.

\begin{lem}\label{lem:uniquecoC}
Let $f'\colon c\to c'$ and $f''\colon c\to c''$ be $p$-coCartesian arrows with the same image $\alpha=p(f')=p(f'').$ Then there is a unique isomorphism $\phi\colon c'\to c''$ in the fiber $\cC_{p(c')}=\cC_{p(c'')}$ such that $\phi\circ f'=f''.$
\end{lem}

This lemma tells us that $p$-coCartesian lifts with a fixed domain are essentially unique if they exist. In particular, for $\alpha\colon d_1\to d_2$ the targets of $p$-coCartesian lifts are (uniquely compatibly) isomorphic as objects in the fiber $\cC_{d_2}.$ Thus, in order to obtain a covariant dependence of the fiber it seems to be a good strategy to ask for a sufficient supply of $p$-coCartesian arrows.

\begin{defn}\label{defn:Grothendieck}
A functor $p\colon\cC\to\cD$ is a \textbf{Grothendieck opfibration} if for all $c_1\in\cC$ and for all morphisms $\alpha$ in $\cD$ with domain $p(c_1)$ there is a $p$-coCartesian lift $f\colon c_1\to c_2$ of $\alpha.$
\end{defn}

Let $p\colon\cC\to\cD$ be a Grothendieck opfibration. Then we can \emph{choose} for each $c\in\cC$ and for each morphism $\alpha\colon p(c)\to d$ a $p$-coCartesian lift. We now fix a morphism $\alpha\colon d_1\to d_2$ in $\cD$ and define
\[
\alpha_!\colon\cC_{d_1}\to\cC_{d_2}\colon c_1\mapsto c_2,
\]
where $c_2$ is the codomain of the \emph{chosen} $p$-coCartesian lift $f\colon c_1\to c_2$ of $\alpha.$ This defines $\alpha_!$ on objects, and we recommend the reader to check that this can be extended to define a functor $\alpha_!\colon\cC_{d_1}\to\cC_{d_2}.$ 

If we now consider an additional morphism $\beta\colon d_2\to d_3$ in $\cD,$ then we obtain associated functors 
\[
\xymatrix{
\cC_{d_1}\ar[r]^-{\alpha_!}&\cC_{d_2}\ar[r]^-{\beta_!}&\cC_{d_3},&
\cC_{d_1}\ar[r]^-{(\beta\circ\alpha)_!}&\cC_{d_3}.
}
\] 
In general, these two functors are not equal since their definitions depend on certain choices of $p$-coCartesian lifts: they send an object $c_1\in\cC_{d_1}$ to the respective targets of two possibly different lifts of $\beta\circ\alpha$ to $p$-coCartesian morphisms with domain $c_1$ (Exercise: The composition of two $p$-coCartesian morphisms is again $p$-coCartesian.). But one can deduce from \autoref{lem:uniquecoC} that there is a unique natural isomorphism
\begin{equation}
\beta_!\circ \alpha_!\cong(\beta \circ \alpha)_! \label{eq:coherence}
\end{equation}
of functors $\cC_{d_1}\to\cC_{d_3},$  i.e., all components of the natural isomorphism are sent to the identity of $d_3$ via $p$.

As an upshot, we essentially succeeded in obtaining a covariant dependence of the fiber by considering Grothendieck opfibrations. To put this in a slightly technical language, we have observed that for such functors the fibers depend \emph{pseudo-functorially} on $d\in\cD$, i.e., that the assignment $d\mapsto\cC_d$ defines a pseudo-functor $\cD\to\cCAT$. One might be disappointed about this lack of strict functoriality, but for our purposes this is very convenient: it allows us to \emph{encode} or, better, \emph{hide a lot of structure} in the natural isomorphisms \eqref{eq:coherence} associated to a Grothendieck opfibration. We illustrate this by the following example.

\begin{eg}
Let $\cM$ be a monoidal category with monoidal pairing $\otimes\colon\cM\times\cM\to\cM$ and monoidal unit $\lS\in\cM.$ We form a new category $\cM^\otimes$ in the following way. The objects of $\cM^\otimes$ are (possibly empty) finite sequences of objects in $\cM,$ 
\[
(M_1,\ldots,M_n),\quad n\geq 0,\;M_i\in\cM.
\] 
Given two such sequences $(M_1,\ldots,M_n)$ and $(L_1,\dots,L_k),$ a morphism 
\[
(\alpha,\{f_i\}_i)\colon (M_1,\dots,M_n)\to(L_1,\dots,L_k)
\] 
consists of a morphism $\alpha\colon[k]\to[n]$ in $\Delta$ together with morphisms
\begin{equation}
f_i\colon M_{\alpha(i-1)+1}\otimes\ldots \otimes M_{\alpha(i)}\to L_i,\quad i=1,\dots,k.\label{eq:components}
\end{equation}
Thus given such a morphism, $\alpha$ encodes the domains of the $f_i.$ In particular, if there is an $i\in [k]$ such that $\alpha(i-1)=\alpha (i),$ then by convention the corresponding map \eqref{eq:components} is to be read as a map $f_i\colon \lS\to L_i.$ The composition of morphisms in $\cM^\otimes$ is defined using the compositions in $\Delta$ and $\cM$ together with the associativity constraints of the monoidal structure on $\cM.$ The identity of an object $(M_1,\ldots,M_n)$ is easily seen to be given by $(\id_{[n]},\{\id_{M_i}\}_i).$

There is an obvious projection functor $p\colon\cM^\otimes\to\Delta\op$ which sends a string $(M_1,\ldots,M_n)$ to $[n]$ and a morphism $(\alpha,\{f_i\}_i)$ to its first component~$\alpha$. One easily checks that $p\colon\cM^\otimes\to\Delta\op$ is a Grothendieck opfibration. Indeed, let us consider a lifting problem given by an object $(M_1,\dots,M_n)$ in $\cM^\otimes_{[n]}$ together with a morphism $\alpha\op\colon[n]\to[k]$ in $\Delta\op,$
\[
\xymatrix@=1.5pc{
(M_1,\ldots,M_n)\ar@{-}[d]&\\
[n]\ar[r]_-{\alpha\op}&[k],\\
[n]&[k].\ar[l]^-{\alpha}
}
\]
Then a $p$-coCartesian lift of $\alpha$ with domain the given string is obtained from any family of isomorphisms 
\[
f_i\colon M_{\alpha(i-1)+1}\otimes\dots\otimes M_{\alpha(i)}\stackrel{\cong}{\to} L_i,\quad i=1,\dots,k.
\]
More precisely, these $L_i$ specify an object $(L_1,\dots, L_k)\in\cM^\otimes_{[k]},$ and the morphism 
\begin{equation}
(\alpha,\{f_i\}_i)\colon (M_1,\dots,M_n)\to(L_1,\dots ,L_k)\label{eq:coCartesianlift}
\end{equation}
is the desired $p$-coCartesian lift. 

Now, this Grothendieck opfibration $p\colon\cM^\otimes\to\Delta\op$ has the property that the fiber $\cM_{[n]}$ is canonically equivalent to the $n$-fold product of $\cM^\otimes_{[1]}\simeq\cM$ in the following sense. Let $\iota_{\{i-1,i\}}\colon[1]\to[n]$ be the inclusion of the $i$-th length one interval, i.e., the unique monomorphism in $\Delta$ with image $\{i-1,i\},$ and let us write $\iota_i=\iota_{\{i-1,i\}}\op\colon [n]\to[1]$ for the opposite morphism in $\Delta\op.$ Since $p\colon\cM^\otimes\to\Delta\op$ is a Grothendieck opfibration, we obtain induced functors
\begin{equation}
(\iota_i)_!\colon\cM^\otimes_{[n]}\to\cM^\otimes_{[1]}=\cM,\qquad i=1,\ldots,n,\label{eq:projection}
\end{equation}
which, taken together, induce the \textbf{Segal maps}
\begin{equation}
\sigma=((\iota_1)_!,\ldots,(\iota_n)_!)\colon\cM^\otimes_{[n]}\stackrel{\simeq}{\longrightarrow}\underbrace{\cM\times\ldots\times\cM}_{n\text{ times}}.\label{eq:Segal}
\end{equation}
Note that the explicit construction of $p$-coCartesian lifts in \eqref{eq:coCartesianlift} implies that these Segal maps are equivalences. We refer to this observation by saying that the Grothendieck opfibration $p\colon\cM^\otimes\to\Delta\op$ satisfies the \textbf{Segal condition}. Let us emphasize that the Segal condition in simplicial degree zero amounts to saying that $\cM^\otimes_{[0]}$ is equivalent to the terminal category $\bbone$.

It turns out that the monoidal product $\otimes\colon\cM\times\cM\to\cM$ can be recovered up to equivalence from $p\colon\cM^\otimes\to\Delta\op$. But something seemingly more general works: Any Grothendieck opfibration $p\colon\cC\to\Delta\op$ satisfying the Segal condition defines a monoidal structure on the fiber $\cM=\cC_{[1]}$. We content ourselves by sketching a proof of this result. As usual, let $d_1\colon[2]\to[1]$ denote the face map in $\Delta^{op}$ which is opposite to the coface map $d^1\colon [1]\to [2].$ \emph{Choosing} an inverse of the equivalence $\sigma\colon\cC_{[2]}\to\cM\times\cM$ given by one of the Segal maps we can define a functor
\[
\otimes\colon\cM\times\cM\stackrel{\simeq}{\ot}\cC_{[2]}\stackrel{(d_1)_!}{\to}\cC_{[1]}=\cM.
\] 
In order to construct an associativity constraint for $\otimes$ we will invoke the cosimplicial identity $d^2\circ d^1=d^1\circ d^1\colon [1]\to[3].$ As a special instance of \eqref{eq:coherence} we thus obtain a natural isomorphism 
\[
(d_1)_!\circ(d_2)_!\cong(d_1)_!\circ(d_1)_!\colon\cC_{[3]}\to\cC_{[1]}=\cM.
\]
It is an instructive exercise to use the Segal condition to translate this into a natural isomorphism 
\[
\alpha\colon M_1\otimes (M_2\otimes M_3)\cong (M_1\otimes M_2)\otimes M_3,\qquad M_1,M_2,M_3\in \cM.
\]
With slightly more effort, one can use the different factorizations of the map $\{0,4\}\colon [1]\to[4]$ to deduce that this associativity constraint satisfies MacLane's pentagon axiom; see \eqref{eq:monoidal-intro}. And finally, in a similar way one checks that a monoidal unit $\lS\in\cM$ and the remaining coherence isomorphisms are also encoded by the Grothendieck opfibration $p\colon\cC\to\Delta\op$ satisfying the Segal condition.
\end{eg}

The point of this lengthy example was to show that there is an equivalent way of encoding monoidal structures. Instead of making a specific choice of a monoidal pairing, a monoidal unit, and coherence isomorphisms one can consider one `global object' which nicely hides all this structure, namely a Grothendieck opfibration $p\colon\cC\to\Delta\op$ satisfying the Segal condition. In particular, one does not have to make precise the coherence axioms of a monoidal category since they will follow from the (co)simplicial identities. (In \S\ref{subsec:algebra} we will briefly mention how monoidal functors and monoid objects fit into this opfibration picture.)

If one only cares about ordinary monoidal categories it might be arguable how large the benefit is by changing the classical perspective on monoidal categories to the Grothendieck opfibration picture. However, as mentioned in the introduction, in the $\infty$-categorical setting it allows us to avoid the complicated combinatorics of the Stasheff associahedra.

\subsection{Monoidal $\infty$-categories via coCartesian fibrations}
\label{subsec:monoidaloo}

We now turn to $\infty$-categorical variants of the above concepts. The main reference for the remainder of \S\ref{sec:monoidal} is Lurie's second book \cite{HA}. For the convenience of the reader we also include references to the former volumes \cite{DAGII,DAGIII} which are now subsumed in \cite{HA} but which might be a bit more accessible as a first reference. 

To begin with, we extend the notion of $p$-coCartesian morphisms to the context of $\infty$-categories. For this purpose it is handy to observe that the following is true (which is a refined version of the statement that \eqref{eq:coCartesian} is a pullback square).

\begin{lem}\label{lem:coCartesian}
Let $p\colon\cC\to\cD$ be a functor between ordinary categories. A morphism $f\colon c_1\to c_2$ in $\cC$ is $p$-coCartesian if and only if the following functor is an isomorphism
\[
\cC_{f/}\to \cC_{c_1/}\underset{\cD_{p(c_1)/}}{\times} \cD_{p(f)/}.
\]
\end{lem}

In \S\ref{subsec:slice} we already introduced $\infty$-categorical versions of slice categories. With this lemma at hand, there is the following formal generalization of $p$-coCartesian arrows to the $\infty$-categorical setting \cite[Definition~2.4.1.1]{HTT}.

\begin{defn}\label{defn:oocoCartesian}
Let $p\colon\cC\to\cD$ be a functor between $\infty$-categories. A morphism $f\colon c_1\to c_2$ in $\cC$ is \textbf{$p$-coCartesian} or a \textbf{$p$-coCartesian lift} of $\alpha=p(f)$ if the following map is an acyclic Kan fibration
\[
\cC_{f/}\to\cC_{c_1/}\underset{\cD_{p(c_1)/}}{\times}\cD_{p(f)/}.
\]
\end{defn}

The $\infty$-categorical concept corresponding to Grothendieck opfibrations is that of a \emph{coCartesian fibration} (Joyal \cite{joyal:I-II} uses the terminology \emph{Grothendieck opfibrations} instead). The main idea is again to axiomatically ask for a sufficient supply of coCartesian morphisms (see \cite[Definition~2.4.2.1]{HTT}).

\begin{defn}
A functor $p\colon\cC\to\cD$ between $\infty$-categories is a \textbf{coCartesian fibration} if the following two properties are satisfied.
\begin{enumerate}
\item The functor $p$ is an inner fibration (\autoref{defn:inner-fibration}).
\item For every object $c _1\in\cC$ and every morphism $\alpha\colon p(c_1)=d_1\to d_2$ in $\cD$, there is a $p$-coCartesian lift $f\colon c_1\to c_2$ of $\alpha.$
\end{enumerate}
\end{defn}

Similar to the case of categories, in this $\infty$-categorical context the fiber $\cC_d$ of a functor $\cC\to\cD$ is defined via a pullback square,
\[
\xymatrix@=1.5pc{
\cC_d\pullbackcorner\ar[r]\ar[d]&\cC\ar[d]^-p\\
\Delta^0\ar[r]_-d&\cD.
}
\]
Lurie proves that a coCartesian fibration gives rise to a covariantly depending family of $\infty$-categories. The fact that the fibers are \emph{$\infty$-categories} is immediate since a coCartesian fibration is an inner fibration and inner fibrations are stable under pullbacks. The hard part is to show that they assemble to a functor; see \autoref{per:Grothendieck}. In particular, any map $\alpha\colon d_1\to d_2$ in $\cD$ induces an essentially unique functor $\alpha_!\colon\cC_{d_1}\to\cC_{d_2}$ defined by means of coCartesian lifts. 

We saw in \S\ref{subsec:monoidal} that monoidal categories can be alternatively encoded by Grothen-dieck opfibrations satisfying the Segal condition. Having introduced the corresponding notion of coCartesian fibrations, we now turn this into a definition of monoidal $\infty$-categories \cite[Definition~1.1.2]{DAGII}.

\begin{defn}\label{defn:monoidaloo}
A \textbf{monoidal $\infty$-category} is a coCartesian fibration $p\colon\cM^\otimes\to N(\Delta\op)$ such that the Segal maps are equivalences, 
\[
\cM^\otimes_{[n]}\stackrel{\sim}{\to}(\cM^\otimes _{[1]})^{\times n},\qquad n\geq 0.
\]
\end{defn}

For simplicity, we refer to the $\infty$-category $\cM=\cM^\otimes _{[1]}$ as a monoidal $\infty$-category. If one wants to be very precise, then one should call the coCartesian fibration $p\colon\cM^\otimes\to N(\Delta\op)$ a \textbf{monoidal structure} on $\cM.$

The interpretation of such a coCartesian fibration $p\colon\cM^\otimes\to N(\Delta\op)$ is now similar to the situation in classical category theory. To give an example we just make the following remark. One immediate consequence of the axioms is that the fiber $\cM^\otimes_{[0]}$ is a contractible space. The unique map $s^0\colon [1]\to[0]$ in $\Delta$ induces a functor 
\[
\eta=(s_0)_!\colon\cM^\otimes_{[0]}\to\cM^\otimes_{[1]},
\]
and we call any object in its image a \textbf{monoidal unit} of $\cM.$ We leave it to the reader to justify this terminology in the classical case. 

Given a monoidal $\infty$-category $\cC$ it can be shown that the homotopy category $\Ho(\cC)$ inherits a monoidal structure. This underlines the idea that the monoidal structure on $\cC$ is associative and unital up to homotopy. Let us however emphasize that \autoref{defn:monoidaloo} really encodes much more structure, namely that of a monoidal product which is associative and unital up to \emph{coherent homotopy} (see \autoref{per:monoids} for a short discussion in the symmetric monoidal context). 

We now turn to important examples of monoidal $\infty$-categories (but see also~\S\ref{sec:stable}). 

\begin{egs}\label{egs:monoidal}
\begin{enumerate}
\item Let $\cM$ be a monoidal category and $p\colon\cM^\otimes\to\Delta\op$ be the associated Grothendieck opfibration. An application of the nerve functor yields a monoidal $\infty$-category
\[
N(p)\colon N(\cM^\otimes)\to N(\Delta\op).
\]
\item A special case of a monoidal category is given by the terminal category $\bbone.$ The associated Grothendieck construction can be identified with the identity functor $\Delta\op\to\Delta\op,$ and we obtain a monoidal $\infty$-category $N(\Delta\op)\to N(\Delta\op).$ This examples is of some interest in the study of algebra objects (see \S\ref{subsec:algebra}).
\item Again, first `honest examples' are obtained from suitable model categorical input, more precisely from suitably compatibly closed monoidal and simplicial model categories, by passing to coherent nerves. Such a category comes with three additional structures, namely a simplicial enrichment, a monoidal structure, and a model structure, which have to be suitably compatible. The assumptions made by Lurie in \cite[Prop.~1.6.5]{DAGII} are the following ones.
\begin{enumerate}
\item The closed monoidal structure is compatible with the enrichment in that $\otimes$ and the adjunction expressing the closedness are simplicial.
\item The monoidal pairing $\otimes\colon\cM\times\cM\to\cM$ is a \emph{left Quillen bifunctor}. 
\item The monoidal unit $\lS\in\cM$ is cofibrant.
\end{enumerate}
Under these assumptions, one can form a simplicial version of the category $\cM^\otimes$ from \S\ref{subsec:monoidal}. More precisely, given two finite strings $(M_1,\dots,M_n)$ and $(L_1,\dots,L_k)$ of objects in $\cM,$ the corresponding simplicial mapping space in $\cM^\otimes$ is given by
\[
\coprod_{\alpha\colon[k]\to [n]}\;\prod_{i=1}^k\Map_{\cM}(M_{\alpha(i-1)+1}\otimes\ldots\otimes M_{\alpha(i)},L_i).
\]
This simplicial category comes with an obvious simplicial functor $\cM^\otimes\to\Delta\op$ (regarding $\Delta\op$ as a discrete simplicial category). In order to obtain an $\infty$-category we consider the full simplicial subcategory 
\[
\cM^\otimes_\ncf\subseteq \cM^\otimes
\]
spanned by the finite strings of fibrant and cofibrant objects. It is then a consequence of the above compatibility assumptions that $\cM^\otimes_\ncf$ is a locally fibrant simplicial category so that $N_\Delta(\cM^\otimes_\ncf)$ is an $\infty$-category (\autoref{cor:fibrantsimplicial}). In this situation, Lurie establishes as \cite[Prop.~1.6.5]{DAGII} that
\[
N_\Delta(\cM^\otimes_\ncf)\to N_\Delta(\Delta\op)=N(\Delta\op)
\]
endows the $\infty$-category $N_\Delta(\cM_\ncf)$ with a monoidal structure. 
\item In ordinary category theory, important monoidal structures are given by categorical products and coproducts. These monoidal structures are referred to as the \textbf{Cartesian} and \textbf{coCartesian} monoidal structures, respectively. Given an $\infty$-category $\cC$ with finite (co)products, then one can construct (co)Cartesian monoidal structures on $\cC,$ given by coCartesian fibrations
\[
\cC^\times\to N(\Delta\op)\qquad\text{and}\qquad\cC^\sqcup\to N(\Delta\op),
\]
respectively. See \cite[\S2.4]{HA} for more details.
\end{enumerate}
\end{egs}

\begin{per}\label{per:Grothendieck}
CoCartesian fibrations are an $\infty$-categorical analogue of Grothen- dieck opfibrations. We saw in \S\ref{subsec:monoidal} that such an opfibration $p\colon\cC\to\cD$ encodes the idea of having a family of categories $\cC_d,d\in\cD,$ which depends covariantly on the object $d\in\cD.$ More precisely, by choosing certain $p$-coCartesian lifts we obtain a pseudo-functor 
\[
F_p\colon\cD\to\cCAT\colon d\mapsto\cC_d.
\] 

There is also a construction in the converse direction called the \textbf{Grothendieck construction}. Given a pseudo-functor $F\colon\cD\to\cCAT,$ one can form a new category~$\cE(F)$ which is defined as follows.
\begin{enumerate}
\item An object in $\cE(F)$ is a pair $(d,x)$ consisting of an object $d\in\cD$ and an object $x\in F(d).$
\item A morphism $(d,x)\to (d',x')$ in $\cE(F)$ is a pair $(\alpha,f)$ consisting of a morphism $\alpha\colon d\to d'$ in $\cD$ and a morphism $f\colon F(\alpha)(x)\to x'$ in $F(d').$
\item Compositions and identities are defined in the obvious way.
\end{enumerate}
The category $\cE(F)$ comes with a forgetful functor 
\[
p_F\colon\cE(F)\to\cD
\]
which projects objects and morphisms onto their respective first components, and one checks that $p_F$ is a Grothendieck opfibration. In fact, given an object $(d,x)$ in $\cE(F)$ and a morphism $\alpha\colon p_F(d,x)=d\to d'$ in $\cD,$ then a $p_F$-coCartesian lift is given by
\[
(\alpha,\id_{F(\alpha)(x)})\colon (d,x)\to (d',F(\alpha)(x)).
\]

It turns out that if we fix a category~$\cD$, then these two constructions tell us that category-valued pseudo-functors defined on $\cD$ and Grothendieck opfibrations with target $\cD$ are essentially the same. References for this theory include \cite[\S8]{borceux:2} and \cite{vistoli:descent}.

In \cite[\S3]{HTT}, Lurie has generalized these constructions to the $\infty$-categorical setting. Roughly speaking, he has established a result saying that a coCartesian fibration $p\colon\cC\to\cD$ is equivalent to giving a functor $\cD\to\cCat_\infty$, where $\cCat_\infty$ is the $\infty$-category of $\infty$-categories introduced in \autoref{per:inftyinfty}. In fact, there is a Quillen equivalence between certain model structures making this idea precise. Let us only mention that the homotopy theory of coCartesian fibrations above $\cD$ is encoded by the \textbf{coCartesian model structure} on $\sSet^+_{/\cD},$ the category of marked simplicial sets above $\cD=(\cD,\cD_1).$ An object $p\colon\cC\to\cD$ in this model structure is fibrant if and only if $p$ is a coCartesian fibration and if the marked edges in $\cC$ are precisely the $p$-coCartesian arrows.

In ordinary category theory, there are variants to these notions for contravariant category-valued pseudo-functors. The $\infty$-categorical analogue is that of a \textbf{Cartesian fibration}. Moreover, in the classical context --- in particular, in the theory of stacks in algebraic geometry or algebraic topology --- one frequently considers \emph{groupoid-valued} pseudo-functors (of either variance). These notions are sometimes also referred to as categories (co)fibered in groupoids. The $\infty$-categorical analogues of these are \textbf{left} and \textbf{right fibrations} and again there are suitable associated model categories in the background. For more details we refer to \cite[\S3]{HTT}.
\end{per}

\subsection{Algebra objects and monoidal functors}
\label{subsec:algebra}

Let $\cM$ be an ordinary monoidal category with monoidal pairing $\otimes$ and unit object $\lS \in\cM.$ An \textbf{algebra} or \textbf{monoid} in $\cM$ is an object $A\in\cM$ together with a multiplication map $\mu\colon A\otimes A\to A$ and a unit map $\eta\colon \lS\to A$ which satisfy obvious associativity and unitality conditions. We also say that $(\mu,\eta)$ specifies an \textbf{algebra structure} on~$A.$ 

In \S\ref{subsec:monoidal} we saw how to encode monoidal structures by means of Grothendieck opfibrations $p \colon\cM^\otimes\to\Delta\op$ satisfying the Segal condition. We now want to describe algebra structures in that picture. As a first guess for a definition of algebra objects we might consider sections of $p\colon\cM^\otimes\to\Delta\op$. Given an arbitrary such section 
\[
A \colon\Delta\op\to\cM^\otimes,
\]
the equivalences given by the Segal maps \eqref{eq:projection} imply that the values $A_{[n]}\in\cM^\otimes_{[n]}$ are determined by $n$ objects of $\cM=\cM^\otimes_{[1]}$,
\begin{equation}\label{eq:algebra}
\cM_{[n]}\ni A_{[n]}\quad\longleftrightarrow\quad A_{[n]}^1,\dots ,A_{[n]}^n\;\in \cM.
\end{equation}
Let us again consider the face map $d_1\colon [2]\to [1]$ which we saw in \S\ref{subsec:monoidal} to encode the monoidal pairing $\otimes\colon\cM\times\cM\to\cM$. The section $A$ evaluated on $d_1$ gives us by means of the identification~\eqref{eq:algebra} a map
\begin{equation}\label{eq:algebra-II}
(A_{[2]}^1,A_{[2]}^2)\to A_{[1]}^1.
\end{equation}
Moreover, by the discussion in \S\ref{subsec:monoidal}, we have a description of the $p$-coCartesian arrows for the Grothendieck opfibration $p\colon\cM^\otimes\to\Delta\op.$ Applied to our situation, a $p$-coCartesian lift of $d_1$ with domain $A_{[2]}$ corresponds under \eqref{eq:algebra} to a morphism
\[
(A_{[2]}^1,A_{[2]}^2)\to A_{[2]}^1\otimes A_{[2]}^2.
\]
The universal property of this $p$-coCartesian lift (as expressed by \eqref{eq:coCartesian} being a pullback square) implies that \eqref{eq:algebra-II} factors uniquely over this lift and we hence obtain an induced map
\begin{equation}
A_{[2]}^1\otimes A_{[2]}^2\to A_{[1]}^1.\label{eq:multiplication}
\end{equation}
If we now want to obtain a classical algebra object, we would like \eqref{eq:multiplication} to be a map of the form $M \otimes M\to M$ for some fixed $M\in\cM.$ Thus we should ensure that, among other things, the objects $A_{[2]}^1,A_{[2]}^2,$ and $A_{[1]}^1$ are isomorphic.

\begin{defn}
A morphism $\alpha\colon [n]\to[k]$ in $\Delta$ is \textbf{convex} if it is injective and if the image $\im(\alpha)\subseteq[k]$ is convex, i.e., the image is given by the interval $[\alpha(0),\alpha(n)].$
\end{defn}

It follows from the construction of the $p$-coCartesian arrows of $p\colon\cM^\otimes\to\Delta\op$ (compare to the discussion around \eqref{eq:coCartesianlift}) that the $p$-coCartesian lifts of convex maps $\alpha\colon [n]\to[k]$ in $\Delta$ induce projection functors $\cM^{\times k}\to\cM^{\times n}.$ In particular, $p$-coCartesian lifts defining the functors $(\iota_i)_!, i=1,2,$ as in \eqref{eq:projection} can be identified with the maps
\[
(A_{[2]}^1,A_{[2]}^2)\to A_{[2]}^1\qquad\text{and}\qquad(A_{[2]}^1,A_{[2]}^2)\to A_{[2]}^2,
\]
respectively. If the images of $\iota_i\colon [2]\to [1],i=1,2,$ under $A$ are $p$-coCartesian arrows, then by \autoref{lem:uniquecoC} we obtain the desired isomorphisms 
\[
A_{[2]}^1\cong A_{[1]}^1\qquad\text{and}\qquad A_{[2]}^2\cong A_{[1]}^1.
\]
With this preparation one establishes the following result.

\begin{prop}
Let $p\colon\cM^\otimes\to\Delta\op$ be a monoidal structure on $\cM=\cM^\otimes_{[1]}.$ Then a section $A\colon\Delta\op\to\cM^\otimes$ of $p$ which sends convex arrows to $p$-coCartesian arrows encodes an algebra structure on $A_{[1]}\in\cM.$ Conversely, any algebra object in $\cM$ determines such a section of $p\colon\cM^\otimes\to\Delta\op.$
\end{prop}

In the world of $\infty$-categories, we turn this observation into a definition (\cite[Definition~1.1.14]{DAGII}).

\begin{defn}
Let $p\colon\cM^\otimes\to N(\Delta\op)$ be a monoidal $\infty$-category. A section $A\colon N(\Delta\op)\to\cM^\otimes$ of $p$ is an \textbf{(associative) algebra object} in $\cM^\otimes$ if $A$ sends convex morphisms to $p$-coCartesian arrows in $\cM^\otimes.$
\end{defn}

Note that it is already a certain abuse of language to speak of algebra objects of $\cM^\otimes$ since the notion of algebra object obviously also depends on the coCartesian fibration. A further comfortable abuse of language is to simply speak of algebra objects in $\cM.$ 

Similar to the case of monoidal structures on $\infty$-categories themselves, an algebra object encodes quite a lot of structure: namely, given an algebra object $A$ in $\cM$, the underlying object $A_{[1]}$ is endowed with a multiplication map which is associative and unital up to \emph{coherent homotopy} (see \autoref{per:monoids} for an explanation of this similarity in the commutative case). In particular, an algebra object in a monoidal $\infty$-category defines an ordinary algebra object in the underlying homotopy category, but not conversely.

Algebra objects in monoidal $\infty$-categories are special cases of lax monoidal functors between monoidal $\infty$-categories. We include the following definition and leave it to the reader to check that in the case of ordinary categories this reduces to the usual concepts.

\begin{defn} 
Let $p\colon\cM^\otimes\to N(\Delta\op)$ and $q\colon\cN^\otimes\to N(\Delta\op)$ be monoidal $\infty$-categories. A \textbf{lax monoidal functor} $F\colon\cM^\otimes\to\cN^\otimes$ is a functor over $N(\Delta\op),$
\[
\xymatrix@=1.5pc{
\cM^\otimes\ar[rr]^-F\ar[dr]_-p\ar@{}[rrd]|{=}&&\cN^\otimes\ar[dl]^-q\\
&N(\Delta\op),&
}
\]
which sends $p$-coCartesian lifts of convex morphisms in $N(\Delta\op)$ to $q$-coCartesian arrows. A \textbf{monoidal functor} $F\colon\cM^\otimes\to\cN^\otimes$ is a functor over $N(\Delta\op)$ which sends arbitrary $p$-coCartesian arrows to $q$-coCartesian ones.
\end{defn}

Of course we also want to consider monoidal transformations between monoidal functors. More generally, (lax) monoidal functors between monoidal $\infty$-categories $\cM^\otimes$ and $\cN^\otimes$ are organized into $\infty$-categories 
\[
\nFun^{\otimes,\nlax}(\cM^\otimes,\cN^\otimes)\qquad\text{and}\qquad\nFun^\otimes(\cM^\otimes,\cN^\otimes),
\]
respectively. The $\infty$-category $\nFun^{\otimes,\nlax}(\cM^\otimes,\cN^\otimes)$ is the full subcategory of the $\infty$-category $\Map_{N(\Delta\op)}(\cM^\otimes,\cN^\otimes)$ of functors over $N(\Delta\op)$ spanned by the lax monoidal functors. Here, $\Map_{N(\Delta\op)}(\cM^\otimes,\cN^\otimes)$ is of course defined as the pullback
\begin{equation}
\vcenter{
\label{eq:maps-over}
\xymatrix@=1.5pc{
\Map_{N(\Delta\op)}(\cM^\otimes,\cN^\otimes)\pullbackcorner\ar[r]\ar[d]&\Map(\cM^\otimes,\cN^\otimes)\ar[d]^-{q_\ast}\\
\Delta^0\ar[r]_-p&\Map(\cM^\otimes,N(\Delta\op)),
}
}
\end{equation}
which is an $\infty$-category because with $q$ also $q_\ast$ is an inner fibration \cite[Corollary~2.3.2.5]{HTT}. Similarly, $\nFun^\otimes(\cM^\otimes,\cN^\otimes)\subseteq\nFun^{\otimes,\nlax}(\cM^\otimes,\cN^\otimes)$ is the full subcategory spanned by the monoidal functors.

As a special case, we see that algebra objects in a monoidal $\infty$-category $\cM^\otimes$ are themselves organized in an $\infty$-category. In fact, the \textbf{$\infty$-category $\AlgA(\cM^\otimes)$ of algebra objects} in $\mathcal{M}^\otimes$ can be defined as
\begin{equation}
\AlgA(\cM^\otimes)=\nFun^{\otimes,\nlax}(N(\Delta\op),\cM^\otimes),\label{eq:algebraoo}
\end{equation}
where $N(\Delta\op)\to N(\Delta\op)$ is the trivial monoidal structure as in \autoref{egs:monoidal}. We will say a bit more about similar $\infty$-categories in the context of symmetric monoidal structures in \S\ref{subsec:symmetric}.

\subsection{Symmetric monoidal $\infty$-categories}
\label{subsec:symmetric}

In order to obtain a theory of \emph{symmetric} monoidal $\infty$-categories one combines the Segal perspective on $\Eoo$-monoids and the $\infty$-categorical Grothendieck construction, and this section is hence inspired by the theory of `special $\Gamma$-spaces' (see \cite{segal:categories, schwede:gamma}). A good part of this section simply amounts to translating parts of \S\S\ref{subsec:monoidal}-\ref{subsec:algebra} to the context of symmetric monoidal $\infty$-categories, and we are hence rather sketchy. 

We again begin with the classical situation in ordinary category theory. In~\S\ref{subsec:monoidal}, we described monoidal categories in terms of Grothendieck opfibrations $\cM^\otimes\to\Delta\op.$ In that picture, the monoidal product was encoded by the induced functor
\[
(d_1)_!\colon\cM^\otimes _{[2]}\to \cM^\otimes_{[1]}=\cM.
\] 
If we want to have a similar description of \emph{symmetric} monoidal categories, we must be able to encode symmetry isomorphisms, thus we have to encode the flip map
\[
t\colon\cM^{\times 2}\to\cM^{\times 2}\colon (X,Y)\mapsto (Y,X)
\]
and, more generally, any permutation of $n$ objects in $\cM.$ It is hence plausible that the role of $\Delta$ is taken by `a category of finite sets with all maps between them'. The details are as follows. 

For a natural number $n\geq 0,$ let $\la n\ra$ be the finite pointed set 
\[
\la n \ra=\{0<1<\ldots <n\}
\] 
with $0\in\la n\ra$ as base point. The category $\Fin$ is the full subcategory of the category of pointed sets spanned by the objects $\la n\ra, n\geq 0$. Note that the natural ordering on $\la n \ra$ does not play a role in the definition of $\Fin$ but it will have its uses in the formation of higher monoidal products. We denote by
\begin{equation}\label{eq:rho}
\rho ^j\colon\la n\ra \to \la 1\ra,\quad n\geq 1,\quad j=1,\ldots, n,
\end{equation}
the unique pointed map $\la n\ra \to \la 1 \ra$ with $(\rho^j)^{-1}(1)=\{j\}.$

Let now $\cM$ be a \emph{symmetric} monoidal category with monoidal product $\otimes$ and monoidal unit $\lS\in\cM.$ Following a pattern similar to \S\ref{subsec:monoidal}, we construct a new category $\cM^\otimes$ as follows. An object in $\cM^\otimes$ is a finite (possibly empty) sequence of objects in $\cM,$
\[
(M_1,\dots ,M_n),\quad M_i\in\cM,\quad n\geq 0.
\]
A morphism $(M_1,\ldots,M_n)\to (L_1,\ldots,L_k)$ between two such sequences is a pair $(\alpha ,\{f_i\}_i)$ consisting of a morphism $\alpha\colon \la n\ra\to \la k\ra$ in $\Fin$ together with morphisms
\[
f_i\colon \bigotimes_{j\in \alpha^{-1}(i)}M_j\to L_i,\quad i=1,\dots,k,
\]
where the tensor product is formed according to the ordering on $\alpha ^{-1}(i).$ Again, if the set $\alpha ^{-1}(i)$ is empty, then this is to be read as a map $f_i\colon \lS\to L_i.$ The composition is obtained from the compositions in $\cM$ and $\Fin$ together with the coherence constraints of $\cM.$ There is an obvious projection functor $p\colon\cM^\otimes\to\Fin$ given by $(M_1,\ldots ,M_n)\mapsto \la n \ra$ and $(\alpha,\{f_i\}_i)\mapsto\alpha$.

\begin{prop}
If $\cM$ is a symmetric monoidal category, then the functor $p \colon\cM^\otimes\to\Fin$ is a Grothendieck opfibration. Moreover, this functor satisfies the \textbf{Segal condition}, i.e., the Segal maps 
\[
(\rho^1_!,\ldots,\rho^n_!)\colon\cM^\otimes_{\la n\ra}\to\cM^{\times n},\quad n\geq 0,
\]
are equivalences. Conversely, any Grothendieck opfibration $p\colon\cC\to\Fin$ satisfying the Segal condition encodes a symmetric monoidal structure on $\cM=\cC_{\la 1\ra}.$
\end{prop}

Here is a sketch of a proof. Given $(M_1,\ldots ,M_n)\in\cM^\otimes_{\la n\ra}$ and $\alpha\colon\la n\ra\to\la k\ra$, an associated $p$-coCartesian is obtained by choosing isomorphisms
\[
f_i\colon\bigotimes_{j\in \alpha^{-1}(i)}M_j\to L_i,\quad i=1,\dots,k.
\]
In the special case of $\rho^j$ as in \eqref{eq:rho}, it follows that such lifts are of the form
\[
(\rho^j,\id_{M_j})\colon (M_1,\ldots,M_n)\to M_j,
\]
and the associated functor $\rho^j_!\colon\cM^\otimes_{\la n\ra}\to \cM^\otimes_{\la 1\ra}$ can hence be identified with a projection functor, implying the Segal condition. 

Conversely, let us consider a Grothendieck opfibration $p\colon\cC\to\Fin$ satisfying the Segal condition. We content ourselves by mentioning the following three steps towards the construction of a symmetric monoidal structure on $\cM=\cC_{\la 1\ra}$. 
\begin{enumerate}
\item Let $m\colon\la 2\ra\to\la 1\ra$ be the map in $\Fin$ determined by $m(1)=m(2)=1.$ By means of the Segal condition, we can then define a functor
\[
\otimes\colon\cM\times\cM\stackrel{\sim}{\ot}\cC_{\la 2\ra}\stackrel{m_!}{\to}\cC_{\la 1\ra}=\cM,\]
which will be the monoidal pairing. 
\item As a special case of the Segal condition we obtain an equivalence $\cC_{\la 0\ra}\simeq\bbone$. Thus, for the unique map $n\colon\la 0\ra\to\la 1\ra$ in $\Fin,$ the induced functor
\[
n_!\colon\cC_{\la 0\ra}\to\cC_{\la 1\ra}=\cM
\]
essentially classifies an object in $\cM$, which will be the monoidal unit~$\lS$.
\item The \textbf{twist map} $t\colon\la 2\ra\to\la 2\ra$ in $\Fin$ is the automorphism that interchanges $1$ and $2$. The equality $m=m\circ t\colon\la 2\ra\to\la 1\ra$ together with \autoref{lem:uniquecoC} yields a unique natural isomorphism
\[
\sigma \colon m_!\cong m_!\circ t_!\colon\cC_{\la 2\ra}\to\cC_{\la 1\ra}=\cM
\]
over $\la 1\ra$, which can be shown to induce a symmetry constraint for $\otimes.$
\end{enumerate}

By similar arguments, one obtains associativity and unitality constraints and establishes the coherence axioms. As in the non-symmetric case, we now turn this observation into a definition \cite{DAGIII-old}.

\begin{defn}
A \textbf{symmetric monoidal $\infty$-category} is a coCartesian fibration $p\colon\cM^\otimes\to N (\Fin)$ such that the Segal maps are equivalences,
\[
(\rho^1_!,\ldots,\rho^n_!)\colon\cM^\otimes_{\la n\ra}\stackrel{\sim}{\to}(\cM^\otimes_{\la 1\ra})^{\times n},\qquad n\geq 0.
\]
\end{defn}

We include a few remarks which are parallel to corresponding statements in \S\ref{subsec:monoidaloo}.

\begin{rmk}
\begin{enumerate}
\item A symmetric monoidal $\infty$-category $p\colon\cM^\otimes\to N(\Fin)$ endows the underlying $\infty$-category $\cM=\cM^\otimes_{\la 1\ra}$ with a monoidal pairing which is associative and commutative \emph{up to coherent homotopies}; see \autoref{per:monoids}. In particular, the homotopy category of a symmetric monoidal $\infty$-category is canonically a symmetric monoidal category.
\item  Applying the nerve construction to Grothendieck opfibrations associated to ordinary symmetric monoidal categories we obtain symmetric monoidal $\infty$-categories. In particular, the identity functor $N(\Fin)\to N(\Fin)$ is a symmetric monoidal $\infty$-category. As a variant, given a closed, symmetric monoidal, simplicial model category satisfying suitable compatibility assumptions, by means of the coherent nerve construction we obtain `honest examples' of symmetric monoidal $\infty$-categories (\cite[\S8]{DAGIII-old}). 
\end{enumerate}
\end{rmk}

Before we turn to algebra objects in the symmetric monoidal context, we expand a bit on the relation between monoidal and symmetric monoidal $\infty$-categories. For this purpose, we consider the following functor 
\begin{equation}\label{eq:circle}
\phi \colon\Delta^{op}\to\Fin
\end{equation}
which on objects is given by $[n]\mapsto \la n\ra$. Given a morphism $\alpha \colon [k]\to [n]$ in $\Delta$, the induced pointed map $\phi (\alpha)\colon \la n\ra\to \la k\ra$ is defined by
$$\phi (\alpha)(j)=\;\left\{ \begin{array}{c@{\;}l}
				i\; & \quad \mbox{if there is an } i\;\mbox{such that}\;j\in[\alpha(i-1)+1,\alpha(i)],\\
				\ast\; & \quad \mbox{otherwise.}
				\end{array} \right.$$
Since $\alpha$ is monotone, such an $i$ is unique if it exists, and we leave it to the reader to check that this defines a functor. In fact, up to a restriction of the codomain, the functor $\phi$ is simply the simplicial circle $S^1\in\sSet_\ast$ defined as the coequalizer
\[
\Delta^0\rightrightarrows \Delta^1\to S^1
\]
and considered as a pointed simplicial set. As a special case we obtain for the image of the opposite $\iota_j$ of $\iota_{\{j-1,j\}}\colon[1]\to[n]$ that $\phi(\iota_j)=\rho^j$ with $\rho^j$ as in \eqref{eq:rho}.

The role of convex maps in the theory of monoidal $\infty$-categories is taken by \emph{inert} or \emph{collapsing} maps in the theory of symmetric monoidal $\infty$-categories.

\begin{defn}
A morphism $\alpha\colon\la n\ra\to\la k\ra$ in $\Fin$ is \textbf{inert} or \textbf{collapsing} if $\alpha^{-1}(i)$ is a singleton for every $1\leq i\leq k$.
\end{defn} 

We note that $\alpha \colon [k]\to [n]$ in $\Delta$ is convex if and only if $\phi (\alpha)\colon \la n\ra\to\la k\ra$ in $\Fin$ is inert, and that the maps \eqref{eq:rho} are inert. Given a symmetric monoidal $\infty$-category $p\colon\cM^\otimes\to N(\Fin)$ we already observed that the induced functors $\rho^j_!$ are projection functors. Similarly, general inert morphisms induce projection and permutation functors. This suggests the following definition.

\begin{defn}
Let $p \colon\cM^\otimes\to N(\Fin)$, $q \colon\cN^\otimes\to N ( \Fin)$ be symmetric monoidal $\infty$-categories and let $F\colon\mathcal{M}^\otimes \to \mathcal{N}^\otimes$ be a functor over $N (\Fin)$.
\begin{enumerate}
\item The functor $F$ is \textbf{symmetric monoidal} if it sends $p$-coCartesian arrows to $q$-coCartesian arrows. 
\item The functor $F$ is \textbf{lax symmetric monoidal} if it sends $p$-coCartesian lifts of inert morphisms to $q$-coCartesian arrows. 
\end{enumerate}
\end{defn}

Symmetric monoidal and lax symmetric monoidal functors respectively are organized in $\infty$-categories, namely the corresponding full subcategories
\[
\nFun^\otimes(\cM^\otimes,\cN^\otimes)\subseteq\nFun^{\otimes,\nlax}(\cM^\otimes,\cN^\otimes)\subseteq \Map_{N(\Fin)}(\cM^\otimes,\cN^\otimes)
\]
where $\Map_{N(\Fin)}(\cM^\otimes,\cN^\otimes)$ is defined in analogy to \eqref{eq:maps-over}. As a special case we obtain \textbf{$\infty$-categories of commutative algebra objects} $\AlgE(\cM^\otimes)$ defined by 
\begin{equation}\label{eq:comm-alg}
\AlgE(\cM^\otimes)=\nFun^{\otimes,\nlax}(N(\Fin),\cM^\otimes).
\end{equation}
More explicitly, a commutative algebra object is a section $E\colon N(\Fin)\to\cM^\otimes$ of $p\colon\cM^\otimes\to N(\Fin)$ sending inert morphisms to $p$-coCartesian ones. Such a section endows the underlying object $E_{\la 1\ra}\in\cM$ with a coherently associative and commutative pairing $E_{\la 1 \ra}\otimes E_{\la 1\ra}\to E_{\la 1\ra}$.

\begin{rmk}
Symmetric monoidal $\infty$-categories have underlying monoidal $\infty$-categories. In fact, given a symmetric monoidal $\infty$-category $\cM^\otimes\to N(\Fin)$, then the underlying monoidal $\infty$-category $U(\cM^\otimes)$ is defined as the pullback
\[
\xymatrix@=1.5pc{
U(\cM^\otimes)\ar[d]\ar[r]\pullbackcorner&\cM^\otimes\ar[d]^-p\\
N(\Delta\op)\ar[r]_-{N(\phi)}&N(\Fin),
}
\]
where $\phi$ is again the simplicial circle~\eqref{eq:circle}. Similarly, one shows that commutative algebra objects have underlying associative algebra objects.
\end{rmk}

\begin{per}\label{per:monoids}
In this section we saw that symmetric monoidal $\infty$-categories are endowed with a coherently associative and commutative monoidal structure and we just claimed that something similar is true for commutative algebra objects in the sense of \eqref{eq:comm-alg}. In this perspective we make the similarity between these two situations precise and refer the reader to \cite[\S2]{DAGIII-old} for details.

Given an $\infty$-category $\cC$ with finite products, there is the Cartesian monoidal structure $\cC^\times\to N(\Fin)$ and there are two associated $\infty$-categories.
\begin{enumerate}
\item As a special case of \eqref{eq:comm-alg} we have the $\infty$-category $\AlgE(\cC^\times)$ of commutative algebra objects in $\cC^\times.$
\item On the other hand, we can consider \textbf{commutative monoids} in $\cC$, i.e., functors $M\colon N(\Fin)\to\cC$ such that the Segal maps
\[
M_{\la n\ra}\stackrel{\sim}{\to} M_{\la 1\ra}\times\ldots\times M_{\la 1\ra}
\]
are equivalences. The $\infty$-category $\MonE(\cC)$ of commutative monoids in $\cC$ is the full subcategory of $\nFun(N(\Fin),\cC)$ spanned by the commutative monoids. 
\end{enumerate}
And it turns out that there is an equivalence $\AlgE(\cC^\times)\simeq\MonE(\cC)$ over~\cC.

We now apply this to the $\infty$-category $\cCat_\infty$ of $\infty$-categories (\autoref{per:inftyinfty}), which is an example of an $\infty$-category admitting finite products. The Grothendieck construction (\autoref{per:Grothendieck}) implies that $\nFun(N(\Fin),\cCat_\infty)$ is equivalent to the $\infty$-category of coCartesian fibrations $p\colon\cM\to N(\Fin)$. And under this equivalence commutative monoids in $\cCat_\infty$ and symmetric monoidal $\infty$-categories correspond to each other since both are defined by similar Segal conditions. Thus, as an upshot, we obtain an equivalence of $\infty$-categories
\[
\AlgE(\cCat_\infty^\times)\simeq\cCat_\infty^{\mathrm{sMon}},
\]
where $\cCat_\infty^{\mathrm{sMon}}$ denotes the $\infty$-category of symmetric monoidal $\infty$-categories and symmetric monoidal functors, explaining why in both cases we obtain similar coherence data. There is a similar equivalence in the case of monoidal $\infty$-categories and monoidal functors, namely
\[
\AlgA(\cCat_\infty^\times)\simeq\cCat_\infty^{\mathrm{Mon}}.
\]
\end{per}

Obviously, having introduced commutative algebra objects, one would now like to study modules over such algebra objects as well as the existence of limits and colimits in related $\infty$-categories. Such a theory exists and we refer to \cite{DAGII,DAGIII} and \cite{HA} for more details. 

We content ourselves by concluding this section with the following result concerning \emph{initial} objects in $\infty$-categories of commutative algebra objects. This result comes up again in the final section~\S\ref{sec:stable}.

\begin{prop}\label{prop:unit-algebra}
For every symmetric monoidal $\infty$-category $\cM^\otimes\to N(\Fin)$ the $\infty$-category $\AlgE (\cM)$ has an initial object. Moreover, a commutative algebra object $E$ is initial if and only if the unit map $\lS\to E_{\la1\ra}$ is an equivalence in $\cM$.
\end{prop}

\section{Stable $\infty$-categories and the universal property of spectra}
\label{sec:stable}

In this final section we give an introduction to stable $\infty$-categories. By definition a finitely complete and finitely cocomplete $\infty$-categories is stable if it admits a zero object and if a square in it is a pullback if and only if it is a pushout. Typical examples of stable $\infty$-categories arise in homological algebra ($\infty$-categories of chain complexes) and stable homotopy theory (the $\infty$-category of spectra). In fact, it turns out the $\infty$-category of spectra is the universal example of a stable $\infty$-category in a certain precise sense.

Similarly to stable model categories, stable $\infty$-categories are an enhancement of triangulated categories. In \S\ref{subsec:stable} we sketch some of ingredients involved in a proof that homotopy categories of stable $\infty$-categories can be turned into triangulated categories. In \S\ref{subsec:spectra} we briefly discuss the \emph{stabilization} of nice $\infty$-categories which is obtained by passing to internal spectrum objects. We conclude this subsection by a precise universal property of this stabilization process. Preparing the ground for the construction of the smash product, in \S\ref{subsec:tensor} we discuss the cocontinuous tensor product of presentable $\infty$-categories. Following Lurie, this allows us in \S\ref{subsec:smash} to give a very conceptual construction of the smash product monoidal structure on spectra and hence to define associative and commutative ring spectra. 

\subsection{Stable $\infty$-categories} 
\label{subsec:stable}

General references for the first two subsections are \cite{DAGI} and \cite[\S1]{HA}. As a first step we collect a few basics concerning \emph{pointed} $\infty$-categories.

\begin{defn}\label{defn:pointed}
An $\infty$-category is \textbf{pointed} if it admits a zero object, i.e., an object which is initial and final.
\end{defn}

Thus, an $\infty$-category~$\cC$ is pointed if there is an object $0\in\cC$ such that for all $x\in\cC$ the mapping spaces $\Map_{\cC}(x,0)$ and $\Map _{\cC}(0,x)$ are contractible. It follows, that for any two objects $x,y\in\cC$ there is a zero map 
\[
0=0_{x,y}\colon x\to y,
\] 
well-defined up to a contractible space of choices. Again by \autoref{prop:final}, if an $\infty$-category $\cC$ is pointed, then the full subcategory spanned by the zero objects is a contractible Kan complex. 

\begin{egs}
\begin{enumerate}
\item Let \cC be an $\infty$-category with a terminal object $\ast\in\cC$. The undercategory $\cC_\ast=\cC_{\ast/}$ (see \autoref{eg:undercategory}) is a pointed $\infty$-category, called the $\infty$-category of \textbf{pointed objects} in \cC. Adding a disjoint base point defines a functor ${}_+\colon \cC\to\cC_\ast$ which is left adjoint to the forgetful functor ${}_-\colon\cC_\ast\to\cC$,
\[
({}_+,{}_-)\colon\cC\rightleftarrows\cC_\ast.
\]
(For the notion of an adjunction between $\infty$-categories we again refer the reader to \cite[p.337]{HTT} and \cite[\S3]{DAGII}.) As a special case we obtain the \textbf{$\infty$-category~$\cS_\ast$ of pointed spaces} and the corresponding adjunction
\[
({}_+,{}_-)\colon\cS\rightleftarrows\cS_\ast.
\]
\item The underlying $\infty$-category of a pointed simplicial model category is pointed; see \autoref{egs:oocats}. An ordinary category is pointed if and only if the nerve is a pointed $\infty$-category.
\end{enumerate}
\end{egs}

The $\infty$-category~$\cS_\ast$ of pointed spaces together with the $0$-sphere $S^0\in\cS_\ast$ enjoys the following universal property (which is a pointed variant of~\autoref{cor:spaces}).

\begin{prop}\label{prop:pointed-spaces}
Let $\cD$ be a pointed, presentable $\infty$-category. Evaluation at the $0$-sphere $S^0\in\cS_\ast$ induces an equivalence of $\infty$-categories
\[
\FunL(\cS_\ast,\cD)\stackrel{\sim}{\to}\cD.
\]
\end{prop}

This result makes precise that $\cS_\ast$ is the \emph{free pointed, presentable $\infty$-category generated by~$S^0$}.

A \textbf{triangle} $\tau$ in a pointed $\infty$-category $\cC$ is a diagram $\tau\colon\square\to\cC$,
\begin{equation}\label{eq:triangle}
\vcenter{
\xymatrix@=1.5pc{
x\ar[r]^f\ar[d]\drtwocell<\omit>{<2>}\ar[dr] & y\ar[d]^g\\
0 \ar[r] & z,\ultwocell<\omit>{<2>}
}
}
\end{equation}
which vanishes at the lower left corner. Thus, a triangle in $\cC$ encodes two composable arrows $f\colon x\to y$ and $g\colon y\to z$, a further arrow $h\colon x\to z$  together with a homotopy $h\simeq g\circ f$ and a \emph{null-homotopy} $h\simeq 0$. Recall the definition of pullback and pushout squares in \autoref{defn:pull-push}.

\begin{defn}
A triangle in a pointed $\infty$-category is \textbf{exact} if it is a pullback square. Dually, a triangle is \textbf{coexact} if it is a pushout square.
\end{defn}

For every finitely complete, finitely cocomplete, and pointed $\infty$-category~$\cC$ we denote by
\[
\cC^\Sigma\subseteq\nFun(\square,\cC)
\]
the full subcategory spanned by the coexact triangles which also vanish on the upper right corner,
\[
\xymatrix@=1.5pc{
x\ar[r]\ar[d]& 0'\ar[d]\\
0\ar[r]&y. \pushoutcorner 
}
\]
There is a dually defined $\infty$-category  $\cC^\Omega\subseteq\nFun(\square,\cC)$ of exact triangles vanishing on the upper right corner. Morally, such diagrams should be determined by the value in the upper left corner in the first and by the value in the lower right corner in the second case. This is made precise by the following result.

\begin{prop}\label{prop:sigma-omega}
Let $\cC$ be a finitely complete, finitely cocomplete, and pointed $\infty$-category. The evaluation maps
\[
\nev_{(0,0)}\colon\cC^\Sigma\to\cC\qquad\text{and}\qquad\nev_{(1,1)}\colon\cC^\Omega\to\cC
\]
are acyclic Kan fibrations.
\end{prop}

This proposition and many similar results in this section are consequences of an $\infty$-categorical version of the calculus of Kan extensions (\autoref{per:derivator}). One of the key facts of constant use is \cite[Prop.~4.3.2.15]{HTT}. In this course we will not pursue this calculus any further, but we only consider results from this calculus which are `similarly plausible' as \autoref{prop:sigma-omega}. 

We briefly recall that given an acyclic Kan fibration $p\colon X\to Y$, then the space of sections $\Gamma(p)\in\sSet$ is a contractible Kan complex. In fact, for every simplicial set~$K$ the induced map $p_\ast\colon\Map(K,X)\to\Map(K,Y)$ between simplicial mapping spaces as defined by \eqref{eq:ssetenriched} is again an acyclic Kan fibrations. Since acyclic Kan fibrations are stable under pullbacks, we conclude that $\Gamma(p)$ is a contractible Kan complex by considering the defining pullback diagram
\begin{equation}\label{eq:setions}
\vcenter{
\xymatrix@=1.5pc{
\Gamma(p)\pullbackcorner\ar[r]\ar[d]&\Map(Y,X)\ar[d]^-{p_\ast}\\
\Delta^0\ar[r]_-{\id_Y}&\Map(Y,Y).
}
}
\end{equation}

Thus, under the assumption of \autoref{prop:sigma-omega}, we can choose sections
\[
s_\Sigma\colon\cC\to\cC^\Sigma\qquad\text{and}\qquad s_\Omega\colon\cC\to\cC^\Omega
\]
of the evaluation maps $\nev_{(0,0)}$ and $\nev_{(1,1)}$, respectively, and these sections are unique up to contractible spaces of choices. Consequently, the following is well-defined.

\begin{defn}\label{defn:susp}
Let $\cC$ be a finitely complete, finitely cocomplete, and pointed $\infty$-category. The \textbf{suspension functor} $\Sigma=\Sigma_{\cC}\colon\cC\to\cC$ and the \textbf{loop functor} $\Omega=\Omega_{\cC}\colon\cC\to\cC$ are respectively defined as
\[
\Sigma\colon\cC\xrightarrow{s_\Sigma}\cC^\Sigma\xrightarrow{\nev_{(1,1)}}\cC\qquad\text{and}\qquad
\Omega\colon\cC\xrightarrow{s_\Omega}\cC^\Omega\xrightarrow{\nev_{(0,0)}}\cC.
\]
\end{defn}

\begin{prop}
Let $\cC$ be a finitely complete, finitely cocomplete, and pointed $\infty$-category. The suspension functor $\Sigma\colon\cC\to\cC$ is left adjoint to the loop functor $\Omega\colon\cC\to\cC$,
\[
(\Sigma,\Omega)\colon\cC\rightleftarrows\cC.
\]
\end{prop}

In a similar way one defines \emph{cofibers} and \emph{fibers} in pointed $\infty$-categories. In fact, in the case of cofibers, starting with a morphism $f\colon x\to y$, suitable combinations of Kan extensions yield a coexact triangle as in \eqref{eq:triangle}. Combining these Kan extensions with a restriction of such triangles to the vertical morphism on the right yields a \textbf{cofiber} functor
$\cof\colon\nFun(\Delta^1,\cC)\to\nFun(\Delta^1,\cC).$ Dualizing this, we obtain a \textbf{fiber} functor $\fib\colon\nFun(\Delta^1,\cC)\to\nFun(\Delta^1,\cC).$

\begin{prop}
Let $\cC$ be a finitely complete, finitely cocomplete, and pointed $\infty$-category. The cofiber functor $\cof\colon\nFun(\Delta^1,\cC)\to\nFun(\Delta^1,\cC)$ is left adjoint to the fiber functor $\fib\colon\nFun(\Delta^1,\cC)\to\nFun(\Delta^1,\cC)$,
\[
(\cof,\fib)\colon\nFun(\Delta^1,\cC)\rightleftarrows\nFun(\Delta^1,\cC).
\]
\end{prop}

One way of defining \emph{stable} $\infty$-categories is as follows.

\begin{defn}\label{defn:stable}
A finitely complete, finitely cocomplete, and pointed $\infty$-category is \textbf{stable} if a triangle in it is exact if and only if it is coexact.
\end{defn}

This is one way of imposing the usual \emph{linearity condition} axiomatizing stability. It turns out that this definition admits the following reformulations.

\begin{thm}\label{thm:char-stable}
The following are equivalent for a finitely complete, finitely cocomplete, and pointed $\infty$-category~\cC.
\begin{enumerate}
\item The adjunction $(\Sigma,\Omega)\colon\cC\rightleftarrows\cC$ is an equivalence.
\item The adjunction $(\cof,\fib)\colon\nFun(\Delta^1,\cC)\rightleftarrows\nFun(\Delta^1,\cC)$ is an equivalence.
\item The $\infty$-category~\cC is stable.
\item A square in $\cC$ is a pullback square if and only if it is a pushout square.
\end{enumerate} 
\end{thm}

\begin{egs}\label{egs:stable}
\begin{enumerate}
\item The underlying $\infty$-category of a stable, simplicial model category (\autoref{egs:oocats}) is stable.
\item Let us recall that there are many Quillen equivalent simplicial model categories of spectra such that the homotopy categories are \emph{the} stable homotopy category $\mathcal{SHC}$ of Boardman; see \cite{vogt:boardman} and \cite[Part~III]{adams:stable} for classical accounts of this category and for example \cite{hss:symmetric,ekmm:rings,mmss:spectra} for good point-set categories of spectra. The underlying $\infty$-categories of these model categories are equivalent, and any of them will be denoted by $\nSp$ and referred to as the \textbf{$\infty$-category of spectra}. We will see in \S\ref{subsec:spectra} that there also is an intrinsic construction of this $\infty$-category.
\item In the context of homological algebra, there are stable \emph{derived $\infty$-categories}; see \cite[\S1]{HTT}.
\end{enumerate}
\end{egs}

Like stable model categories also stable $\infty$-categories provide an enhancement of triangulated categories. In particular, the homotopy category of a stable $\infty$-category can be endowed with a triangulation. To this end, we define a \textbf{cofiber sequence} in a pointed $\infty$-category~$\cC$ to be a diagram $\Delta^2\times\Delta^1\to\cC$ looking like
\begin{equation}\label{eq:cof-seq}
\vcenter{
\xymatrix@=1.5pc{
x\ar[r]^-f\ar[d]&y\ar[r]\ar[d]^-g&0'\ar[d]\\
0\ar[r]&z\ar[r]_-h\pushoutcorner&w\pushoutcorner
}
}
\end{equation}
and such that both squares are pushout squares. (A cofiber sequence is essentially obtained by two iterations of the passage to the cofiber of a morphism.) As in the case of ordinary category theory, it follows that also the composite square is a pushout square, and, by definition of the suspension functor (\autoref{defn:susp}), we obtain a canonical equivalence $\phi\colon w\simeq\Sigma x$. Thus, if we pass to homotopy classes of morphisms, then associated to \eqref{eq:cof-seq} we obtain by means of this equivalence an \emph{incoherent} cofiber sequence or \emph{triangle}
\[
\xymatrix{
T_f\colon&x\ar[r]^-f&y\ar[r]^-g&z\ar[r]^-{\phi\circ h}&\Sigma x,
}
\]
which is an ordinary diagram in the homotopy category $\Ho(\cC)$. If \cC is a stable $\infty$-category, then we say that a triangle in $\Ho(\cC)$ is \textbf{distinguished} if it is isomorphic to $T_f$ for some $f\colon\Delta^1\to\cC$. 

In the following result (\cite{DAGI}) we also denote by $\Sigma\colon\Ho(\cC)\to\Ho(\cC)$ the functor induced by $\Sigma\colon\cC\to\cC$. We assume that the reader is familiar with the notion of a triangulated category; see the original references of Puppe~\cite[Satz~3.5 and \S4.1]{puppe:stabil} or Verdier \cite{verdier:derived} (a reprint of his 1967 thesis) as well as the monographs~\cite{neeman:triangulated,hjr:triangulated}.

\begin{thm}\label{thm:triang}
Let $\cC$ be a stable $\infty$-category. The functor $\Sigma\colon\Ho(\cC)\to\Ho(\cC)$ and the above class of distinguished triangles endow the homotopy category $\Ho(\cC)$ with the structure of a triangulated category.
\end{thm}

A natural class of functors between stable $\infty$-categories is the class of exact functors in the sense of the following definition.

\begin{defn}\label{defn:exact}
Let \cC and \cD be finitely complete and finitely cocomplete $\infty$-categories.
\begin{enumerate}
\item A functor $\cC\to\cD$ is \textbf{left exact} if it preserves pullbacks and terminal objects.
\item A functor $\cC\to\cD$ is \textbf{right exact} if it preserves pushouts and initial objects.
\item A functor $\cC\to\cD$ is \textbf{exact} if it is left exact and right exact.
\end{enumerate}
\end{defn}

Clearly, limit-preserving functors (and hence, in particular, right adjoint functors) are left exact, and dually. If \cC and \cD are stable $\infty$-categories, then the three classes in \autoref{defn:exact} agree. It turns out that the triangulations of \autoref{thm:triang} are natural with respect to exact functors in the following sense. Given an exact functor $F\colon\cC\to\cD$ of stable $\infty$-categories, then the functor $F\colon\Ho(\cC)\to\Ho(\cD)$ can be endowed with the structure of an exact functor.

\begin{rmk}\label{rmk:stable-moral}
Note that, by the very definition, a stable $\infty$-category is obtained from the general notion of an $\infty$-category by imposing certain (easily motivated) exactness \emph{properties}; namely, we ask that finite limits and finite colimits exist and that certain limit type constructions are colimit type constructions and conversely. Similarly, the good notion of morphisms of stable $\infty$-categories, namely, exact functors, are defined by asking for the property that certain (co)limits are preserved.

This is in contrast to the more classical notion of a triangulated category which addresses the bad categorical properties of derived categories of abelian categories or of homotopy categories of stable model categories or stable $\infty$-categories by imposing additional \emph{structure}. Given an additive category, the axioms of a triangulated category ask for the existence of some non-canonical additional structure (the suspension functor and the class of distinguished triangles) satisfying certain properties. Similarly, given two triangulated categories $\cT$ and $\cT'$, a morphism should be an additive functor $F\colon\cT\to\cT'$ which sends distinguished triangles to distinguished triangles. In order to make this precise, one has to ask for the existence of an exact structure, i.e., a natural isomorphism $F\Sigma\cong\Sigma F$.

Despite their great successes in many areas of pure mathematics, it was obvious from the very beginning on (see for example already the introduction to \cite{heller:shc}) that the axioms of a triangulated category suffer certain defects (non-functoriality of the cone construction, no good theory of homotopy limits and homotopy colimits, diagram categories of triangulated categories do not admit canonical triangulations). 

There are more traditional attempts to improve the axioms of a triangulated category and the basic idea goes back at least to \cite{beilinson:perverse}. The idea is to ask for more structure, leading to \emph{higher triangulations}. This use of the word `higher' is meant to indicate that one asks for higher octahedron axioms, i.e., that one also encodes iterated (co)fibers associated to longer strings of composable morphisms (see \cite{maltsiniotis:higher} for a precise definition). It can be shown that homotopy categories of stable $\infty$-categories or stable model categories can be naturally endowed with higher triangulations \cite[\S13]{gst:dynkin-A}, illustrating the slogan that `these enhancements encode all the triangulated structure'.
\end{rmk}

\subsection{Stabilization and the universality of spectra}
\label{subsec:spectra}

In this section we discuss the stabilization process which can be realized by passing to $\infty$-categories of internal spectrum objects. Similar constructions were also carried out in the language of model categories (for example by Schwede~\cite{schwede:spectra-model} and Hovey~\cite{hovey:spectra}) as well as in the framework of derivators by Heller~\cite{heller:stable}.

\begin{defn}
Let \cC be a finitely complete, finitely cocomplete, and pointed $\infty$-category. A \textbf{prespectrum object} in \cC is a functor
\[
X\colon N(\lZ\times \lZ)\to \cC
\]
such that for all $i\neq j$ the value $X(i,j)$ is a zero object. The full subcategory of $\nFun(N(\lZ\times\lZ),\cC)$ spanned by the prespectrum objects is denoted by $\nPSp(\cC).$
\end{defn}

Here, we consider the poset $\lZ$ with the natural ordering as a category. Since only the diagonal entries are possibly non-trivial, we use the shorthand notation $X_m=X(m,m)$. Thus, a part of a prespectrum object $X\in\nPSp(\cC)$ looks like 
\[
\xymatrix@=1.5pc{
&0\ar[r]&X_{m+1}\\
0''\ar[r]&X_m\ar[r]\ar[u]&0'\ar[u]\\
X_{m-1}\ar[r]\ar[u]&0'''.\ar[u]&
}
\]
By definition of the suspension and loop functors $(\Sigma,\Omega)\colon\cC\rightleftarrows\cC$, given $X\in\nPSp(\cC)$ we obtain induced morphisms
\begin{equation}\label{eq:alpha-beta}
\alpha_{m-1}\colon\Sigma X_{m-1}\to X_m\qquad\text{and}\qquad \beta_m\colon X_m\to\Omega X_{m+1}.
\end{equation}

\begin{defn}
Let \cC be a finitely complete, finitely cocomplete, and pointed $\infty$-category and let $X\in\nPSp(\cC)$.
\begin{enumerate}
\item The prespectrum~$X$ is a \textbf{spectrum below n} if $\beta_m\colon X_m\stackrel{\sim}{\to}\Omega X_{m+1}$ is an equivalence for all $m< n$. The full subcategory of $\nPSp(\cC)$ spanned by the spectra below $n$ is denoted by $\nSp_n(\cC)$.
\item The prespectrum~$X$ is a \textbf{spectrum object} if $\beta_m\colon X_m\stackrel{\sim}{\to}\Omega X_{m+1}$ is an equivalence for all $m\in\lZ$. The full subcategory of $\nPSp(\cC)$ spanned by the spectrum objects is denoted by $\nSp(\cC).$
\end{enumerate}
\end{defn}

\begin{eg}
An important special case is obtained if we start with the pointed $\infty$-category $\cS_\ast$ of pointed spaces. In that case, we simplify notation and write $\nSp=\nSp(\cS_\ast)$ for the \textbf{$\infty$-category of spectra}.
\end{eg}

\begin{thm}
The $\infty$-category $\nSp$ of spectra is stable and presentable.
\end{thm}

We discuss further below that $\nSp$ is the \emph{stabilization} of the $\infty$-category~$\cS$ of spaces. More generally, given an $\infty$-category~\cC, we refer to 
\[
\nStab(\cC)=\nSp(\cC_\ast)
\]
as the \textbf{stabilization} of~\cC.

Under certain assumptions on a pointed $\infty$-category \cC we will now construct a \emph{spectrification functor}, i.e., a left adjoint $L\colon\nPSp(\cC)\to\nSp(\cC)$ to the fully faithful inclusion functor $\iota\colon\nSp(\cC)\to\nPSp(\cC)$, exhibiting $\nSp(\cC)$ as a localization of $\nPSp(\cC)$. To begin with there is the following result.

\begin{prop}\label{prop:n-spectrify}
Let \cC be a finitely complete, finitely cocomplete, and pointed $\infty$-category. The fully faithful inclusion $\iota_n\colon\nSp_n(\cC)\to\nPSp(\cC)$ admits a left adjoint
$L_n\colon\nPSp(\cC)\to\nSp_n(\cC)$,
\[
(L_n,\iota_n)\colon\nPSp(\cC)\rightleftarrows\nSp_n(\cC).
\]
\end{prop}

The idea is of course that spectra below a certain level are somehow determined by the higher levels. And in fact, the left adjoint $L_n$ can be constructed as follows. Given a prespectrum $X\in\nPSp(\cC)$ we restrict it to the full subcategory of $N(\lZ\times\lZ)$ spanned by
\[
Q_n=\{(i,j)\in\lZ\times\lZ\mid i\neq j\;\mbox{or}\; i=j\geq n\}
\]
and then set
\[
L_n(X)=\nRKan_{Q_n\hookrightarrow N(\lZ\times\lZ)}(X\!\mid_{Q_n}).
\]
Here $\nRKan$ stands for an $\infty$-categorical variant of the usual right Kan extension; see \cite[\S4.3]{HTT} or \autoref{per:derivator}. Under suitable completeness assumptions on the $\infty$-categories involved, right Kan extensions can again be calculated pointwise, i.e., are given by limits over certain slice categories. In our situation, the corresponding slice categories are cofinally finite and the above right Kan extensions hence exist.  The essential image of $L_n$ consists of the spectra below $n$. 

With a bit more care, one can show that there is a sequence of functors
\[
\id\to L_0\to L_1\to L_2\to\ldots\colon\nPSp(\cC)\to\nPSp(\cC),
\]
and it is hence tempting to simply set $L:=\colim_n L_n.$ This in fact works if one imposes the following conditions on the $\infty$-category~\cC. 

\begin{prop}\label{prop:spectrify}
Let \cC be a finitely complete, countably cocomplete, and pointed $\infty$-category. If the loop functor $\Omega_{\cC}\colon\cC\to\cC$ commutes with sequential colimits, then 
\[
L:=\ncolim _n L_n\colon\nPSp(\cC)\to\nPSp(\cC)
\]
is a localization with essential image $\nSp(\cC).$ We refer to $L$ as the \textbf{spectrification functor}.
\end{prop}

An important example of an $\infty$-category satisfying these assumptions is the $\infty$-category of pointed spaces. In this case, let $\cD_n\subseteq\nSp_n$ be the full subcategory spanned by those $X$ such that $\alpha_m\colon\Sigma X_m\to X_{m+1}$ defined in \eqref{eq:alpha-beta} is an equivalence for $m\geq n$. Thus, morally, such a prespectrum~$X$ is essentially determined by its value $X_n$. And in fact, the evaluation map  $\nev_n\colon\cD_n\to\cS_\ast$ is an acyclic Kan fibration. Let us choose a section $s_{\tilde{\Sigma}^{\infty-n}} \colon\cS_\ast\to\cD_n$ of $\nev_n$ and set 
\[
\tilde{\Sigma}^{\infty-n}\colon\cS_\ast\xrightarrow{s_{\tilde{\Sigma}^{\infty-n}}}\cD_n\to\nPSp.
\]
(We note again that this is well-defined since the space of sections is a contractible Kan complex; see the discussion around \eqref{eq:setions}.) Denoting the $n$-th evaluation functor $\nPSp\to\cS_\ast$ by $\tilde{\Omega}^{\infty-n}$ we obtain an adjunction
\[
(\tilde{\Sigma}^{\infty-n},\tilde{\Omega}^{\infty-n})\colon\cS_\ast\rightleftarrows\nPSp.
 \]
Combining this with \autoref{prop:n-spectrify} and \autoref{prop:spectrify} we deduce the following result.

\begin{prop}
There is the following sequence of adjunctions
\[
\xymatrix{
(\Sigma_+^{\infty-n},\Omega_-^{\infty-n})\colon\cS\ar@<0.5ex>[r]^-{{}_+}&\cS_\ast\ar@<0.5ex>[r]^-{\tilde{\Sigma}^{\infty-n}}\ar@<0.5ex>[l]^-{{}_-}&\nPSp\ar@<0.5ex>[r]^-{L}\ar@<0.5ex>[l]^-{\tilde{\Omega}^{\infty-n}}&\nSp\ar@<0.5ex>[l]^-{\iota}.
}
\]
\end{prop} 

Dropping the first adjunction in the proposition we obtain the adjunction
\[
(\Sigma^{\infty-n},\Omega^{\infty-n})\colon\cS_\ast\rightleftarrows\nSp,
\]
and it turns out that a similar adjunction exists for arbitrary pointed, presentable $\infty$-categories. In fact, the evaluation functor $\Omega^{\infty-n}\colon\nSp(\cC)\to\cC$ clearly makes sense for every finitely complete, finitely cocomplete, and pointed $\infty$-category~\cC. If \cC is moreover presentable, then it can be shown that $\Omega^{\infty-n}$ satisfies the assumptions of the special adjoint functor theorem (\autoref{thm:SAFT}). Thus, there is a left adjoint $\Sigma^{\infty-n}\colon\cC\to\nSp(\cC)$, the \textbf{suspension spectrum functor} or the \textbf{$n$-th free spectrum functor},
\[
(\Sigma^{\infty-n},\Omega^{\infty-n})\colon\cC\rightleftarrows\nSp(\cC),
\] 
although in this generality the functor $\Sigma^{\infty-n}$ does not admit such a nice explicit description as in the case of $\cS_\ast$. In the following important result \cite{DAGI} we denote by $\FunL(\cC,\cD)$ the $\infty$-category of colimit-preserving functors between presentable $\infty$-categories $\cC$ and $\cD$.

\begin{thm}\label{thm:stabilization}
Let $\cC,\cD$ be pointed, presentable $\infty$-categories and let $\cD$ be stable. Restriction along the suspension spectrum functor $\Sigma^{\infty}\colon\cC\to\nSp(\cC)$ induces an equivalence of $\infty$-categories
\[
\FunL(\nSp(\cC),\cD)\stackrel{\sim}{\to}\FunL(\cC,\cD).
\]
\end{thm}

As a central special case, let us again consider the $\infty$-category $\cS_\ast$. We refer to the image of the zero sphere $S^0=\Delta^0_+$ under $\Sigma^{\infty}\colon\cS_\ast\to\nSp$ as the \textbf{sphere spectrum}. 

\begin{cor}\label{cor:spectra}
Let $\cD$ be a stable, presentable $\infty$-category. Evaluation at the sphere spectrum induces an equivalence of $\infty$-categories
\[
\FunL(\nSp,\cD)\stackrel{\sim}{\to}\cD.
\]
\end{cor}

In fact, this follows from \autoref{thm:stabilization} by observing that the evaluation map factors as 
\[
\FunL(\nSp,\cD)\stackrel{\sim}{\to}\FunL(\cS_\ast,\cD)\stackrel{\sim}{\to}\cD
\]
where the second equivalence is given by evaluation on $S^0$; see \autoref{prop:pointed-spaces}. This corollary makes precise that the $\infty$-category $\nSp$ of spectra is the \emph{free stable} $\infty$\emph{-category on one generator}, namely on the sphere spectrum.

\subsection{Tensor products of presentable $\infty$-categories}
\label{subsec:tensor}

We now turn to the tensor product of presentable $\infty$-categories which plays a key role in the construction of the smash product on the $\infty$-category of spectra; see \S\ref{subsec:smash}. For the convenience of the reader we begin by some heuristics indicating the corresponding construction in ordinary category theory. 

Let $\cC_1$ and $\cC_2$ be locally presentable categories. The \textbf{(cocontinuous) tensor product} of $\cC_1$ and $\cC_2$ is a locally presentable category $\cC_1\otimes\cC_2$ together with a universal bilinear map, i.e., a functor 
\begin{equation}\label{eq:tensor}
\cC_1\times\cC_2\to\cC_1\otimes\cC_2
\end{equation}
which preserves colimits in both variables separately and which is universal in the following sense: For any locally presentable category $\cD$ restricting along \eqref{eq:tensor} induces an equivalence of categories
\[
\FunL(\cC_1\otimes\cC_2,\cD)\stackrel{\sim}{\to}\FunLL(\cC_1\times\cC_2,\cD)
\]
where $\FunLL(\cC_1\times\cC_2,\cD)\subseteq\nFun(\cC_1\times\cC_2,\cD)$ is the full subcategory spanned by the functors which preserve colimits in both variables separately.

It can be shown that the cocontinuous tensor product always exists and that it admits the explicit description
\begin{equation}\label{eq:model-tensor}
\cC_1\otimes\cC_2=\FunR(\cC_1\op,\cC_2),
\end{equation}
where $\FunR(-,-)\subseteq\nFun(-,-)$ denotes the full subcategory spanned by the \emph{limit-preserving} functors. 

As a first step one has to show that the category $\FunR(\cC_1\op,\cC_2)$ is again locally presentable. Here we content ourselves by verifying this in the case that $\cC_1$ is a category of presheaves, i.e., that $\cC_1\simeq\nFun(A\op,\cSet)$ for a small category~$A$. Using that presheaf categories are cocompletions (\autoref{thm:yoneda}) we obtain equivalences
\begin{equation}\label{eq:tensor-locpres}
\begin{split}
\nFun(A\op,\cSet)\otimes\cC_2&=\FunR(\nFun(A\op,\cSet)\op,\cC_2)\\
&\simeq\FunL(\nFun(A\op,\cSet),\cC_2\op)\op\\
&\simeq\nFun(A\op,\cC_2\op)\op\\
&\simeq\nFun(A,\cC_2).
\end{split}
\end{equation}
Since $\nFun(A,\cC_2)$ is again locally presentable, we deduce that the tensor product is locally presentable in this case.

A lengthy formal calculation using the duality between left adjoint and right adjoint functors as well as \autoref{thm:catSAFT} implies that for locally presentable categories $\cC_1,\cC_2,$ and $\cD$ there is an equivalence of categories
\[
\FunL(\FunR(\cC_1\op,\cC_2),\cD)\simeq\FunLL(\cC_1\times\cC_2,\cD),
\]
showing that \eqref{eq:model-tensor} provides a model for the cocontinuous tensor product. 

The symmetry isomorphism $\cC_1\times\cC_2\cong\cC_2\times\cC_1$ together with the Yoneda lemma implies that the tensor product is symmetric. In a similar way one easily observes that the tensor product is suitably associative. If we specialize the calculation \eqref{eq:tensor-locpres} to the case of the terminal category $A=\bbone$, then we obtain an equivalence $\cSet\otimes\cC_2\simeq\cC_2$, showing that the locally presentable category $\cSet$ of sets behaves as a monoidal unit for the cocontinuous tensor product.

Finally, for locally presentable categories $\cC_1,\cC_2$, and $\cD$ we observe that there are equivalences of categories
\begin{equation}\label{eq:tensor-closed}
\begin{split}
\FunL(\cC_1\otimes\cC_2,\cD)&\simeq\FunLL(\cC_1\times\cC_2,\cD)\\
&\simeq\FunL(\cC_1,\FunL(\cC_2,\cD)),
\end{split}
\end{equation}
showing that the cocontinuous tensor product is closed with internal hom given by categories of colimit-preserving functors.

Thus, as an upshot, these heuristics suggest that the cocontinuous tensor product yields some kind of a closed symmetric monoidal structure. Slightly more precisely, there is a symmetric, closed $2$-multicategory such that
\begin{enumerate}
\item objects are locally presentable categories,
\item morphisms of arity~$n$ are functors $\cC_1\times\ldots\times\cC_n\to\cD$ which preserve colimits in each variable separately, 
\item and $2$-cells of arity~$n$ are natural transformations between such functors,
\end{enumerate}
and the above discussion suggests that this $2$-multicategory is \emph{representable}, i.e., that it comes from a symmetric monoidal closed bicategory.

In \cite{HA} Lurie establishes a variant of this monoidal structure for presentable $\infty$-categories. Let $\Prl$ be the (very large) \textbf{$\infty$-category of presentable $\infty$-categories} and colimit-preserving functors.

\begin{thm}\label{thm:tensor}
The $\infty$-category $\Prl$ admits a closed symmetric monoidal structure $\Prlo\to N(\Fin)$ such that the following properties are satisfied.
\begin{enumerate}
\item The tensor product $\cC_1\otimes\cC_2$ corepresents the functor which sends $\cD$ to the $\infty$-category $\FunLL(\cC_1\times\cC_2,\cD)$ of functors which preserve colimits separately in both variables.
\item The $\infty$-category $\cC_1\otimes\cC_2$ is equivalent to $\FunR(\cC_1\op,\cC_2)$.
\item The $\infty$-category $\cS$ of spaces is the monoidal unit.
\item The internal hom is given by $\FunL(\cC_1,\cC_2)$.
\end{enumerate}
\end{thm}  

We denote by $\Prlsmon$ the $\infty$-category of presentable, symmetric monoidal closed $\infty$-categories and symmetric monoidal, colimit-preserving functors. The following is a variant of \autoref{per:monoids}.

\begin{prop}
There is an equivalence of $\infty$-categories 
\[
\AlgE(\Prlo)\simeq\Prlsmon.
\]
\end{prop}

Thus, \autoref{prop:unit-algebra} applied to the monoidal structure of \autoref{thm:tensor} implies that the $\infty$-category~$\cS$ of spaces endowed with a certain symmetric monoidal structure is an initial object in $\Prlsmon$. It turns out that the monoidal structure is the usual Cartesian monoidal structure and that it can be characterized by the properties that
\begin{enumerate}
\item the pairing $\times\colon\cS\times\cS\to\cS$ preserves colimits separately in both variables and
\item the point $\Delta^0\in\cS$ is a monoidal unit.
\end{enumerate}

\subsection{The smash product}
\label{subsec:smash}

In this subsection we briefly discuss an $\infty$-categorical version of the smash product monoidal structure on spectra. The stabilization of presentable $\infty$-categories can be summarized by the following theorem. We denote by $\Prlst\subseteq\Prlpt\subseteq\Prl$ the full subcategories spanned by stable, presentable and pointed, presentable $\infty$-categories, respectively. 

\begin{thm}
The stabilization $\nStab\colon\Prl\to\Prlst$ of presentable $\infty$-categories factors as a composition of adjunctions 
\[
\Prl\rightleftarrows\Prlpt\rightleftarrows\Prlst.
\]
\end{thm}

It turns out that the symmetric monoidal structure $\Prlo\to N(\Fin)$ has similar variants in the pointed and in the stable context. Here we collect the variant of \autoref{thm:tensor} for stable, presentable $\infty$-categories.

\begin{thm}\label{thm:tensor-st}
The $\infty$-category $\Prlst$ admits a closed symmetric monoidal structure $\Prlsto\to N(\Fin)$ such that the following properties are satisfied.
\begin{enumerate}
\item The tensor product $\cC_1\otimes\cC_2$ corepresents the functor which sends $\cD$ to the $\infty$-category $\FunLL(\cC_1\times\cC_2,\cD)$ of functors which preserve colimits separately in both variables.
\item The $\infty$-category $\cC_1\otimes\cC_2$ is equivalent to $\FunR(\cC_1\op,\cC_2)$.
\item The $\infty$-category $\nSp$ of spectra is the monoidal unit.
\item The internal hom is given by $\FunL(\cC_1,\cC_2)$.
\end{enumerate}
\end{thm}  

And there is a similar variant for $\Prlpt$ the only difference being that the $\infty$-category $\cS_\ast$ of pointed spaces is the monoidal unit.

Let $\Prlstsmon$ denote the $\infty$-category of stable, presentable, symmetric monoidal closed $\infty$-categories and symmetric monoidal, colimit-preserving functors. As a consequence of \autoref{per:monoids} there is the following result.

\begin{prop}
There is an equivalence of $\infty$-categories 
\[
\AlgE(\Prlsto)\simeq\Prlstsmon.
\]
\end{prop}

An application of \autoref{prop:unit-algebra} to $\Prlsto$ implies that the $\infty$-category~$\nSp$ of spectra can be endowed with a certain symmetric monoidal structure, the \textbf{smash product}, such that the resulting symmetric monoidal $\infty$-category $\nSp^\otimes$ is an initial object in $\Prlstsmon$. It turns out that the monoidal structure can be characterized by the properties that
\begin{enumerate}
\item the monoidal structure $\otimes\colon\nSp\times\nSp\to\nSp$ preserves colimits separately in both variables and
\item the sphere spectrum is the monoidal unit.
\end{enumerate}

Having the smash product at our disposal we can finally make the following definition.

\begin{defn}
\begin{enumerate}
\item The \textbf{$\infty$-category of $\Eoo$-ring spectra} is the $\infty$-category of commutative algebra objects $\AlgE(\nSp^\otimes)$. 
\item The \textbf{$\infty$-category of $\Aoo$-ring spectra} is the $\infty$-category of (associative) algebra objects $\AlgA(\nSp^\otimes)$.
\end{enumerate}
\end{defn}

These are $\infty$-categorical versions of the more classical model categories of commutative or associative ring spectra. In fact, for a precise statement along these lines using \emph{symmetric spectra} see \cite[Example~8.21]{DAGIII}. Having these key notions in place one could now begin with an $\infty$-categorical study of stable homotopy theory \cite{HA} which together with the theory of $\infty$-topoi \cite[\S\S6-7]{HTT} provides the foundations for Lurie's $\infty$-categorical approach to derived algebraic geometry. For this we refer the reader to the literature. 

We conclude this section by a short discussion of monoidal aspects of the stabilization of presentable $\infty$-categories.

\begin{thm}
Let $\cC^\otimes$ be a closed symmetric monoidal structure on a presentable $\infty$-category. The $\infty$-categories $\cC_\ast,\nSp(\cC)$ admit closed symmetric monoidal structures, which are uniquely determined by the requirement that the respective free functors from $\cC$ are symmetric monoidal. Moreover, also the remaining left adjoint in 
\[
\cC\to\cC_\ast\to\nSp(\cC)
\]
 is uniquely symmetric monoidal. 
\end{thm}

The monoidal structures in the previous theorem enjoy the following universal properties. Given closed symmetric monoidal, presentable $\infty$-categories $\cC,\cD$ we denote by $\FunLo(\cC,\cD)$ the $\infty$-category of symmetric monoidal, colimit-preserving functors from $\cC$ to $\cD$.

\begin{thm}
Let $\cC,\cD$ be closed symmetric monoidal, presentable $\infty$-categories.
\begin{enumerate}
\item  If $\cD$ is pointed, then the symmetric monoidal functor $\cC\to\cC_\ast$ induces an equivalence of $\infty$-categories
\[
\FunLo(\cC_\ast,\cD)\stackrel{\sim}{\to}\FunLo(\cC,\cD).
\]
\item  If $\cD$ is stable, then the symmetric monoidal functor $\cC\to\nSp(\cC)$ induces an equivalence of $\infty$-categories
\[
\FunLo(\nSp(\cC),\cD)\stackrel{\sim}{\to}\FunLo(\cC,\cD).
\]
\end{enumerate}
\end{thm}

One way to summarize some of these results is as follows. 

\begin{thm}\label{thm:stab-monoidal}
The stabilization $\Prl\to\Prlst$ of presentable $\infty$-categories admits a monoidal refinement which factors as a composition of adjunctions
\[
\Prlo\rightleftarrows\Prlpto\rightleftarrows\Prlsto.
\]
\end{thm}

We conclude this course by the following perspective on a refined picture of the stabilization process \cite{gepner-groth-nikolaus}.

\begin{per}\label{per:stabilization}
Let us recall that stable $\infty$-categories are obtained from pointed $\infty$-categories by imposing an additional exactness property, asking that pushouts and pullbacks agree. There are two more intermediate steps given by preadditive and additive $\infty$-categories, respectively, both of which are obtained by adding certain exactness properties to pointed $\infty$-categories. A \textbf{preadditive} $\infty$-category is a pointed $\infty$-category with finite biproducts. It follows from these axioms that every object can be canonically turned into an $\Eoo$-monoid object and this monoid structure is given by the fold map. We say that a preadditive $\infty$-category is \textbf{additive} if these canonical $\Eoo$-monoid structures actually are $\Eoo$-group structures.

Focusing again on the context of presentable $\infty$-categories it can be shown that there are universal examples of such $\infty$-categories. In fact, the $\infty$-category $\MonE(\cS)$ of $\Eoo$-spaces is a preadditive, presentable $\infty$-category and the left adjoint $\cS\to\MonE(\cS)$ to the forgetful functor $\MonE(\cS)\to\cS$ exhibits $\MonE(\cS)$ as the free preadditive $\infty$-category on one generator, namely the free $\Eoo$-space generated by $\Delta^0$. More specifically, for every preadditive, presentable $\infty$-category~\cD, there are equivalences of $\infty$-categories 
\[
\FunL(\MonE(\cS),\cD)\stackrel{\sim}{\to}\FunL(\cS,\cD)\stackrel{\sim}{\to}\cD.
\]
A similar universal property is enjoyed by $\cC\to\MonE(\cC)$ where $\cC$ is a general presentable $\infty$-categories. And if we pass to the context of additive presentable $\infty$-categories instead, then the universal example is given by the $\infty$-category $\GrpE(\cS)$ of grouplike $\Eoo$-spaces.

It turns out that the stabilization of presentable $\infty$-categories factors through the $\infty$-category  of preadditive, presentable $\infty$-categories and also through the $\infty$-category  of additive, presentable $\infty$-categories. More precisely, the obvious forgetful functors admit left adjoints and the stabilization hence factors as
\[
\Prl\rightleftarrows\Prlpt\rightleftarrows\Prlpre\rightleftarrows\Prladd\rightleftarrows\Prlst.
\]
Finally, let us also mention that this factorization admits a monoidal refinement parallel to the one in \autoref{thm:stab-monoidal}. For details and applications we refer the reader to \cite{gepner-groth-nikolaus}.
\end{per}

\bibliographystyle{alpha}
\bibliography{infinity}

\end{document}